%&amstex

\documentstyle{amsppt}
\magnification=1200
\vsize=9 true in
\hsize=6.5 true in
\voffset= -.35truein %%% 0pt
\hoffset=0pt
\overfullrule=0pt
%%%
\parskip=4pt
\mathsurround=2pt
\input amssym.def
\input amssym
\input psfig.sty

\def\IMSmarkvadjust{0 pt}
\def\IMSmarkhadjust{0 pt}
\def\IMSmarkhpadding{0 pt}
\def\SBIMSMark#1#2#3{
 \font\SBF=cmss10 at 10 true pt
 \font\SBI=cmssi10 at 10 true pt
 \setbox0=\hbox{\SBF \hbox to \IMSmarkhpadding{\relax}
                Stony Brook IMS Preprint \##1}
 \setbox2=\hbox to \wd0{\hfil \SBI #2}
 \setbox4=\hbox to \wd0{\hfil \SBI #3}
 \setbox6=\hbox to \wd0{\hss
             \vbox{\hsize=\wd0 \parskip=0pt \baselineskip=10 true pt
                   \copy0 \break%
                   \copy2 \break% 
                   \copy4 \break}}
 \dimen0=\ht6   \advance\dimen0 by \vsize \advance\dimen0 by 8 true pt
                \advance\dimen0 by -\pagetotal
	        \advance\dimen0 by \IMSmarkvadjust
 \dimen2=\hsize \advance\dimen2 by .25 true in
	        \advance\dimen2 by \IMSmarkhadjust

%
%   Check for publication info
%
%  \newread\jref
  \openin2=publishd.tex
  \ifeof2\setbox0=\hbox to 0pt{}
  \else 
     \setbox0=\hbox to 3.1 true in{
                \vbox to \ht6{\hsize=3 true in \parskip=0pt  \noindent  
                {\SBI Published in modified form:}\hfil\break
                \input publishd.tex 
                \vfill}}
  \fi
  \closein2
  \ht0=0pt \dp0=0pt
 \ht6=0pt \dp6=0pt
 \setbox8=\vbox to \dimen0{\vfill \hbox to \dimen2{\copy0 \hss \copy6}}
 \ht8=0pt \dp8=0pt \wd8=0pt
 \copy8
 \message{*** Stony Brook IMS Preprint #1, #2. #3 ***}
}

\SBIMSMark{1998/9}{September 1998}{revised June 1999}
\def\maybebreakhere{%
	\vskip 0pt plus .15\vsize  % add a bunch of glue
        \penalty-250               % encourage a break
        \vskip 0pt plus -.15\vsize % take away the glue
}
\def\boldsig{\sigma\!\!\!\!\sigma}
\def\wideQP{\smallskip\leftskip=.25in\rightskip=.25in\noindent}
\def\cl{\centerline}
\def\min{\roman{min}}
\def\max{\roman{max}}
\def\adm{\roman{Adm}}
\def\C{{\Bbb {C}}}
\def\N{{\Bbb {N}}}

\def\R{{\Bbb {R}}}

\def\A{{\frak {A}}}
\def\cf{{\frak {C}}}
\def\I{{\Cal {I}}}
\def\K{{\Cal {K}}}
\def\O{{\bold o}}
\def\p{{\bold p}}
\def\q{{\bold q}}
\def\v{{\bold v}}

\def\Neg{\roman{Neg}}
\def\>{>\!\!>}
\def\<{<\!\!<}
\def\={~=~}
\def\sgn{\roman{ sgn}}
\def\QED{  \rlap{$\sqcup$}$\sqcap$ \smallskip}
\def\QP{\smallskip\leftskip=.4in\rightskip=.4in\noindent}
\def\b1{{\pmb{1}}}
\def\ss{\smallskip}
\font\bit=cmssi12 at 12truept
\font\tenmsy=msym10

 at 10 truept
\textfont8=\tenmsy
\def\ssm{\smallsetminus}
\mathsurround = 1pt
\abovedisplayskip=4pt
\belowdisplayskip=4pt
\parskip=5pt
\def\ref{\hangindent=1pc \hangafter=1 \noindent}

%%% to let amstex generate the proper headline, without the usual title stuff
\author
John Milnor and Charles Tresser
\endauthor
\title
On Entropy and Monotonicity for Real Cubic Maps
\endtitle
%%% to supress pagenumber on first page, just dont use footlines.
%%% this won't work with running heads off, though.

\centerline{\bf On Entropy and Monotonicity for
Real Cubic Maps} %\quad[\S\S~1--4]}
\medskip

\centerline{John Milnor and Charles Tresser\footnote
{Partially supported by NSF grant DMS-97-04867.}
}\ss

\centerline{with an appendix by Adrien Douady and
Pierrette Sentenac}\ss
%\centerline{[draft of 6-16-98; last hand  = jm]}
\bigskip

\centerline{\bf Contents}

1. Introduction

2. Piecewise monotone maps and kneading theory

3. Parametrization of polynomial maps

4. Topological entropy and periodic orbits

5. The stunted sawtooth family

6. Contractibility of isentropes for the stunted sawtooth family

7. Bones in $P^2$

8. Monotonicity, intersections of bones, and the period $n$-skeleton

9. From connected bones to connected isentropes

%9. Intersections of bones, and the $n$-skeleton

%10. Bones, negative orbits, and the connected isentrope conjecture

Appendix A. (by A. Douady and P. Sentenac): Characterization of a
polynomial
\break
\indent\indent by its critical values

Appendix B: Tight symbol sequences and Thurston's theorem

Appendix C. Monotonicity versus antimonotonicity

References
\bigskip
\centerline{\bf \S 1. Introduction}
\medskip

Consider a continuous map $f$ from a closed interval $I$ to itself. In the
simplest cases, $I$ is the union (in a not necessarily unique way)
of sub-intervals where $f$ is monotone. The minimal number of
sub-intervals required (called the {\bit lap number\/})
is a rough measure of the complexity
of the map $f$. If now we think of $f$ as generating a dynamical
system, a quantitative description of the dynamic complexity
of $f$ is obtained by measuring the rate of exponential
growth of the lap numbers
of the successive iterates $f^{\circ 1}=f\,,\,f^{\circ 2}=f\circ f\,,\,
\dots$ of the map $f$. This growth rate is clearly invariant under
a continuous change of coordinate, and turns out to be equal to a measure
of dynamic complexity used in a much more general context and
known as the {\bit topological entropy\/} $h(f)$.

Maps on the interval provide the simplest examples to illustrate
the problem of understanding how dynamic complexity
evolves under deformations of a dynamical system. However, even the
basic problem of comparing the topological entropies of two
smooth interval maps which are close to each other is
only partly understood. The
best results so far have been obtained in the more specific case
of polynomial maps, where one can use complex analysis
methods in the study of the associated complex polynomial maps.
%The simplest example is
The theory is most complete in the particular case of
the family of quadratic maps $Q_v(x)=4vx(1-x)$, parametrized by the
critical value $0\le v \le 1$. It has been known for some time
that the topological entropy $h(Q_v)$
is a (non-strictly) increasing continuous function of the parameter $v$,
with $h(Q_0)=0$  and $h(Q_1)=\log 2$.
All known proofs of the monotonicity of topological
entropy for quadratic maps use complex analytic methods. (Compare [D],
[DH2],
[MvS2], [MTh]\footnote{%
  In contrast with these older proofs, a recent proof by M.~Tsujii does not
  depend on holomorphic dynamics, but does depend on complex analysis
  [Ts]. Notice that Tsuji indicates that his proof can be understood as
a
  local version of the argument in [MTh].}%
.) The particular proof which we
will generalize is based on the following ideas.
These maps $Q_v$ have a single critical point $c=\frac{1}{2}$
for $v\not =0$, and there
are a countable infinity of values of $v\in (0,1]$ such that $c$
belongs to a periodic orbit. The restriction of
$Q_v$ to any periodic orbit $x_1<\cdots<x_p$ can be described
combinatorially by the cyclic permutation $\O$ of $\{1,\ldots,p\}$
which we call its {\bit order type\/},
defined by the property that $Q_v(x_i)=x_{\O(i)}$. Thurston showed
that each order type for a periodic critical orbit which occurs in
this way for some $Q_v$, occurs for a single value of $v$ only.
(Compare \S8 and Appendix B.)
This fact is used in [MTh] to prove monotonicity. It has been recognized
for a long time that for each quadratic map with a periodic
critical orbit of period
$p$ there exists a countable infinity of quadratic maps with the same
topological entropy and with a critical orbit whose period is some multiple
of $p$. Thus the
monotonicity of the correspondence $v\mapsto h(Q_v)$ cannot be strict.

On the other hand, it has been proved more recently that any value
of the topological entropy in $[0,\log 2]$ which cannot
be achieved by a map having a periodic critical orbit, is realized
for a single value of $v$. ([Gs], [L2].)
This result is a special case of
the Generic Hyperbolicity Conjecture,
which states that every rational
map can be approximated arbitrarily closely by a rational
map (real if the original map is real, and a polynomial if the original
is a polynomial) such that
the orbits of all critical points converge to attracting periodic
orbits. (For an early version on this conjecture, see Fatou [F].)
%Similar definitions of the generic hyperbolicity property
%can be formulated for other families of maps, and one such formulation
%will play an important role below.

One would like to understand as much as possible of the global
bifurcation theory for polynomials of higher degree. But there,
even the relevant concepts are harder to isolate because families
of degree $d$ polynomials depend naturally
on $d-1$ parameters. (In fact we will usually work with the integer
$m=d-1$.) For each fixed
degree $d$ and for each fixed sign $\pm 1$ we consider the family
%${\Cal P}_\pm^{d-1}$
of all degree $d$ polynomial maps from the interval $I=[a,b]$
to itself,
% we mean the set of all polynomial maps from $I$ to itself of degree $d$,
with all $d-1$ critical points in $I$, which send
the boundary $\partial I=\{a,b\}$ into itself, and with
leading coefficient
of specified sign, compactified
for $d$ even by the constant map with the same boundary behavior.
(Compare Figures 5, 6.) For these families %${\Cal P}_\pm^{d-1}$,
it has been proposed to generalize the monotonicity
property of the topological entropy when $d=2$ to the connectedness
of the topological entropy level sets or {\bit isentropes\/}.
Numerical computations in the case when $d=3$
have suggested the Connected Isentrope Conjecture according to
which all isentropes are connected [M1]. (Compare the earlier report
[DGMT].
In fact, for any $d$,
one can ask the sharper question as to whether
isentropes are contractible, or cellular.)

In the degree 3 case, the ``{\bit bones}'', or sets of parameters such
that one critical point is periodic with a specified order type,
are no longer points, as when $d=2$, but
are smooth curves. We conjectured [DGMT] that these bones cannot have any
connected component which is a simple closed curve or ``{\bit
bone-loop\/}''.
(It follows that every bone is a connected simple arc, except in
a few exceptional cases with very low periods.) Furthermore, we
sketched a proof of the following implications:
\medskip

\centerline{\bit Generic Hyperbolicity for Real Cubic Maps}

\hskip 1.55in{\bit $\Longrightarrow$~ No Bone-Loops} %\hskip 1.4in{$(*)$}

\hskip 1.65in{\bit  $\Longrightarrow$~ Isentropes are Connected.}\smallskip

\noindent
The generic hyperbolicity question remains open. However, C. Heckman
[He] has been able to prove a weaker version which is enough to show
that there are no bone-loops. The present paper will assume Heckman's
result, and provide a more detailed
exposition for the last implication. Thus, assuming
that there are no bone loops, we will prove
that isentropes, in either of the two families of real cubic
maps, are indeed connected.

(In fact we
will derive the slightly more precise statement that isentropes are
cellular sets. However, even assuming generic hyperbolicity, we do
not know whether isentropes are contractible. Compare [FT] which shows,
for the family of analytic circle maps
$$\theta~\mapsto~\theta+ a+b\sin(\theta)\quad(\text{mod}~ 2\pi)~,$$
that the zero isentrope is
non locally connected, with a ``comb'' structure.)\smallskip

The paper is organized as follows. \S2 contains a short discussion
of kneading theory, generalized so as to allow maps with ``plateaus''
or intervals of constancy. \S 3 describes
the parametrization of families of polynomial maps of the interval by
their critical values.
An alternate approach, due to Douady and Sentenac, is given in
Appendix A. We discuss in \S 4 the main aspects of topological
entropy that we need, in particular continuity properties and
relations between topological entropy and
kneading information. In \S 5, we describe families of continuous
maps of the interval closely related to kneading theory, that we
call {\bit stunted sawtooth maps\/}: Corresponding to the family of degree
$d$ polynomial maps with leading coefficient of specified sign,
%To each family ${\Cal P}_\pm^{d-1}$ of polynomial maps
there is an essentially unique family
of stunted sawtooth maps which mimic their behavior. Our basic hope is that
most of the essential features
of a family of polynomial maps of the interval are more or less faithfully
mirrored in the corresponding stunted sawtooth family, where they are
much easier to verify. As an example, the
Generic Hyperbolicity Conjecture is
easily verified for stunted sawtooth families. In \S 6, for each family of
stunted sawtooth maps, we show that all isentropes are contractible,
and therefore connected. This might be a step toward the Connected
Isentrope Conjecture, or some more precise form of monotonicity property,
for polynomial families of any degree.
The rest of the paper is devoted to the proof that isentropes are connected
for the two families of cubic maps. In order to use
what we know for sawtooth maps
in the study of cubic maps, we recall
the definition of bones in \S 7, discuss their elementary
properties in \S 8, and finish the proof in \S9.

\maybebreakhere
\bigskip
\centerline{\bf \S 2. Piecewise Monotone Maps and Kneading Theory.}
\medskip

(Compare  [My], [MSS], [MTh].)
Let $I=[a,b]$ be a closed interval of real numbers. A map $f:I\to I$
will be called {\bit piecewise monotone\/}
if $I$ can be covered by finitely many
closed intervals on which $f$ is monotone (but not necessarily strictly
monotone\footnotemark\footnotetext
{We will use the term piecewise {\it strictly\/} monotone
in the special case where $f$ has no intervals of constancy. However,
it will be important to allow limiting cases where $f$ may
have intervals of constancy. Because of this there is
some choice of appropriate definitions for
kneading theory (compare [BMT]), but
the following will be convenient for our purposes.}).

Let $\boldsig~=~(\sigma_0\,,\,\ldots\,,\,\sigma_m)$ be an alternating
sequence of $m+1$ signs. That is, we assume that
$\sigma_j=(-1)^j\sigma_0$ with $\sigma_0=\pm 1$. By
an $m${\bit -modal map of shape\/} $\boldsig$
will be meant a piecewise monotone map $f:I\to I$ as
above, together with a sequence of points
$$  a~=~c_0~<~c_1~<~\cdots~<~c_{m}~<~c_{m+1}~=~b \eqno (1)$$
in the interval $I=[a,b]$ satisfying the following
condition.\footnotemark\footnotetext
{We will sometimes need to consider the more general
case where
$$  a~=~c_0~\leq ~c_1~\leq ~\cdots~\leq ~c_{m}~
\leq ~c_{m+1}~=~b~. \eqno (1')$$
(Compare \S3.) However, this is awkward for kneading theory,
and we avoid it in this section.} For each
$0\le j\le m$ the restriction of $f$ to the interval $c_j\le x\le
c_{j+1}$ should be either monotone increasing or monotone decreasing
according as $\sigma_j$ equals $+1$ or $-1$.
Thus the $c_j$ are uniquely determined by $f$ in the
special case that $f$ is piecewise strictly monotone; but
not in the general case when $c_j$ may belong to an
interval of constancy. One also uses terms such
as {\bit unimodal} for an $m$-modal map with $m=1$, {\bit bimodal}
when $m=2$, and so on. We will use the notation
$$   (f,\bold{c})\qquad\text{with}\qquad \bold{c}\=(c_1\,,\,\ldots\,,\,c_m)
 ~ \in~ I^m $$
for an $m$-modal map when it is important to specify the precise choice
of the $c_i$. We will call $c_1\,,\,\ldots\,,\,c_m$
the {\bit folding points\/} of $f$,
and their images $v_i=f(c_i)$ the {\bit folding values\/}.
In the case of a $C^1$-smooth
function, note that the first derivative necessarily vanishes at our
folding
points, but may vanish at other points also.
In the important case of $C^2$-smooth maps without critical
inflection points, the $c_i$ are precisely the critical points of $f$,
and their images $v_i$ are the critical values. However, we introduce a
new name for these points in the general case in order to avoid
confusion. In that special case, note that $\sigma_j$
can be identified with the sign of the first derivative $f'(x)$ for
$c_j<x<c_{j+1}$. (If the second derivative $f''(c_j)$ is non-zero,
then $\sigma_j$ can also be identified with
the sign of $f''(c_j)$ for $j=1\,,\,\ldots\,,\, m$.)
When it is necessary to emphasize that some choice is involved
(that is, when some of the  $c_i$  lie in intervals of constancy)
we may refer to the  $c_i$  more explicitly as the
{\bit designated folding points} of the  $m$-modal map  $(f ,\bold{c})$.

The $m$-tuple
$\v=(v_1\,,\,\ldots\,,\,v_m)\in I^m$ will be called the {\bit folding
value vector\/} for $(f,\bold{c})$. Note that
$$\eqalign{ v_j~\le~ v_{j+1}\qquad &\text{when}
\qquad \sigma_j=+1\cr
    v_j~\ge ~v_{j+1}\qquad &\text{when}\qquad
\sigma_j=-1\cr}\eqno (2) $$
for $0<j<m$. In fact, if we set $v_0=f(a)\,,\,v_{m+1}=f(b)$, then this same
inequality (2) will be true for $0\le j\le m$.
If strict equalities hold in (2) for $0\le j\le m$,
or equivalently if $m$ is minimal,
then $f$ will be called a {\bit strictly $m$-modal map\/}, and
$m+1$ will be called the {\bit lap number} of $f$, denoted by $\ell(f)$.
Thus $\ell(f)$ is the minimum number of intervals of monotonicity
needed to cover the interval $I$.\ss
%\eject

\midinsert
\centerline{\psfig{figure=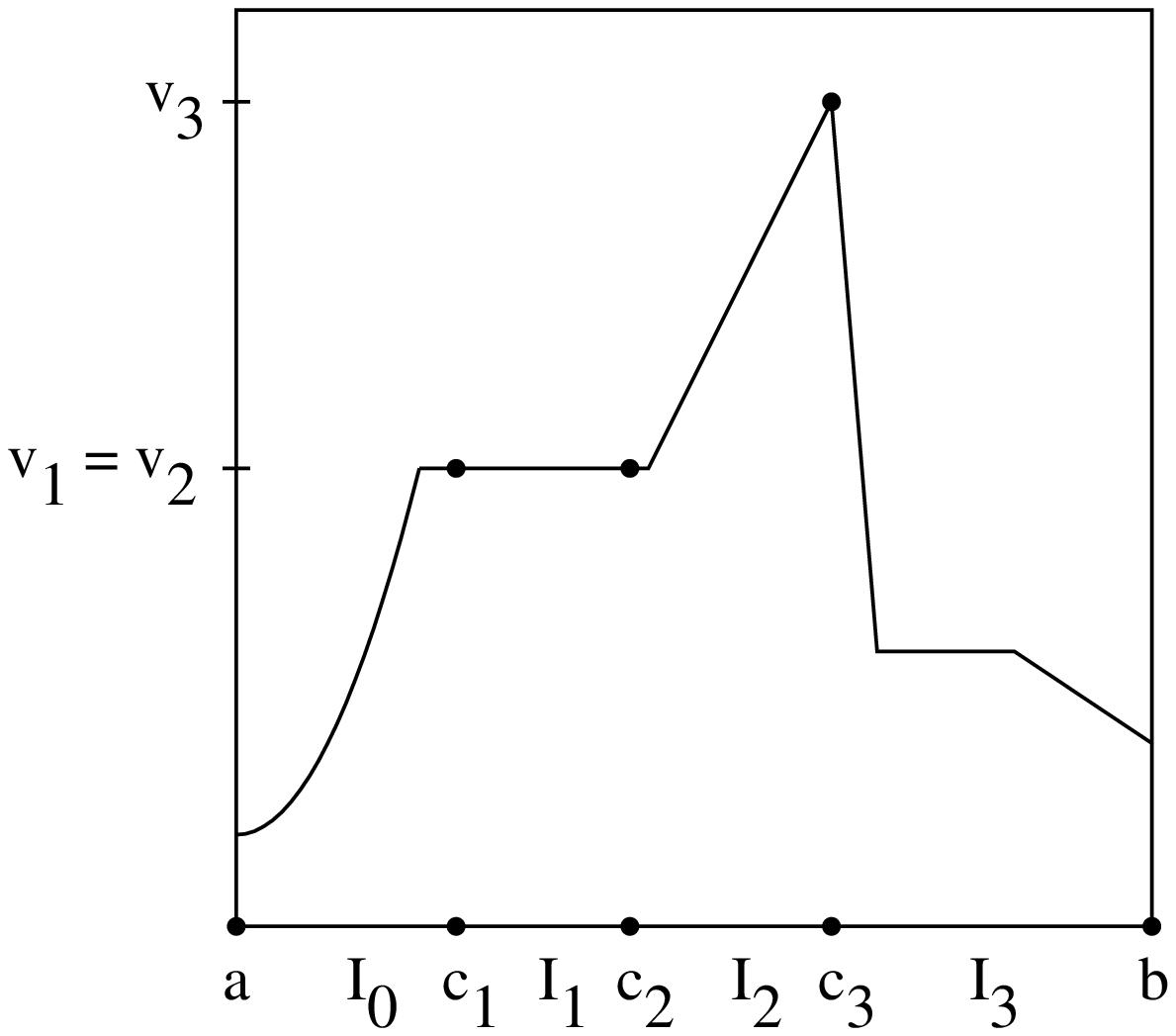,height=2.in}}\smallskip
{\QP\bit Figure 1. Graph of a $3$-modal map of shape $(+-+-)$, with
$v_1=v_2<v_3$. $($Depending on the choice of designated folding points
within the intervals of constancy,
this same map could equally well be considered as strictly unimodal, or as
$m$-modal for any odd value of $\,m\,.)$\medskip}
\endinsert

Let $f:I\to I$ be an $m$-modal map of shape $\boldsig$.
Partition the interval $I$ into disjoint subsets
$$   I~=~I_0~\cup~C_1~\cup~I_1~\cup~C_2~\cup~\cdots~\cup~C_m~\cup~ I_m $$
where $C_j$ is the set $\{c_j\}$ consisting
of a single folding point of $f$, and where the open sets
$$   I_0=[c_0\,,\,c_1)~,~I_1=(c_1\,,\,c_2)~,~\cdots~,~
     I_m=(c_m\,,\,c_{m+1}] $$
are the connected components of the complement $I\ssm(C_1\cup\cdots\cup
C_m)$.

We will also think of the $C_j$ and $I_j$ as
abstract symbols, forming the letters of an alphabet
%Let $    {\Cal A}={\Cal A}(n)$ be the ordered alphabet consisting
%of the $2n+1$ abstract symbols
$${\A}~=~{\A}(m)~=~\{I_0~,~C_1~
     ,~I_1~,~C_2~,~\cdots~,~C_m~,~ I_m\}~. $$
We order the $2m+1$ symbols in $\A(m)$
as they lie along the real line, so that
$$   I_0~<~C_1~<~I_1~<~C_2~<~\cdots~<~C_m~<~ I_m~. $$
Each point $x\in I$ has a unique {\bit address\/}
$A(x)\in{\A}$ defined by the condition that
$x\in A(x)$. Note that $A(x)\le A(y)$ whenever $x< y$.

Let $\A^\N$ be the set of all infinite sequences
$(A_0\,,\,A_1\,,\,\ldots)$ of symbols $A_i\in{\A}$. Each point
$x\in I$ has a well defined {\bit itinerary\/}
$$      {\Cal I}(x)~=~\big(A(x)\,,\,A(f(x))\,,\,
A(f^{\circ 2}(x))\,,\,\ldots  \big)~\in~\A^\N~.$$
Note that the map $f$ on the interval $I$
corresponds to the shift map
$$  \roman{shift\,}(A_0\,,\,A_1\,,\,\ldots)~
=~(A_1\,,\,A_2\,,\,\ldots) $$
on the space $\A^\N$ of symbol sequences, in the sense that the following
diagram is commutative
$$  \matrix
I&{\buildrel \I\over\longrightarrow}&\A^\N~~\\
        ~~\downarrow f&&~~~\downarrow\roman{shift}\\
        I&{\buildrel \I\over\longrightarrow}&\A^\N~~.
\endmatrix\eqno (3) $$

Some aspects of the structure of the orbit of $x$ reflect in
an obvious way in the itinerary.
For instance, if $x$ belongs to a periodic orbit, its itinerary must
be periodic, and since the only periodic orbits of monotone maps
have period 1 or 2, the period of the itinerary of $x$ is either
the period of the orbit or half of it. This periodic
itinerary (or any shift of it) will
be called the {\bit kneading type} of the periodic orbit.

For each fixed shape $\boldsig$, we introduce a partial ordering
of $\A^\N$ as follows. Define the sign function
$\epsilon:{\A}\to \{-1\,,\,0\,,\,1\}$
by the formula $\epsilon(I_j)=\sigma_j$ and $\epsilon(C_j)=0$.\break
(In the special case of a $C^1$-map whose derivative vanishes only at the
folding points, note that $\epsilon(A(x))$ can be identified with the
sign of the derivative $f'(x)$.)
% (Intuitively, this is the sign of the derivative at the points
%of $I_j$ or at $c_j$.)
By definition, $(A_0\,,\,A_1\,,\,\ldots)~<~(B_0\,,\,B_1\,,\,\ldots)$~~
if and only if there is a index $k\ge 0$ so that\break
$A_i=B_i$ for $i<k$, but
$$\eqalign{A_k<B_k \quad &\text{if the product}\quad
\epsilon(A_0)\,\epsilon
 (A_1)\cdots\epsilon(A_{k-1})\quad\text{is equal to}~~ +1\cr
A_k>B_k\quad & \text{if this product equals}~~ -1~. }$$
% $A_j=B_j$ for $j<m$
%$A_0=B_0\,,\,A_1=B_1\,,\,\ldots\,,\,A_{m-1}=B_{m-1}$,
%but either $A_m<B_m$ or $A_m>B_m$ according as the product
%$\epsilon(A_0)\,\epsilon(A_1)\cdots\epsilon(A_{m-1})$ is $+1$
%or $-1$.
(If the product is zero, or in other words if some $A_i=B_i$
with $i<k$ is equal to one of the folding point sets
$C_1\,,\,\ldots\,,\,C_n$, then the ordering is not defined. However, when
comparing two different itineraries for the same map,
this case occurs only when the two symbol sequences are
identically equal.) It is not difficult to check that
$$  \I(x)\le\I(y)\qquad\text{whenever}\qquad x\le y~.$$
Note that there is a unique smallest element $\I_\min$ and a unique
largest element $\I_\max$ in the space $\A^\N$ with this
ordering, so that
$\I_\min\le\I\le\I_\max$ for all $\I\in\A^\N$. For example:
$$   \I_\min\=\cases (I_0\,,\,I_0\,,\,I_0\,,\,I_0\,,\,\ldots)
&\text{if\quad
      $\boldsig=(+\cdots\pm)$}\\
     (I_0\,,\,I_m\,,\,I_m\,,\,I_m\,,\,\ldots) & \text{if\quad $\boldsig
     =(-\cdots +)$}\\
     (I_0\,,\,I_m\,,\,I_0\,,\,I_m\,,\,\ldots) & \text{if\quad $\boldsig=
     (-\cdots -)$.} \endcases $$

Sometimes we will need to truncate the infinite sequences
in $\A^\N$ and
consider the space $\A^k$ of finite sequences of some fixed length
$k$. There is a completely analogous ordering of such finite sequences.

The itineraries
$$   {\Cal K}_j~=~{\Cal I}(f(c_j))~\in~ \A^\N~. $$
of the folding values $f(c_j)$ will play a special role, and will
be called the {\bit kneading sequences\/} of $f$. To each
$m$-modal map $f$ there is associated a vector
$$   \bold{K}(f)~=~({\Cal K}_1\,,\,\ldots\,,\,{\Cal K}_m)
     ~\in~(\A^\N)^m
      $$
of kneading sequences. This vector of kneading sequences, together with
the shape $\boldsig$, will be called the {\bit kneading data\/}
%or the {\bit kneading invariant\/}
for the piecewise monotone map $f$.

The choice of kneading data imposes sharp restrictions on which
itineraries can actually occur. Suppose that
$\I(x)=(A_0\,,\,A_1\,,\,\ldots)$. Evidently:

{\bit Compatibility Condition 1.}
If some symbol $A_k$ is equal to $C_j$, then the sequence
$(A_{k+1}\,,\,\discretionary{}{}{}A_{k+2}\,,\,\ldots)$ which follows it
must
be equal to the kneading sequence $\K_j$.

In view of this condition, it is often convenient to terminate the
sequence\break $(A_0\,,\,A_1\,,\,\ldots)$
at the first $C_j$ which occurs in it, since the subsequent
symbols give no further information.

{\bit Compatibility Condition 2.} If $A_k=I_{j-1}$
or $A_k=I_j$, then the sequence which follows must satisfy either
$$   (A_{k+1}\,,\,A_{k+2}\,,\,\ldots)~\le~ \K_j $$
or
$$   (A_{k+1}\,,\,A_{k+2}\,,\,\ldots)~\ge ~\K_j $$
according as $\sigma_{j}=-1$ (so that $c_j$ is a local maximum point)
or $\sigma_{j}=+1$ (so that $c_j$ is a local minimum point).
\ss

{\bf Definition.} The symbol sequence $(A_0\,,\,A_1\,,\,\cdots)\in\A^\N$
will be called {\bit admissible\/} for the kneading data
$(\K_1\,,\,\ldots\,,\,
\K_m)$ if it satisfies these two compatibility conditions. (This
terminology
will be justified in 5.2.)
\smallskip

{\bf Remark.} Evidently the orbit of the folding point $c_j$ is periodic
if and only if the itinerary $\I(c_j)=(C_j\,,\,\K_j)$ is periodic. Using
Condition 1, a completely equivalent condition is that the kneading
sequence
$\K_j$ contains the symbol
$C_j$ as one of its entries.

{\bf Example 2.1.} The map $f(x)=(4x-1)^2(1-x)$ of shape $(-+-)$
on the unit interval has periodic kneading
sequences $\K_1=\overline{I_0I_2}$ and $\K_2=\overline{I_2I_0}$
(where the overline indicates a sequence which is to be repeated
infinitely often),
yet the critical points themselves are not periodic.
In fact $f$ maps both critical points $1/4$ and $3/4$ to the boundary
period two orbit $\{0,1\}$. On the other hand, for the map $f(x)=1-3x^2+2x^3$
of shape $(+-+)$ on the interval bounded by $(1\pm\sqrt 5)/2$,
the critical points $0\leftrightarrow
1$ are periodic, and hence their itineraries $\I(0)=\overline{C_1C_2}$
and $\I(1)=\overline{C_2C_1}$ are also periodic.\smallskip

In general, one cannot hope that the set
of all itineraries for a given map
will be completely determined by its kneading data. The first problem
is that the kneading data does not usually determine the itineraries of
the two endpoints $a,b\in\partial I$. Yet every itinerary $\I(x)$
which actually occurs must certainly satisfy
$$\I(a)~\le~\I(x)~\le~\I(b)~.$$
If these two boundary itineraries have been specified, then setting
$\K _0=\I(f(a))$ and $\K _{m+1}=\I(f(b))$ we can introduce
the following slightly sharper version of Condition 2 for an itinerary
$\I(x)=(A_0\,,\,A_1\,,\,\ldots)$:

\noindent {\bit Compatibility Condition 2$^\sharp$.}
If $A_k=I_j$, then
$$\K_j ~\le~ (A_{k+1}\,,\,A_{k+2}\,,\, \ldots ) ~\le~ \K_{j+1}\qquad
\text{whenever}\qquad\sigma_j=+1,$$
$$\K_j ~\ge~ (A_{k+1}\,,\,A_{k+2}\,,\, \ldots   ) ~\ge~ \K_{j+1}\qquad
\text{whenever}\qquad \sigma_j=-1~.$$

However, it is still not true that
every sequence satisfying Conditions 1 and $2^\sharp$ necessarily
occurs as itinerary for the given $m$-modal map. The following helps
to indicate the more serious difficulties.

{\bf Example 2.2.} Consider the following three families of unimodal maps
of shape $(+-)$ on the unit interval, with folding point $1/2$,
all parametrized by the folding value $v$,

\hskip .5in the tent family: $\quad\qquad T_v(x)\=2v\,\roman{Min}(x\,,\,1-x)$,

\hskip .5in
the quadratic family: $\quad Q_v(x)\=4vx(1-x)$, %\qquad and

\hskip .5in
and the family: $\qquad\qquad S_v(x)\=\roman{Min}(2x\,,\, v\,,\, 2-2x)$
\qquad (compare \S5).

\noindent In each case there is a unique parameter value $v$ so that
the orbit of the folding point has period 3, necessarily with
kneading sequence $\K_1=\overline{I_1 I_0 C_1}$.
(The corresponding parameter values are respectively
$$	v\=(1+\sqrt 5)/4\=0.80901\cdots~,\quad v'\=0.95796\cdots\,,
 \quad\text{and}\quad 	v''\=7/8\=0.875 $$
in the three cases.)
However, the itineraries which can actually occur for these three      
maps are different. For the tent case, the orbit of the folding point
is the only period 3 orbit. (The graph of $T_v^{\circ 3}$ touches the
diagonal without crossing it at the three points of this orbit.)
In the other two cases the corresponding graph definitely crosses the
diagonal and must cross back,
so there is a period 3 point close to the folding point
with itinerary $\overline{I_1I_1I_0}$. In the third case, the corresponding
graph is horizontal near the folding point, so that there
are also infinitely many nearby points with an itinerary $\I(x)$ equal to
$$  I_0\K_1=I_0\overline{I_1I_0C_1}\qquad\roman{or}\qquad I_1\K_1=
  I_1\overline{I_1I_0C_1}~. $$
Evidently such itineraries can occur only if the map
is constant on the entire interval between $x$ and the folding point.
Hence they can never occur for maps like $T_v$ or $Q_{v'}$
which are piecewise {\it strictly\/} monotone.\ss

Thus the compatibility conditions are not sufficient to guarantee the
existence of itineraries, or even of finite truncations of itineraries
which contain folding point symbols. However, the next lemma
provides a useful existence statement by
working only with finite sequences containing no folding point symbols.
Note that the analogues of Conditions 1, 2, and $2^\sharp$ for
finite sequences in $\A^k$ make perfect sense. In particular, the concept
of ``admissibility'' for finite sequences also
makes sense. A finite or infinite
symbol sequence $\{A_i\}$ will be called {\bit acritical\/}
if the $A_i$ all belong to the smaller alphabet
$$   \A_0\=\A_0(m)\=\{I_0\,,\,I_1\,,\,\ldots\,,\,I_m\}~\subset~\A(m)~, $$
with no folding point symbols.

{\QP{\bf Lemma 2.3.} \it Let $f$ be an $m$-modal map, and let
$(I_{\alpha_0}\,,\,I_{\alpha_1}\,,\,\ldots\,,\,I_{\alpha_k})$ be a finite
sequence of intervals $I_{\alpha_j}\in\A_0(m)$. There exists an orbit
$$x_0~\mapsto~ x_1~\mapsto~\cdots~\mapsto~
x_k~\mapsto~\cdots$$ with $x_i\in I_{\alpha_i}$ for all $i\le k$
if and only if this %finite acritical
sequence satisfies Condition $2^\sharp$, modified
so as to apply to sequences of finite length.\ss}

{\bf Proof.} If there exists such an orbit, then clearly the
sequence $(I_{\alpha_0}\,,\,\ldots\,,\,I_{\alpha_k})$ must satisfy the
suitably
modified Condition $2^\sharp$. Conversely, we will prove by induction
on $k$ that every sequence $(I_{\alpha_0},I_{\alpha_1},\ldots
,I_{\alpha_k})$ which satisfies this condition can
be realized by an orbit of $f$. This statement is
certainly true when $k=0$. Suppose then that
it is known for all sequences of shorter length. Suppose also,
to fix our ideas, that $\sigma_{\alpha_0}=+1$, so that
$$   \K_{\alpha_0}^{(k)}~\le~ (I_{\alpha_1}\,,\,I_{\alpha_2}\,,\,\ldots
\,,\,I_{\alpha_k})~\le~
     \K_{\alpha_0+1}^{(k)}  \eqno (4)$$
by Condition $2^\sharp$. (Here the superscript $(k)$ indicates
truncation to length
$k$.) By the induction hypothesis, there exists a point $x_1\in I$
with itinerary
$$   (I_{\alpha_1}\,,\,I_{\alpha_2}\,,\,\ldots
\,,\,I_{\alpha_k}\,,\,\ldots)~. $$
The proof will now be divided into two cases
according as strict inequalities do or do not hold in (4). If strict
inequalities hold, then it follows that $v_{\alpha_0}
<x_1<v_{\alpha_0+1}$.
It then follows from the intermediate value theorem that
there exists at least
one point $x_0\in I_{\alpha_0}$ with
$f(x_0)=x_1$, so that $x_0$ has
the required itinerary. On the other hand, if for example
$$   \K_{\alpha_0}^{(k)}~=~ (I_{\alpha_1}\,,\,I_{\alpha_2}\,,\,\ldots
\,,\,I_{\alpha_k})~, $$
then the point $x_0=c_{\alpha _0}+\epsilon$ with $\epsilon>0$
sufficiently small will have the required
itinerary, since the $I_{\alpha_j}$ are open subsets
of $I$. The remaining cases are completely analogous.\QED

It is often practical to first consider maps
whose boundary behavior is specified in the simplest way.
Suppose that $f$ maps the interval $I=[a,b]$ into itself.
%It will be convenient to introduce the notation
%$$  \partial_-(I)=a\,,\qquad\partial_+(I)=b $$
%for the two endpoints of the interval $I=[a,b]$.
Again choose some fixed shape $\boldsig$.

{\bf Definition.} The $m$-modal map $f:I\to I$ will be called
{\bit boundary anchored\/} for the shape $\boldsig$ if $f$ maps
the boundary $\partial I=\{a,b\}$ into itself by the rule
$$   f(a)\=\cases a & \text{if}\quad \sigma_0=+\cr
     b &\text{if}\quad \sigma_0=-~,\endcases\qquad\qquad
     f(b)\=\cases b & \text{if}\quad \sigma_m=+\cr
     a &\text{if}\quad \sigma_m=-~.\endcases \eqno(5)$$
%$$\eqalign{   f(a)\=\cases a & \text{if}\quad \sigma_0=+\cr
%    b &\text{if}\quad \sigma_0=-~.\endcases\cr
%    f(b)\=\cases b & \text{if}\quad \sigma_m=+\cr
%    a &\text{if}\quad \sigma_m=-~.\endcases}\eqno(5) $$
In the sequel, we will nearly always work with boundary anchored maps, so
that the distinction between Conditions 2 and $2^\sharp$ will disappear.

%$$  f(a)\=\left\{\eqalign{a &\qquad\text{if}\qquad \sigma_0=+1~\cr
%    b &\qquad \text{if}\qquad \sigma_0=-1~, \cr}\right . $$
%with
%$$  f(b)\=\left\{\eqalign{a &\qquad\text{if}\qquad \sigma_n=+1~\cr
%    b &\qquad \text{if}\qquad \sigma_n=-1~. \cr}\right . $$

\midinsert
\centerline{\psfig{figure=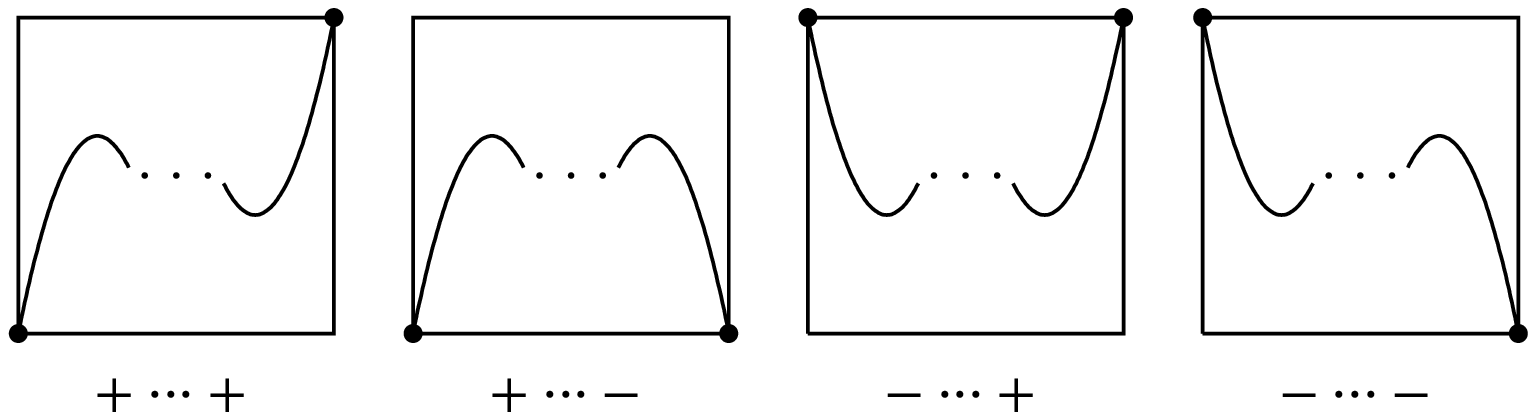,height=1.2in}}\medskip
{\QP\bit Figure 2. Boundary anchored maps: graphs illustrating the four
possible cases. These maps are $m$-modal with $m$ odd for the middle
pictures, and with $m$ even for the end pictures.\smallskip}
\endinsert

For a boundary anchored map, note that the itinerary $\I(a)$ is
precisely equal to $\I_\min$ and that $\I(b)=\I_\max$. It follows
that the sequences $\K _0=\I(f(a))$ and $\K _{m+1}=\I(f(b))$
are determined by the shape,
and that each one is equal to either $\I_\min$ or $\I_\max$.
Lemma 2.3 has the following consequence. (Compare 4.5 below.)

{\QP{\bf Corollary 2.4.} \it Let $f$ be a boundary anchored $m$-modal
map. Then a finite acritical sequence $(I_{\alpha_0}\,,\,\ldots\,,\,
I_{\alpha_k})$ is actually realized by some orbit\break
$x_0\mapsto\cdots\mapsto x_k$
if and only if it is admissible. Hence,
for any $k$, the set of all such sequences in $\A_0^k$
is determined by the kneading data for $f$.\ss}

The kneading sequences $\K_j$ are themselves itineraries, and hence
must satisfy these two compatibility
conditions. Evidently they must also satisfy the following.

{\bit Compatibility Condition 3 (for kneading data).}
$$  \eqalign{\K_j~\le~\K_{j+1}\qquad\text{if}\qquad \sigma_j=+1~,\cr
    \K_j~\ge~\K_{j+1}\qquad\text{if}\qquad
\sigma_j=-1~.\cr} $$
\ss

{\bf Definition.} The kneading data $\bold{K}=(\K_1\,,\,\ldots\,,\,\K_m)$
will be called {\bit admissible\/} for the specified shape $\boldsig$
if it satisfies these three compatibility conditions.\ss

As we shall see in 5.2, the three compatibility
conditions characterize the $m$-tuples of symbolic sequences
in $\A(m) ^\N$ that can actually occur as kneading data for some
$m$-modal map. However, if we restrict to polynomial (or analytic) maps,
then further restrictions are needed. (Compare [MaT1, p.179], as well as
2.2 and Appendix B, Example 4.)

%\eject
\bigskip
\maybebreakhere
\centerline{\bf \S 3. Parametrization of polynomials.}\medskip
\medskip

First consider a polynomial map $f:\R\to\R$ of degree $m+1$
with distinct real critical points %$f'(c_i)=0$, where
$$   c_1~<~c_2~<~\cdots~<~c_m~.$$
As in the previous section, we can form the
{\bit critical value vector\/} ($=$ folding value vector)
 $\v=(v_1\,,\,\ldots\,,\,v_m)
\in\R^m$, where $v_i=f(c_i)$.
Setting $\sigma_m$ equal to the sign of the leading coefficient, or
equivalently the
sign of the $(m+1)$-st derivative, and setting $\sigma_i=(-1)^{m-i}
 \sigma_m$, note the following strict form of (2):
$$   \sigma_i\,(v_{i+1}-v_i)~>~0\qquad\text{for}\qquad
 1\le i<m~. \eqno (2') $$
Again we will refer to $\boldsig=(\sigma_0\,,\,\ldots\,,\,\sigma_m)$ as
the shape. We will first prove the following.

{\QP{\bf Lemma 3.1.} \it
Given any $m$-tuple $(v_1\,,\,\ldots\,,\,v_m)\in\R^m$
satisfying the inequalities $(2')$, there exists a polynomial $f$
of degree $m+1$ with these critical values, listed in order of the
corresponding critical points as above. Furthermore $f$ is unique up to
precomposition with a positive affine transformation
$$   f(x)~\mapsto~f(px+q)\qquad\text{with}\qquad p>0~.$$\smallskip}

{\bf Proof of 3.1.} We will give a proof based on complex analysis.
Alternatively, a purely real
proof may be extracted from [MvS2, p.~120.], and a different real
proof, due to Douady and Sentenac, is provided in Appendix A.

To begin the construction, start with $m+1$ copies of the
Riemann sphere $\hat\C$, labeled as $\{k\}\times\hat\C$
for $0\le k\le m$. Slit the first two
copies $\{0\}\times\hat\C$ and $\{1\}\times\hat\C$ along the real
axis from $v_1$ to $\pm\infty$, taking the plus sign or the minus
sign according as the sign $\sigma_0$ is $+1$ or $-1$.
That is, if $\sigma_0 =+1$ so that $v_1$ is a local maximum,
we remove the
open interval consisting of real $z$ with $v_1<z<+\infty$, while
if $\sigma_0=-1$ so that $v_1$
is a local minimum we remove the open interval $-\infty<z<v_1$.
Similarly, for each $1\le k\le m$ we slit both $\{k-1\}\times
\hat\C$ and $\{k\}\times\hat\C$ from $v_k$
to $\pm\infty$, now choosing the sign according as $v_k$ is a local
maximum or minimum. The hypothesis $(2')$ guarantees that these
various slits can never meet.
\bigskip

\midinsert
\centerline{\psfig{figure=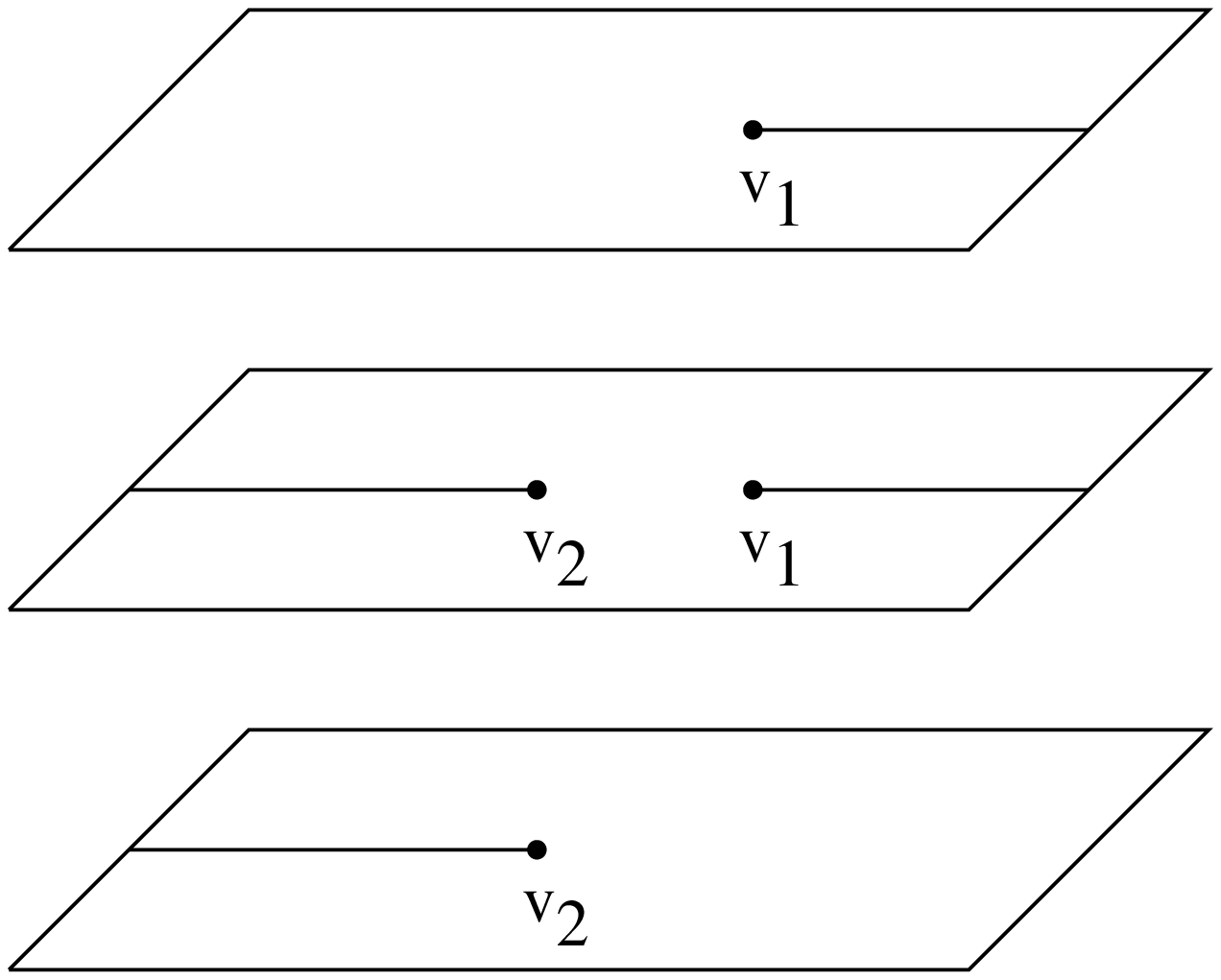,height=2in}}\medskip
{\QP\bit Figure 3. Three copies of the complex numbers, slit for the
construction of a bimodal map of shape $(+-+)$, with critical values
$v_1>v_2$.
\medskip}
\endinsert

Now, for each $1\le k\le m$, sew together the pair of corresponding
slits in $\{k-1\}\times\hat\C$ and $\{k\}\times\hat\C$ so
that the top edge of either one is matched with the
bottom edge of its mate. The result will be a compact simply connected
Riemann surface $S$. (Note that there is a canonical way
of assigning a conformal structure to this surface $S$, even at
the $m+1$ ramification points.) Furthermore, the natural projection maps
$(k,z)\mapsto z$ from $\{k\}\times\hat\C$ to $\hat \C$ fit together
to yield a holomorphic map $\eta:S\to\hat\C$ of degree $m+1$, with
the $v_k$
and the point at infinity as critical values. By the Uniformization
Theorem, $S$ must be conformally isomorphic to the
standard Riemann sphere, under some isomorphism
$u:\hat\C\to S$. (See for example [FK] or [Be].) If we choose
the isomorphism $u:\hat\C\to S$ in such a way that the points at
infinity correspond, then the composition $\eta\circ u:\hat\C\to
\hat\C$ is clearly a polynomial map with the required critical values.

To recover the real mapping, note that the complex conjugation
operations on the various copies $\{k\}\times\hat\C$ fit
together to yield an anti-holomorphic involution of $S$
whose fixed point set $F$ is the union
of the $m+1$ real axes, with the various slits removed. There is a
preferred ordering of $F\smallsetminus\{\infty\}$
so that its intersections with
the various $\{k\}\times\hat\C$ occur in the order of increasing $k$.
Now choose the conformal isomorphism $u:\hat\C\to S$ so that $\R\subset
\hat\C$ corresponds to $F\smallsetminus\{\infty\}$,
preserving orientation. In this way,
we obtain a real polynomial map $\eta\circ u|_\R$, with the specified
critical values occurring in the specified order. Evidently the conformal
isomorphism $u$, with these restrictions, is unique
up to an affine change of coordinates
$$   u(z)~\mapsto u(pz+q) $$
with $p,q\in\R$ and $p>0$.\QED
\smallskip

In practice, we are interested in polynomial maps which carry some closed
interval $I=[a,b]\subset\R$ into itself, and which are boundary
anchored, so that $\partial I=\{a,b\}$ maps into itself.
Let us fix the shape $\boldsig$.
It is sometimes convenient to set
$$   v_0\=v_0(\boldsig)\=f(a)~,\qquad\qquad
v_{m+1}\=v_{m+1}(\boldsig)\=f(b)~, $$
where now $f(a)$ and $f(b)$ are to be defined by formula (5).
%[IS THIS CLEAR ENOUGH?]
%$$  v_0\=\partial_{-\sigma_0}(I)\,,\qquad v_{n+1}\=\partial_{\sigma_n}(I)
%\,, $$
%so that every boundary anchored map $g$ satisfies
Thus a map $g$ is boundary anchored for the shape $\boldsig$
if and only if it satisfies
$g(a)= v_0$ and $g(b)= v_{m+1}$. Note that the inequality $(2')$
can be sharpened to include $v_0$ and $v_{m+1}$:
$$   \sigma_i\,(v_{i+1}-v_i)~>~0\qquad\text{for}\qquad
 0\le i\le m~. \eqno (2'') $$

{\QP{\bf Theorem 3.2.} \it Given an $m$-tuple
$(v_1\,,\,\ldots\,,\, v_m)
\in I^m$ satisfying the inequalities $(2'')$, %and satisfying one
% additional inequality $(6)$ in the unimodal case,
there exists one and only one boundary anchored polynomial map
$g:I\to I$ of degree $m+1$
which has distinct critical points
$$   a\,<\,c_1\,<\,\cdots\,<\,c_m\,<\,b~, $$
with $g(c_k)=v_k$.\smallskip}

{\bf Proof of 3.2.}
Let $f:\R\to\R$ be the polynomial whose existence is promised
by 3.1. Then we claim that there are unique real
numbers $p>0$ and $q$ so that the polynomial $g(x)=f(px+q)$
satisfies the additional conditions
that $g(\partial I)\subset \partial I$. It is then easy to check that
$g(I)\subset I$.

Without loss of generality, we may assume that $I$ is the unit
interval $[0,1]$. Suppose first that $\boldsig$ has the form $(+\cdots-)$.
Then $f(c_1)=v_1$ is a local maximum point.
Evidently $f$ maps the closed interval $[-\infty\,,\,c_1]$
homeomorphically onto the closed interval $[-\infty\,,\,v_1]$. Since
$v_1>a=0$, there is a unique point $-\infty<q<c_1$ with $f(q)=0$. Similarly,
$f$ maps the closed interval $[c_m\,,\,+\infty]$ homeomorphically onto
$[-\infty\,,\,v_m]$, reversing orientation. Since $v_m>a=0$,
there is a unique $p+q$ with
$c_m<p+q<+\infty$ and $f(p+q)=0$. The map $g(x)=f(px+q)$ will then have
the required properties.

The proof for the other three possible shapes is similar.
Details will be left to the reader.\QED

Thus far, we have assumed that our polynomials have distinct critical
points,
however it is often convenient to relax this condition
in order to obtain a compact parameter space.
First consider polynomials $f(x)$ with derivative
$$   f'(x)\= \kappa(x-c_1)\,\cdots(x-c_m)\qquad\text{where}
\qquad a\le c_1\le\cdots\le   c_m \le b$$
for some real constant $\kappa\ne 0$. The corresponding critical value
vector $(v_1,\ldots
,v_m)$ then satisfies the weaker inequalities (2):
%equivalent to (2)
$$   (v_{i+1}-v_i)\,\sigma_i~\ge~ 0\qquad\text{for}\qquad 1\le i<m~, $$
%\eqno (2) $$
where $\sigma_i=(-1)^{m-i}\sgn(\kappa)$. Just as in
Lemma 3.1, the polynomial
$f$ is uniquely determined, up to precomposition with a positive affine
transformation, by its critical value vector, which is required
to satisfy
only $(2)$. Either the proof above or the proof in Appendix A can be
adapted to show this. (Compare the Addendum to Appendix A.)
Details will be omitted.

There is a corresponding generalization of 3.2. However, when $m$
is odd, we must now allow for the case
of a constant function with $f(x)$ identically equal to $a$ or $b$.
In this limiting case, note that the critical
value vector $v_1=\cdots=v_m$ is still well defined,
although we can no longer distinguish $m$ critical points. With this
understanding, we have the following.

{\QP{\bf Theorem 3.3.} \it Given a specified shape $\boldsig$,
and given $(v_1\,,\,\ldots\,,\,
v_m)\in I^m$ satisfying $(2)$, there is one and only one
boundary anchored polynomial map $f:I\to I$ of degree $m+1$
having $(v_1\,,\,\ldots\,,\,v_m)$ as critical value vector.
\smallskip}

The proof is quite similar to the proof of 3.2.\QED

\midinsert
\centerline{\psfig{figure=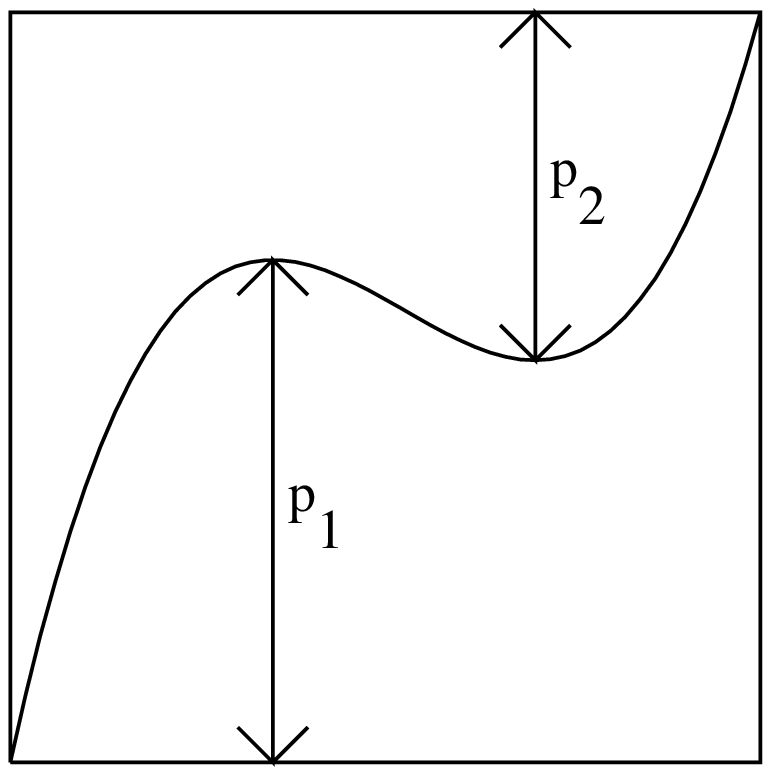,height=1.3in}\qquad\qquad\qquad
     \psfig{figure=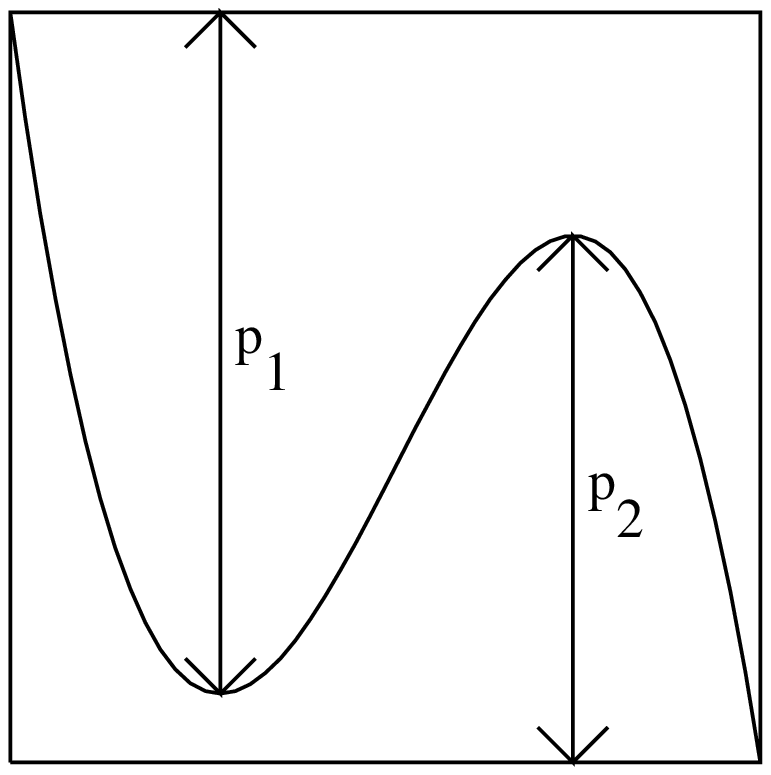,height=1.3in}}\ss
{\QP\bit Figure 4. The parameters $p_1$ and $p_2$ for a cubic map of
shape $+-+$ on the left, and $-+-$ on the right, on an interval of length
$1$.
\ss}
\endinsert

It would be quite natural to parametrize the set of all boundary anchored
polynomial maps $f:I\to I$of shape $\boldsig$ by the polyhedron consisting
of all vectors $\v=
(v_1\,,\,\ldots\,,\,v_m)\in I^m$ satisfying the inequalities
$(2)$ for the shape $\boldsig$. However, we will
rather work with a slightly different but affinely equivalent
polyhedron which is independent of the shape $\boldsig$,
and is invariant under affine reparametrization of the interval $I$ (as
used
in renormalization theory).
As a first attempt in this
direction, note that our interval $I$ is the union of non-overlapping
subintervals $I_j$, which map onto intervals $f(I_j)\subset I$. Let
$$
        \Delta_j = {\roman{length}(f(I_j))\over\roman{length}(I)} \in [0,1]
  $$
be the relative length of this image interval.
%  Intuitively, it seems that
%an increase in any one of these invariants $\Delta_j$ is likely to lead to
%more complex dynamic behavior.

Since our maps are boundary anchored, these invariants
$\Delta_0,\ldots,\Delta_m$ are not independent. (Their alternating sum is
constant, and various inequalities are needed to guarantee that
$f(I)\subset I$.) However,
we can obtain a set of invariants which are independent
and more manageable by setting
$$
        \Delta_j = p_j + p_{j+1} -1     \quad\roman{ for }\quad 0 < j < m ,
  $$
with
$$
         \Delta_0 = p_1 ,\quad           \Delta_m = p_m  .
  $$
(Note that the ``normalized total variation'' $\sum\Delta_i$ of our
$m$-modal
map is equal to a constant plus $2(p_1+\cdots+p_m)$.)
%$2\sum p_i$.)

If $I$ is the interval $[a,b]$, then we can write
$$
        p_i = \cases  (v_i-a)/(b-a)  &\roman{ if}\quad \sigma_i = +1 \cr
                             (b - v_i)/(b-a) &\roman{ if}\quad \sigma_1 =
-1~.
              \endcases\eqno (6)
  $$
(Compare Figures 4, 9.) In particular, in the special case of the unit
interval $I=[0,1]$ this formula simplifies to
$$
        p_i = \cases  v_i  &\roman{ if}\quad \sigma_i = +1 \cr
                             1 - v_i &\roman{ if}\quad \sigma_1 =
-1.\endcases
  $$
It is easy to check that the parameters
$p_1 , \ldots , p_m\in[0,1]$ satisfy only the relations
$$ %        p_i \in [0,1],\quad
                    p_i + p_{i+1} \ge 1~ .\eqno(7)
  $$
It may seem that
an increase in any $p_i$ should lead to more complex dynamical
behavior, but this is not quite true for the cubic family. (Compare
Figure 8.) However, in 5.6 we will see that the analogous statement
for the family of ``stunted sawtooth'' maps % of either shape is actually
is true.\ss

\midinsert
\centerline{\psfig{figure=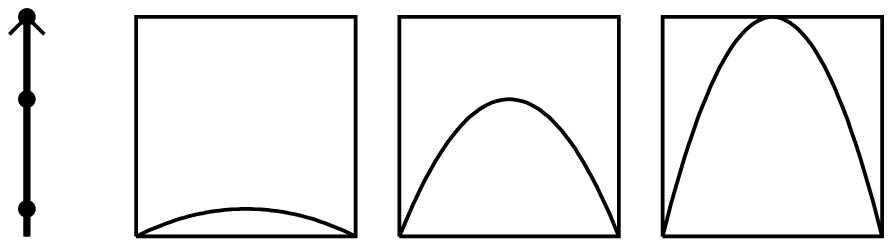,height=.8in}}

{\QP\bit Figure 5. Picture of the polyhedron $P^1$ on the left,
followed by the graphs of the quadratic maps of shape $+-$ corresponding
to three representative points of $P^1$.\smallskip}
\endinsert

%\centerline{\psfig{figure=p1.ps,height=1.9in}~~
%        \psfig{figure=p2p.ps,height=1.9in}\qquad\qquad\qquad
%        \psfig{figure=p2m.ps,height=1.9in}}\bigskip
%{\QP\bit Figure 4. Pictures of the polyhedra $V^1_{+-}$,
%        $V^2_{+-+}$ and $V^2_{-+-}$,

\midinsert
\centerline{\psfig{figure=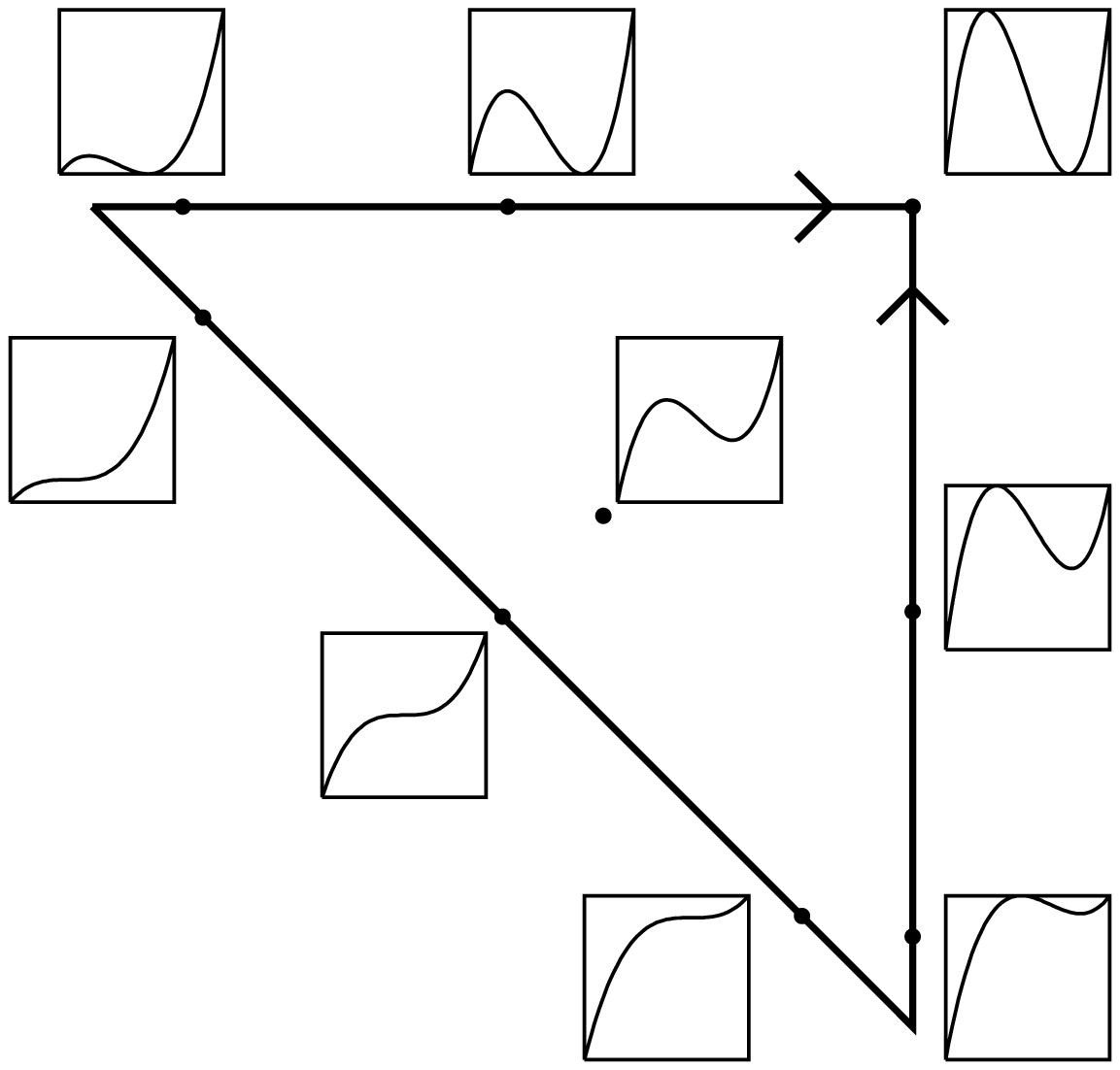,height=1.7in}\qquad
\psfig{figure=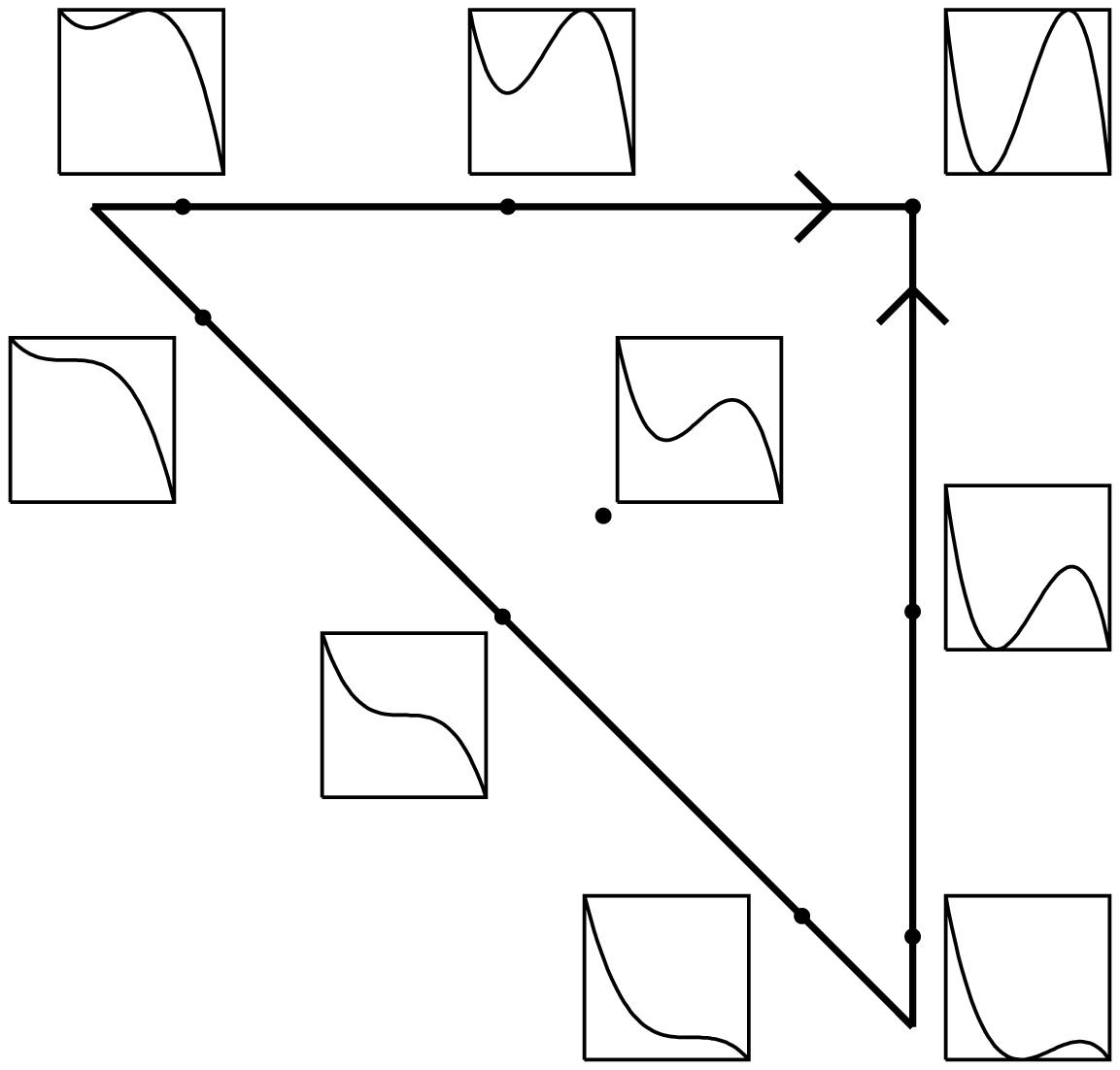,height=1.7in}}\ss
{\QP\bit Figure 6. The triangle $P^2$, showing cubic maps of shape $+-+$
on the left, and cubic maps of shape $-+-$ on the right,
for nine representative points of $P^2$. Arrows point in the direction of
the $p_1$ and $p_2$ axes.\smallskip}
\endinsert

\midinsert
\centerline{\psfig{figure=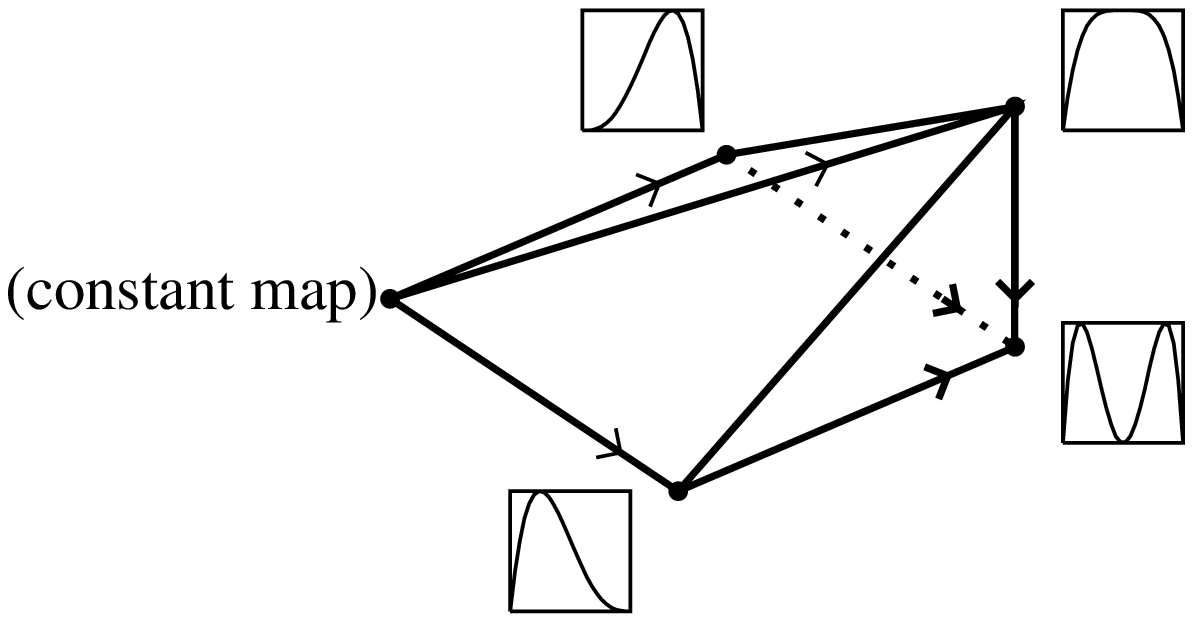,height=1.5in}}\medskip
{\QP\bit Figure 7. The polyhedron $P^3$ together with graphs
of the 4-th degree maps of shape $+-+-$ corresponding to its vertices.
Arrows point in the direction of increasing entropy, with heavy arrows
parallel to the $p_j$ axes.\bigskip}
\endinsert

\comment
We denote by $V^m_{\boldsig}=V^m_{\boldsig}(I)$ the parameter space
for the boundary anchored polynomial maps $f:I\to I$ of degree $d=m+1$
and shape $\boldsig$ parameterized by their critical value vector.
>From the discussion in this section, it follows that $V^m_{\boldsig}$
is the interval $I$ when $m=1$, and more generally,
the convex sub-polyhedron of the $m$-dimensional cube $I^m$ defined
by the inequalities $(2)$ when $m>1$. (Compare Figures 4, 5.) When
$m>3$, the polyhedron $V^{m}_{\pm}$ contains no full face
of $I^{m}$. (Here we abbreviate $\boldsig$ by
$\sigma _0$ when convenient, writing for example $V^3_+$ instead
of the more explicit notation $V^3_{+-+-}$.)
Note that the polyhedron $V^{m}_{+}$ is always
isomorphic to the polyhedron $V^{m}_{-}$, since composition with
the reflection
$r(x)=a+b-x$ carries a polynomial $g$ of shape $\boldsig$ to a polynomial
$r\circ g$ of shape $-\boldsig$. Of course this transformation is
not compatible with the dynamics. However, when $m$ is odd, the two
shapes are completely equivalent since the conjugation $g\mapsto
r\circ g\circ r$ also reverses shape.

To exploit the isomorphism between the polyhedrons $V^{m}_{+}$
and $V^{m}_{-}$, we use an alternate parametrization where
the vector $\v$ is replaced by $\p=(p_1,p_2,\dots ,p_m)$ with
$$  \eqalign{v_j\=a+(b-a)\,p_j\qquad\text{if}\qquad \sigma_j=+1~,\cr
    v_j\=b-(b-a)\,p_j\qquad\text{if}\qquad
\sigma_j=-1~.\cr} \eqno (6)$$
Thus $0\le p_j\le 1$, and the inequalities $(2)$ reduce to
% $$  \eqalign{p_j=\frac{v_j-a}{b-a}\qquad\text{if}\qquad \sigma_j=+1~,\cr
%    p_j=\frac{b-v_j}{b-a}\qquad\text{if}\qquad
% \sigma_j=-1~.\cr} \eqno (6)$$
$$p_j+p_{j+1}~\ge~1\qquad\text{for}\qquad 1\le j<m~. \eqno (2''')$$
\endcomment

As parameter space for boundary anchored polynomial maps of degree $m+1$
and shape $\boldsig$ (or indeed for any family of maps which can be
parametrized by their folding value vectors) we take the polyhedron
$P^m$ consisting of all vectors $\p=(p_1\,,\,\ldots\,,\,p_m)\in [0,1]^m$
which satisfy the inequalities (7). This
is a convex polyhedral subset of the $m$-dimensional unit cube. More
exactly,
$P^m$ can be described
as the convex hull of the set of vectors $(p_1,p_2,\dots ,p_m)$
with $p_i\in\{0,1\}$, such that there is
no consecutive pair of $0$'s. From this characterization, we get:

{\QP{\bf Proposition 3.4.} \it The polyhedron
$P^m$ has $F_{m+1}$ vertices, where $F_i$ stands for the
$i^{\text{th}}$ Fibonacci number.
\smallskip}

{\bf Proof of 3.4.} Recall that the sequence of Fibonacci numbers is
defined by $F_{n+2}=F_n+F_{n+1}$ for $n>0$, with $F_1=F_2=1$. Notice
first that $P^2$ is a segment with $F_1=1$ vertex labeled $0$ and $F_2=1$
vertex labeled $1$. Next, if $P^m$ has $F_{m-1}$ vertices whose
label terminates with $0$ and $F_m$ vertices whose label terminates
with $1$, the ``not two consecutive zeros'' rule implies that $P^{m+1}$
has $F_m$ vertices whose label terminates by $0$ and $F_{m-1}+F_m=F_{m+1}$
vertices whose label terminates with $1$.
\QED

\comment
More generally, $P^m$ appears as a ``standard parameter polyhedron''
for any family of $m$-modal maps parameterized by their critical values.
The $p_i$ might be called {\bit jaggedness parameters} or {\bit
complexity parameters} since an increase in any $p_i$ will produce
a function whose graph is more jagged, and which will have a more
complex dynamics in many cases. (Compare \S 6.)
\endcomment

While the abstract polyhedron $P^m$ is completely determined by $m$, the
way it is partitioned according to dynamical properties depends of course
on
the particular family of maps it parameterizes. Consider some family
of $m$-modal maps parametrized by $P^m$:
we will write $\p=\p(\v)$ (which can be inverted to $\v=\v(\p)$) so
that, for any chosen shape, $f_{\p}$ is the boundary anchored map
on $I$ in that family, which
corresponds to the parameter vector $\p$ and to the
folding value vector $\v$. For example, in
the polynomial case, the map $f_{(1,\ldots,1)}$
corresponding to the vector $\p=(1,1,\dots ,1)$ is, up
to sign, a Chebyshev polynomial: it is the unique boundary anchored
polynomial map of shape $\boldsig$ with maximal topological entropy.

{\bf Note:} When $m$ is odd, the two possible shapes for an $m$-modal map
are not really different, since a map of shape $(+-\cdots+-)$ with
parameters
$(p_1,\ldots,p_m)$ is topologically conjugate, under an orientation
reversing
reflection of its interval, to a map of shape $(-+\cdots-+)$ with
parameters
$(p_m,\ldots,p_1)$. However, when $m$ is even the two shapes are
dynamically essentially different. (Even the dynamical behavior of $f$
on the boundary of $I$ is enough to distinguish the two.)\bigskip

%Although the $\p$-parametrization is
%more transparent from the point of view of dynamics, we preferred
%to use first the $\v$-parametrization since
%results of this chapter belong to function theory where
%the $\v$-parametrization seems more natural.\bigskip
\bigskip

\maybebreakhere
\centerline{\bf \S 4. Topological Entropy and Periodic Orbits.}\medskip

A particularly useful measure of the dynamic ``complexity'' of a continuous
self-map is provided
by the topological entropy $h$ which was defined in [AKM] as an invariant
of topological conjugacy. The definition requires some notations.
Let $X$ be a non-vacuous compact space.
By a {\bit cover}  $\cf$  we mean simply a collection
of subsets with union $X$.
The ``covering number''  ${\bold n}(\cf)$  is defined to be the smallest
cardinality of a subcollection of $\cf$
with union  $X$. Given a continuous map $f:X\to X$  and
$k\ge 1$, define  $\cf _f^k$  to be the
cover consisting of all intersections of the form
$$C_0 \cap f^{-1}C_1 \cap \cdots \cap f^{1-k}C_{k-1}$$
where each $C_i$ is a set belonging to the collection $\cf$.
(Here $f^{-i}(C)$ is the set of all $x\in X$ with
$f^{\circ i}(x)\in C$.)
%For a continuous map $f:X\to X$ and an open or finite cover $\cf$ of $X$,
Using the inequality $\bold{n}(\cf_f^{k+\ell})\le\bold{n}(\cf_f^k)\cdot
\bold{n}(\cf_f^\ell)$, one can check that the limit
$$h(f,\cf)=\lim_{k\to \infty}\,\frac{1}{k}\log\,{\bold{n}}
(\cf ^k_f)$$
exists, with  $0\le h(f,\cf)\le\infty$. Furthermore, this limit
is equal to the infimum
$$\inf_{k>0}~{1\over k}\log\,{\bold n}(\cf^k_f)~.$$
In particular, $h(f,\cf)\le{\bold n}(\cf)$, so if
${\bold n}(\cf)$ is finite (for example if
$\cf$ is a covering by finitely many sets, or if $\cf$ is a covering of
the compact space $X$ by open subsets), then it follows
that $h(f,\cf)$ must also be finite.

The {\bit topological
entropy\/} of $f$ is defined to be the the supremum of $h(f,\cf)$ over
all open covers of $X$. (This may be infinite, even if each $h(f,\cf)$
is finite, since covers by smaller sets may yield larger values of
$h(f,\cf)$.)

In the case of maps of the interval, combining theorems of
Misiurewicz and Yomdin, we have the following.
Let $C^\infty(I,I)$ be the space of $C^\infty$
maps from a closed interval $I$ to itself, with the $C^\infty$ topology.

{\QP{\bf Theorem 4.1.} \it
The topological entropy function $$h:C^\infty(I,I)\to[0,\infty)$$
is continuous.\medskip}

{\bf Proof.} Lower semi-continuity of the entropy, for interval maps,
has been proved by Misiurewicz [Mis2]. (See also [ALM].) In fact he even
proved
lower semi-continuity for $C^0$-maps with the $C^0$-topology. Upper
semi-continuity, for $C^\infty$-maps in any dimension, has been proved by
Yomdin [Y]. Finally, the statement that $h(f)<\infty$ is true for any
$C^1$ map of a compact manifold, and follows from an easily verified
bound which takes the form
$$   h(f)~\le~\log^+\roman{max}_x|f'(x)| $$
in the 1-dimensional case.
(Here $\log^+(s)$ is defined to be the maximum of $\log(s)$ and zero.)\QED

The following immediate consequence of this theorem is essential for our
purpose:

{\QP{\bf Corollary 4.2.} \it For any $d$, the topological entropy
function is continuous on the finite dimensional compact space
consisting of all polynomial maps of the interval
with degree $\leq d$.\medskip}

{\bf Remarks.} It is essential for Yomdin's Theorem
that we work with the $C^\infty$ topology. In fact Misiurewicz and Szlenk
have shown that $h:C^r(I,I)\to[0,\infty)$ is not upper semi-continuous for
any $r<\infty$. However, their example involves a sequence of $m$-modal
maps with $m\to\infty$, converging to a bimodal limit. Misiurewicz [Mis4]
has recently proved the following statement, which is much sharper than
4.2:
{\it If ${\Cal M}^1_\ell$ is the space of $C^1$-smooth maps of the interval
$I$
which are $m$-modal for some $m\le \ell-1$, with the $C^1$-topology, then $h:
{\Cal M}^1_\ell\to[0,\log(\ell)]$ is continuous.\/}
Note that Misiurewicz' s lower semi-continuity result for entropy
is true only
in dimension one. As an example ([Mis1]), consider the family of maps
$$   M_t(x,y)\= (4xy(1-x)\,,\,ty) $$
of the unit square. Here $h(M_1)=\log(2)$, but $h(M_t)=0$ for $t<1$.
In the case of diffeomorphisms of class $C^{1+\epsilon}$,
Katok [Ka] has proved an analogous lower semi-continuity theorem in
dimension 2,
but again there are easy counterexamples in higher dimensions.
For further information about
continuity properties of the topological entropy function,
see also [Mis3] and [N], as well as [D].\medskip

Consider a family of maps $f_\p$ parameterized by some compact space $P$,
and suppose that the topological entropy function $\p\mapsto h(f_\p)$
is continuous, with values in some interval $[\alpha, \beta]$.

{\bf Definition.} For each fixed $h_0\in [\alpha, \beta]$, the
{\bit $h_0\!$-isentrope\/} for this family is defined to be the
set consisting of all parameter values $\p\in P$ for which the
topological entropy $h(f_\p)$ is equal to $h_0$.

Evidently the various
isentropes are disjoint compact subsets, with union equal to $P$.
Our goal is to show that all isentropes in the cubic family are connected.
(Compare Figure 8. In fact countably many of these isentropes
are connected regions with interior, while one, with $h_0=\log(3)$, is a
single point. It seems possible that all of the rest are simple arcs.)
\medskip

The rest of this section will outline basic results about the
topological entropy of multimodal maps, which we will need in order to
describe the structure of isentropes. First a fundamental
formula proved by Misiurewicz and Szlenk and also by Rothschild:
{\it For any piecewise monotone map, we have
$$ h(f)~=~\lim_{k\to\infty}{\log\ell(f^{\circ k})\over
k}~=~\inf_{k>0}{\log\ell(f^{\circ k})\over k}~,\eqno (8) $$
where $\ell$ is the lap number, as defined in \S1 or \S 2.\/}
(Compare [MSz], [Ro], [ALM].) It follows that
$$   0~\le~h(f)~\le~\log\,\ell(f)~. $$
In particular, $~h(f)\le\log(m+1)~$ if $f$ is $m$-modal.
%(For computation of $h$, see [BK], [BST].)

In practice, it is often more convenient to work with the quantity
$$   \gamma~=~\exp(h)~=~\lim_{n\to\infty} \root n \of {\ell(f^{\circ n})}~, $$
known as the {\bit growth number\/} of $f$. For a polynomial
map of degree $d>0$, note that this number $\gamma$ lies in the closed interval
$[1\,,\, d]$.
%an $m\!$-{\bit modal\/} map, that is for a map with $m+1$ laps,
%this number $\gamma$ lies in the closed interval
%$[1\,,\, m+1]$.
In the special case
of a piecewise linear map with $~|{\roman {slope}}|
={\roman {constant}}\ge 1$, the
growth number $\gamma$ is precisely equal to this
constant $|{\roman {slope}}|$. (Compare [MSz].\footnote
{Here is is essential that $f$ be piecewise monotone. The map
$~F(x)=\inf\{\;3\,|x-1/2^n|\;;\;n\ge 0\,\}~$ on the unit interval
%$\ell(f)<\infty$. The map $F(x)=3\,\roman{Min}_{n>0}|x-1/2^n|$
has $|\roman{slope}|=3$ almost everywhere, and yet has $h(F)=0$
since $F(x)\le x$ for all $x$.})

\midinsert
\centerline{\psfig{figure=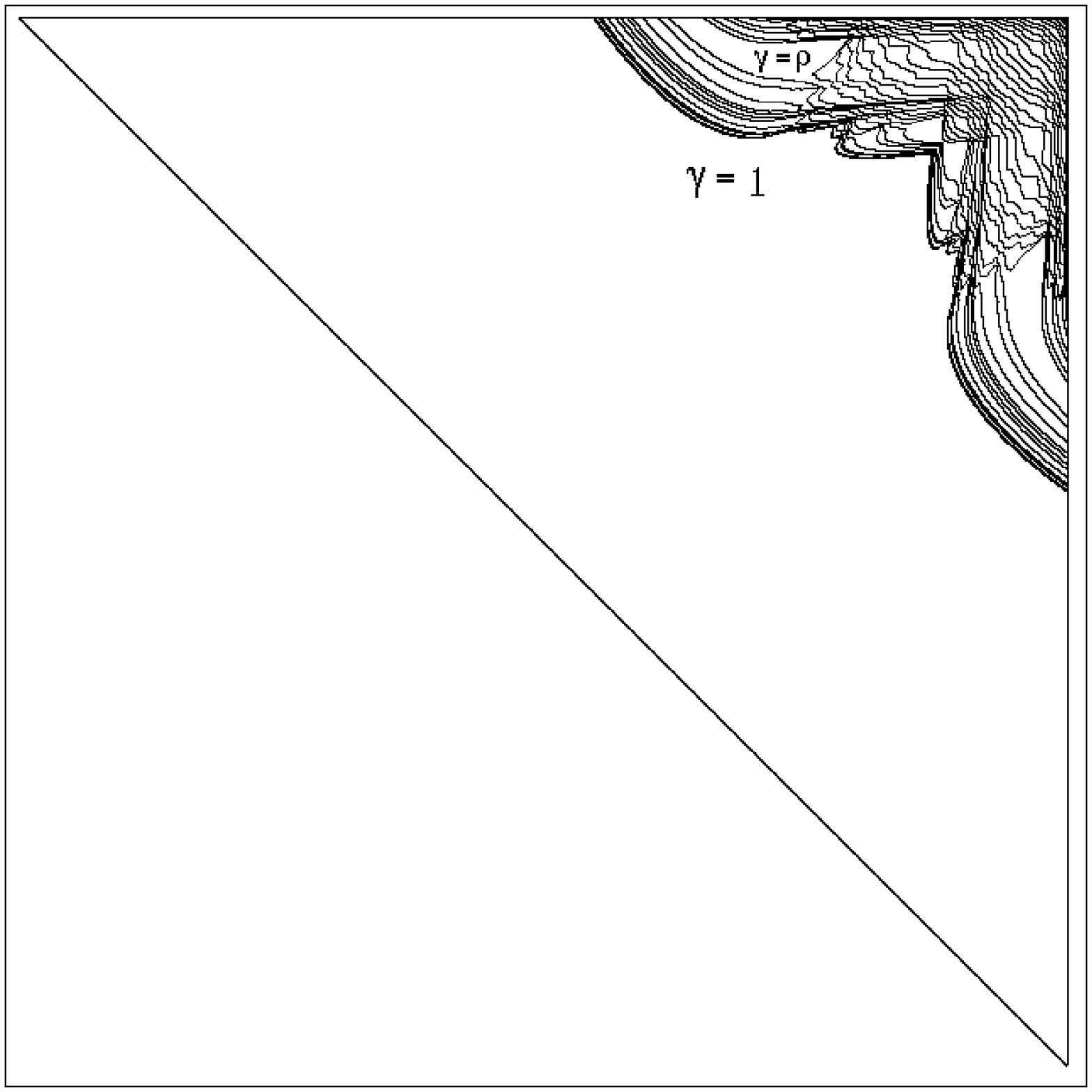,height=2.7in}\quad
     \psfig{figure=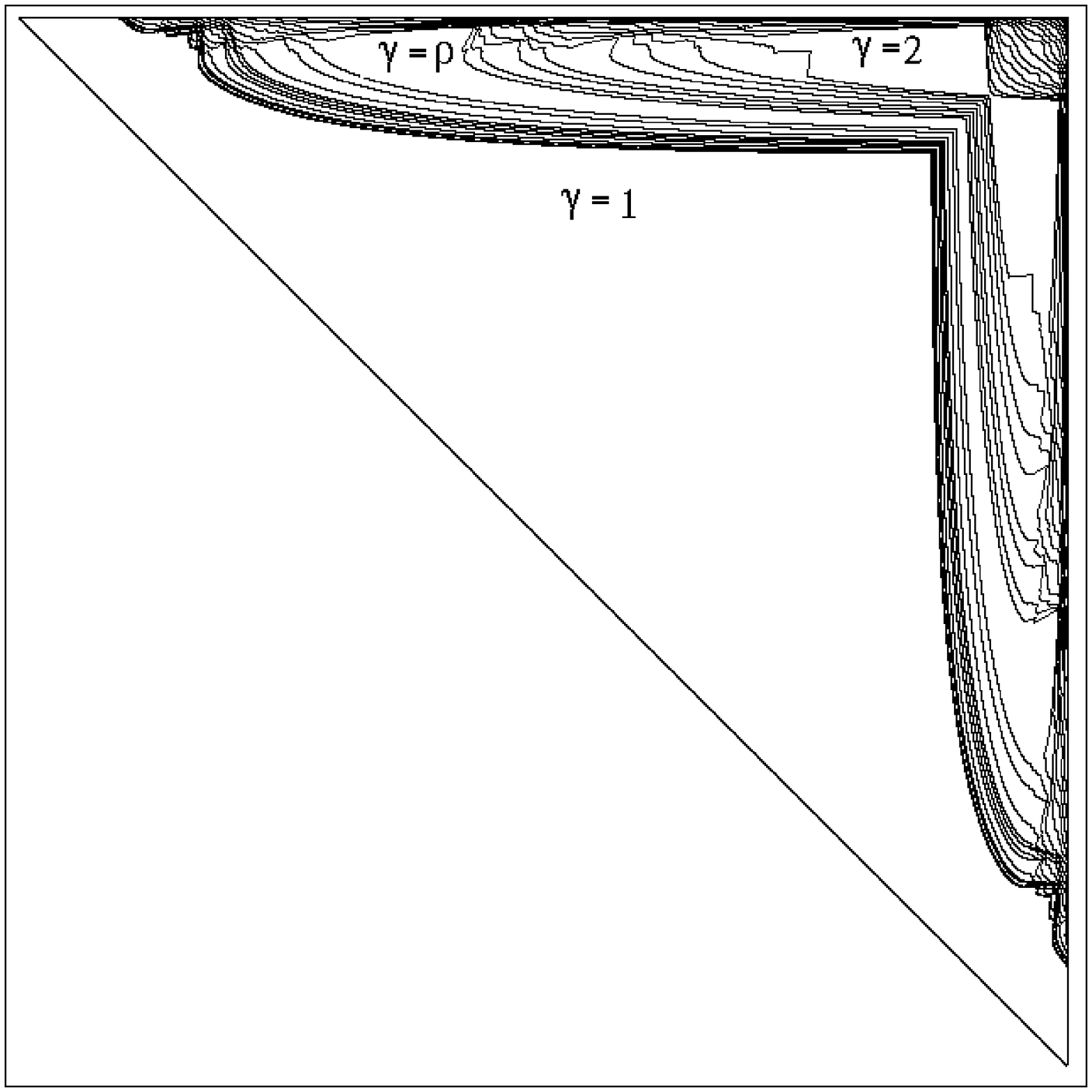,height=2.7in}}\smallskip
{\QP\bit Figure 8. Isentropes $h=\log(\gamma)=\roman{constant}$
in the parameter triangle $P^2$
for real cubic maps of shape $(+-+)$ on the left and $(-+-)$
on the right. Here $\rho$ stands for the golden ratio
$(1+\sqrt 5)/2$, associated with an attracting period 3 orbit with
both critical points in its immediate basin.
Similarly, the value $\gamma=2$ is associated with a
``capture'' component, with say $c_1$ and $f^{\circ 2}(c_2)$ in the
immediate basin of an attracting fixed point. Compare [M1, App. B].\par}
\endinsert

The proof of (8) actually yields
a slightly sharper statement as follows. Let us say that
a cover $\cf$ of the interval $I$ is {\bit $f$-mono} if $\cf$
consists of finitely many (possibly degenerate) subintervals,
and if $f$ is monotone on each of these subintervals. Evidently
such a cover exists if and only if $f$ is piecewise monotone. According to
[MSz] (or [ALM], Prop.4.2.3), we have
$$  h(f)~=~h(f,\cf)\qquad\text{ whenever the cover}~~~ \cf~~~ \text{is}~~
 f\text{-mono} ~. \eqno (9) $$
%for any $f$-mono cover, we have $h(f,\cf)=h(f)$. 
In particular, interpreting the alphabet
${\A}$ of \S 2 as an $f$-mono cover of the interval $I$, we have
$~h(f)=h(f,{\A})$. (Compare 4.4 below.)

One well known consequence of (9) is an easy computation
of $\gamma$, or of $h=\log(\gamma)$, in the
case of a map where the orbits of all of the folding points are finite.
%{\bit post-critically finite map.\/}
For such maps, the forward
orbits of the folding points cut the interval $I$ into finitely
many subintervals. Numbering these subintervals in their natural order as
$J_1\,,\,\ldots\,,\,J_n$, the associated $n\times n$ {\bit Markov transition
matrix\/} $M$ is defined by setting
$$	M_{ij}~=~\cases +1 & \roman{if}\quad f(J_i)\supset J_j~,\cr
	0 & \roman{otherwise}. \endcases $$
Recall that a real or complex number $\gamma$ is
called an {\bit algebraic integer\/} if it is a root of a monic polynomial
equation with integer coefficients, and an {\bit algebraic unit\/} if both
$\gamma$ and $1/\gamma$ are algebraic integers.

{\QP{\bf Lemma 4.3.} \it If $f$ is a piecewise monotone map such that
the orbits of the folding points are
eventually periodic, then the topological entropy of $f$ has the form
$h=\log(\gamma)$ where $\gamma$ is the largest real eigenvalue
of the associated Markov transition matrix.
In particular, $\gamma$ is an algebraic integer. Furthermore
$$	h~\le~\limsup_{k\to\infty}{\log\big(\#\roman{fix}(f^{\circ k})\big)
 \over k}~, $$
where $\#\roman{fix}
(f^{\circ k})$ denotes the number of fixed points of $f^{\circ k}$.
If all of the folding points
are actually periodic, then $\gamma$ is an algebraic unit.\smallskip}

(Compare 4.8 and 4.10.) See [BST] as well as 5.10 for methods
of computing $h$ in more general cases, based on this formula $h=\log\,\gamma$.

{\bf Outline Proof of 4.3.} It is convenient to use the sum of absolute values
norm $\|A\|=\sum|A_{ij}|$ for any $n\times n$ real or complex matrix.
If $\lambda_1\,,\,\ldots\lambda_n$ are the eigenvalues of $A$, and if $\gamma$
is the maximum of $|\lambda_i|$, then careful matrix estimates
show that
$$	\log\,\gamma~=~\lim_{k\to\infty}\big(\log \|A^k\|\big)/k~
	=~\limsup_{k\to\infty}\big(\log|\roman{trace}(A^k)|\big)/k~. $$
We will apply these equations to the Markov matrix $M$. In this case, it
follows from the Perron-Frobenius Theorem, that the number $\gamma
=\max|\lambda_j|$ is itself an eigenvalue, and hence can be described as the
largest real eigenvalue of $M$. Let us say that
a finite or infinite sequence $(J_{i_0}\,,\,J_{i_1}\,,\,\ldots)$ of
subintervals $J_{i_k}$ is $M$-{\bit admissible\/} if $f(J_{i_k})\supset
 J_{i_{k+1}}$ for all appropriate $k$.
It is easy to check that the number of sequences
$(J_{i_0}\,,\,J_{i_1}\,,\,\ldots\,,\, J_{i_\ell})$ of length $\ell+1$
which are $M$-admissible is precisely equal to the norm $\|M^\ell\|$.
These facts, together with equation (9) for the covering $\{J_i\}$,
imply that $h=\log\,\gamma$. Similarly, the number of infinite $M$-admissible
sequences which are periodic of period dividing $k$ is equal to the
trace of $M^k$. Since every $M$-admissible sequence of $J_i$ with period
dividing $k$ corresponds to at least one periodic point $x=f^{\circ k}(x)$,
and since this sequence is usually uniquely determined by $x$, the required
inequality for periodic points follows.

Now suppose that every folding point of $f$ belongs to a periodic orbit.
%, then the growth number $s$ of $f$ is actually an algebraic
%unit; that is, both $s$ and $1/s$ are algebraic integers. To prove this,
List the points $x_1<x_2<\cdots<x_{n+1}$ which belong to the orbits of
folding points. The action
of $f$ on these points is described by an $(n+1)\times( n+1)$
permutation matrix $A$, where $A_{ij}=1$ if and only if $f(x_i)=x_j$.
The associated
Markov transition matrix $M$ between the intervals $[x_i\,,\,x_{i+1}]$ can be
contructed out of $A$ in four steps, as follows. Let $A'$ be the matrix
obtained from $A$ by adding ones to the left of every one.
(For a rough bar-graph of the map, rotate this matrix $90^\circ$
counterclockwise.)
Let $A''$ be obtained from $A'$ by replacing every row $R_i$ except the first
by $R_i-R_{i-1}$, and let $A'''$ be the smaller
matrix which is obtained from $A''$ by deleting its first row and first
column. Then it is not difficult to check that all four of these
matrices have the same determinant $\pm 1$. Note that the $i$-th row of $A'''$
is negative if and only if $f$ is decreasing on the interval
$[x_i\,,\,x_{i+1}]$. The required Markov matrix $M$, with
$M_{ij}=1$ if and only if the image of the interval
$[x_i\,,\,x_{i+1}]$ contains $[x_j\,,\,x_{j+1}]$, can now be obtained
from $A'''$ by changing the sign of every negative row,
so that $M_{ij}=|A'''_{ij}|$. Since $M$ is a matrix of
integers with determinant $\pm1$, it follows that every eigenvalue
is an algebraic unit. In particular, the largest eigenvalue
$\gamma$ must be an algebraic unit. This computes the topological entropy
of $f$ restricted to the interval $[x_1\,,\,x_{n+1}]$. The actual interval
of definition for $f$ may be strictly larger than $[x_1\,,\,x_{n+1}]$.
However,
it is not difficult to check that the portion of this interval to the left
of $x_1$ or to the right of $x_{n+1}$ makes no contribution to $h$.
(Compare the proof of 4.4 below.)\QED

{\bf Remark.} Let $(L_1\,,\,\ldots\,,\,L_n)$ be an eigenvector for $M$,
so that $\sum_j M_{ij}L_j=\gamma L_i$. If the $L_i$ are strictly positive,
then we can construct a piecewise linear ``model'' for $f$, having
slope $\pm\gamma$ everywhere, as follows. Replace each $J_i$ by an interval of
length $L_i$, so that
$f(J_i)$ will be replaced by an interval of length $\gamma L_i$. Interpolating
linearly we obtain the required piecewise linear map,
which has the same kneading data as $f$.

As another application of (9), we will prove the following.

{\QP{\bf Lemma 4.4.} \it The topological entropy of an
$m$-modal map $(f,\bold{c})$ is determined by its kneading data.\smallskip}

{\bf Proof.} Assume first that $f:I\to I$ is boundary anchored.
Recall that an itinerary $\I(x)=(A_0\,,\,A_1\,,\,\ldots)$ is called
{\bit acritical\/} if each $A_i$ belongs to the subalphabet
$\A_0=\{I_0\,,\,\ldots\,,\,I_n\}$.

{\bf Definition:} For each $k>0$, let $\adm(f,k)$, or more precisely
$\adm(f,\bold{c},k)$, be the number of acritical
sequences of length $k$ which are admissible for $f$. According to
Corollary 2.4, these numbers $\adm(f,k)$ are determined by the kneading
data for $f$. On the other hand, it is not difficult to check that
$$   \ell(f^{\circ k})~\le~\adm(f,k)~\le~\roman{Card}(\A^k_f)~~. $$
(In fact $\ell(f^{\circ k})=\adm(f,k)$ in the piecewise strictly monotone
case.) Taking the logarithm of these quantities, dividing by $k$, and then
letting $k\to\infty$, we see that
$$   h(f)~\le~\lim_{k\to\infty}{1\over k}\log\big( \adm(f,k)\big)~\le~ h(f)
 $$
by (8) and (9), hence
$$   h(f)\=\lim_{k\to\infty}\,{1\over k} \log\big(\adm(f,k)\big)~. \eqno
(10) $$
Thus $h(f)$ is determined by kneading data when $f$ is boundary anchored.

Now consider the more general case where $f:I\to I$ is not boundary
anchored.
It is not difficult to extend $f$ to a boundary anchored map
$g:J\to J$ on a strictly larger interval $J\supset I$, where $g$ has the
same
shape and the same kneading data. Let $L$
and $R$ be the two connected components of $J\ssm I$. (One of these two
may be empty.) Since $g(I)=f(I)\subset I$, we see that any orbit
for $g$ either lies completely in $I$, or lies completely in $L\cup R$, or
else consists of a finite initial segment in $L\cup R$ followed by a
terminal
segment in $I$. Note also that neither $g(L)$ nor $g(R)$ can intersect
both $L$ and $R$. Hence there are only two possible sequence of any given
length in $L\cup R$. It follows that
$$   \roman{Card}(\A^k_g)~\le~2\sum_{i=0}^k\roman{Card}(\A^i_f)~. $$
%\eqno (11)$$
Now for any constant $\log(c)>h(f)$ we have $(\log\roman{Card}(\A^k_f))/k<
\log(c)$ or equivalently $\roman{Card}(\A^k_f)<c^k$, for large $k$.
Hence there is a constant $a$ so that
$~\roman{Card}(\A^k_f)<a\,c^k~$ for all $k$. It then follows that
$$   \roman{Card}(\A^k_g)~<~2a(1+c+\cdots+c^k)~<~a'\,c^k$$
for some constant $a'$. Therefore $h(g)<\log(c)$ for all such $c$. This
implies that\break $h(g)\le h(f)$, hence $h(g)=h(f)$, as
required.\QED\smallskip

Using these results, we can define a useful partial order on the possible
kneading data $\bold{K}(f)$ for $m$-modal maps of a given shape.
Let us say that $\bold{K}(f)\gg\bold{K}(g)$ if and only if
$$\eqalign{ \K_i(f)~\ge~\K_i(g)\quad&\text{when}\quad \sigma_i=-1~,\cr
    \K_i(f)~\le~\K_i(g)\quad&\text{when}\quad \sigma_i=+1~,}$$
%$$\sigma _{i-1}\K_i (f)~\geq~\sigma _{i-1}\K_i (g)$$
for $1\le i\le m$.
%(Here multiplication of a symbol sequence by $-1$ means that we
%reverse the ordering.) %\Longrightarrow h(f)\geq h(g)\,.$$

{\QP{\bf Corollary 4.5.} \it If $~\bold{K}(f)\gg\bold{K}(g)~$, then $~h(f)
\ge h(g)$.\medskip}

{\bf Proof.} If $\bold{K}(f)\gg\bold{K}(g)$, then it is
easy to check that every
admissible sequence for $g$ is also admissible for $f$. Therefore
$\adm(f,k)\ge\adm(g,k)$ (with notation as in the proof of 4.4).
Making use of formula (10), %the equality
%$h(f)=\lim_k\log(\adm(f,k))/k$, which was derived in the proof of 4.4,
the conclusion follows.\QED

This corollary will help us to reduce questions about isentropes
to questions about kneading data. More precisely, our strategy will
develop as follows: we will first verify the Connected Isentrope
Conjecture for the special family of ``stunted sawtooth''
maps, and then transfer as much as we can of parameter space
information from that family to the cubic family.\smallskip

We will see in 5.8 that topological entropy depends continuously on the
kneading data. Here is a preliminary result in that direction. We give the
space $\A^\N$ of itineraries the usual infinte product topology.

{\QP{\bf Lemma 4.6.} \it Topological entropy is upper
semi-continuous as a function of kneading data.\smallskip}

{\bf Proof.} Let $f:I\to I$ be a map with given kneading data,
and let $g_k:I\to I$ be any sequence of maps with the
same shape such that the truncations
of corresponding kneading sequences of $f$ and $g_k$ agree up to
length $k$. It follows that $\adm(f,k')=\adm(g_k\,,\,k')$ for $k'\le k$.
Given $\epsilon>0$, it follows from (10) that we can
choose $k_0$ large enough so that $(\log\adm(f, k_0))
/k_0<h(f)+\epsilon$. Then for $k\ge k_0$ we have
$$   h(g_k)~\le~(\log\adm(g_k\,,\, k_0))/k_0\=(\log\adm(f,k_0))/k_0
     ~<~h(f)+\epsilon~, $$
as required.\QED

%{\bf Remark.} It follows from [MTh] that continuity also holds
%true for kneading data of strictly monotone maps.

{\QP{\bf Corollary 4.7.} \it
Let $f$ be a piecewise monotone map whose kneading sequences are all
acritical, so that the orbit of a folding value never meets a folding
point. Then the topological entropy is continuous at $f$ under
$C^0$-deformations which move the folding points continuously, keeping
these points separated and keeping their number fixed. \medskip}

{\bf Proof.} Upper semi-continuity follows easily from Lemma 4.6.
Lower semi-continuity holds true at all continuous interval maps
by [Mis2]. (In fact here, the same result for piecewise monotone maps, as
proved in [MSz], would suffice.)
\QED\smallskip

For one-dimensional maps,
there is a close relationship between topological entropy and the existence
of periodic orbits. Let us define
$$  h_{\roman {per}}(f)~=~\limsup_{k\to
\infty}{\log\#{\roman {fix}}(f^{\circ k})\over k}~, $$ %\eqno (11) $$
where $\#{\roman {fix}}$ is the number of fixed points. (Compare 4.3.)

{\bf Caution.}
It is essential to take the {\it lim\thinspace sup\/} since the limit
may well not exist. As an example, consider the tent map $T_s(x) =
s\,\roman{max}(x,1-x)$, with slope $s=\sqrt 2$ chosen so that
the two subintervals $[\sqrt 2-1\,,\,2-\sqrt 2]$
and $[2-\sqrt 2\,,\,\sqrt 2/2]$  map to each other.
Then it is not difficult to check that
$$   \#\roman{fix}(T_s^{\circ k})\=\cases 2\quad\qquad\text{for $k$ odd}\\
     2^{k/2+1}\quad\text{for $k$ even.}\endcases $$
Hence $(\log\#\roman{fix}(T_s^{\circ k}))/k$ tends to zero as $k$ tends to
infinity through odd integers, but tends to the limit
$h_{\roman{per}}(T_s)=\log\sqrt 2$ as $k$ tends to infinity through even
integers.

{\QP{\bf Lemma 4.8 (Misiurewicz and Szlenk).} \it The inequality
$~h_{\roman{per}}(f)\ge h(f)$ is valid for any piecewise monotone map.
\medskip}

Compare [MSz]. In fact Misiurewicz later showed that this inequality is
valid for any continuous interval map (compare [Mis2] or [ALM, 4.3.14]).
However, the piecewise monotone case will suffice for our purposes.
(Compare 4.12.)
%A proof based on 4.3 will be sketched in 5.10. **

In fact, equality holds in many important cases. See 4.10 below.

We will need to subdivide periodic orbits into three classes.
Recall that every periodic point, $f^{\circ p}(x)=x$ has an
itinerary $\I(x)=(A_0\,,\,A_1\,,\,\ldots)$ which is also periodic,
with $A_i=A_{i+p}$.
Define the {\bit sign\/} of this fixed point of $f^{\circ p}$ to be the
product $\epsilon(A_0)\epsilon(A_1)\cdots\epsilon(A_{p-1})$, where
$$   \epsilon(I_j)=\sigma_j\in\{\pm 1\}\qquad\text{and}
 \qquad\epsilon(C_j)=0~,$$
as in \S2. (In the special case where $f$ is piecewise strictly monotone,
this sign is either $+1\,,\,-1$, or zero according as $f^{\circ p}$ is
increasing, decreasing, or has a folding point at $x$.)
%$C^1$-smooth, with derivative vanishing only at
% $c_1\,,\,\ldots\,,\,c_m$, this is just the sign of the
% first derivative of $f^{\circ p}$ at $x$.) 
A periodic orbit will be
called of {\bit positive, negative\/}, or {\bit critical\/} type according
as its sign is $+1\,,\,-1$, or $0$. Similarly, a periodic
sequence in $\A^\N$ has either positive, negative or folding type.
We will be particularly interested in
periodic points of negative type, since they are quite stable under
perturbation of the $m$-modal map $(f,\bold{c})$. In particular,
we have the following, with notation as in 4.5.

{\QP{\bf Lemma 4.9.} \it  Given any admissible
symbol sequence $\{A_i\}$ which is periodic of negative type, with
$A_i=A_{i+p}$, there is one and only one fixed point of $f^{\circ p}$
which has this symbol sequence as itinerary.
Hence the number $\roman{Neg}(f^{\circ p})$ of fixed
points of negative type for each iterate $f^{\circ p}$ is completely
determined by the kneading data for $(f,\bold{c})$, and satisfies
$\roman{Neg}(f^{\circ p})\le\adm(f,p)$. Furthermore, if
$~\bold{K}(f)\gg\bold{K}(g)~$, then $~\roman{Neg}(f^{\circ p})\ge
\roman{Neg}(g^{\circ p})~$ for every $p\ge 1$.
\medskip}

(Note that a fixed point of negative type for $f^{\circ p}$ counts also
as a fixed point of negative type for the odd iterates $f^{\circ 3p}\,,\;
f^{\circ 5p}\,,\;\ldots$, but as a fixed point of positive type
for $f^{\circ 2p}\,,\,f^{\circ 4p}\,,\,\ldots$.)

{\bf Proof of 4.9.} Let $J=(\alpha,\beta)\subset I$ be the subinterval
consisting of all $x\in I$ with\break
$A(f^{\circ i}(x))=A_i$ for $0\le i< p$. Then the restriction $f^{\circ
p}|_J$
is monotone decreasing with $f^{\circ p}(J)\cap J\ne\emptyset$. This
implies
that $f(\alpha)>\alpha$ and $f(\beta)<\beta$. Hence by the Intermediate
Value
Theorem,  $f^{\circ p}|_J$ has a fixed point, which is unique since this
restriction is monotone decreasing. Thus  $\roman{Neg}(f^{\circ p})$ is
equal to the number of admissible sequences which are fixed points of
negative
type for the $p$-fold iterate of the shift. It follows that
$\roman{Neg}(f^{\circ p})\le\adm(f,p)$.
In the boundary anchored case, the statement that the numbers
$\roman{Neg}(f^{\circ p})$ are completely determined by the kneading data
follows immediately from
Corollary 2.4. The general case follows since, extending $f$ to a boundary
anchored map as in the proof of 4.4, it is easy to see that no periodic
orbit of negative type can involve either of
the two intervals which are adjoined to the ends of $I$.\QED

{\QP{\bf Lemma 4.10.} \it If $f$ is piecewise monotone, with
at most finitely many non-repelling periodic orbits, then
$$	h(f)~\ge~h_{\roman{per}}(f)~=~\limsup_{k\to\infty}
	{\log^+\roman{Neg}(f^{\circ k})\over k}~. $$
\medskip}

Combining this inequality with 4.8, it follows of course that
$$	h(f)~=~h_{\roman{per}}(f)~=~\limsup_{k\to\infty}
	{\log^+\roman{Neg}(f^{\circ k})\over k} \eqno(11) $$
whenever $f$ has at most finitely many non-repelling periodic orbits.
This hypothesis is satisfied in many important cases. For example,
for a polynomial map of degree $d>1$, the
classical theory of Fatou and Julia shows that the number of such
orbits is at most $d-1$.
For a smooth\footnote {Note also the following theorem of Martens, de Melo
and van Strien [MMS]: For any $C^1$-smooth interval map with non-flat
critical points, every orbit of sufficiently high period is repelling.}
$m$-modal map with negative Schwarzian derivative, the number
of such orbits is at most $m+2$. For an $m$-modal stunted sawtooth map
with $m>1$, as studied in \S5, the number of such orbits is at most $m$.

{\bf Remark.} This close relationship between topological entropy and
periodic orbits exists only in low dimensions.
%\footnote{It is possible that the equation
%$h=h_{\roman {per}}$ is true for a $C^r\!$-{\it generic\/}
%map in any dimension. However, no proof is known. Compare [B, p. 23].}
Katok [Ka] has proved the inequality $h_{\roman{per}}\ge h$ for
2-dimensional diffeomorphisms which are $C^{1+\alpha}$-smooth.
However, Kaloshin [K] has shown that $h_{\roman {per}}$ is
infinite (and hence strictly greater that $h$) for $C^r$-generic maps
inside a ``Newhouse region'' in parameter space.
As soon as we go to higher dimensions or allow non-smooth maps, there
can be maps of positive
entropy with no periodic orbits at all. As an example, the cartesian
product
of an irrational rotation of the circle with
an arbitrary dynamical system for which $h>0$ will have $h>0$, but
no periodic orbits. For the case of non-smooth surface homeomorphisms,
see Rees [Re].

{\bf Outline Proof of 4.10.}
First suppose that all of the periodic orbits of $f$
are non-folding and strictly repelling. Then evidently the fixed points
of $f^{\circ p}$,
that is the places where the graph of $f^{\circ p}$ crosses the diagonal,
must be alternately of positive and negative type, and it follows easily
that
$$	|\#\roman{fix}(f^{\circ p})~-~ 2\,\roman{Neg}(f^{\circ p})|~\le~1~. $$
If there are $\ell$ periodic points which are either of folding type
or are non-repelling,
then a similar argument shows that
$$	|\#\roman{fix}(f^{\circ p})~-~ 2\,\roman{Neg}(f^{\circ p})|~\le~
 2\ell+1~. $$
Whenever this number $\ell$ is finite, it follows immediately that
$$	h_{\roman{per}}(f)~=~\limsup_{k\to\infty}
	\big(\log^+\roman{Neg}(f^{\circ k})\big)/ k~, $$
and since $\roman{Neg}(f^{\circ p})\le
\adm(f,p)$ by 4.9, it follows by (10) that $~h_\roman{per}~\le~ h~.$\QED

{\QP{\bf Main Theorem 4.11.} \it
For any $m$-modal map, we have
$$   h(f)\=\limsup_{k\to\infty}{1\over k}\,\log^+(\roman{Neg}(f^{\circ
k}))~.
 $$\medskip}

This is proved in [MTh], and also in [Pr]. (Either one of these
references describes explicitly how to compute these numbers $\roman{Neg}
(f^{\circ k})$ %		, and also this {\it lim sup\/}, 
from the kneading data.)

{\bf Remark 4.12.} These results are closely related. Thus
4.8 is an immediate corollary of 4.11, since
$\#\roman{fix}(f^{\circ p})\ge\roman{Neg}(f^{\circ p})$. On the other hand,
a proof that 4.11 follows from 4.8
can be sketched as follows. If $f$ has only finitely many
non-repelling periodic orbits, then 4.11 follows immediately from 4.8 and
4.10. However, we know by 4.4 that the
topological entropy $h(f)$ is completely determined by the kneading data for
$f$, and we know by 4.9 that the numbers $\roman{Neg}(f^{\circ p})$ are
completely determined by the kneading data for $f$. Hence, if we can find
just one example of a map $g$ which has the same kneading data as $f$, and
which has only finitely many non-repelling periodic orbits, then 4.11 will
follow. In fact such an example will be provided in the next section
(5.3 together with 5.4).

Another proof of both 4.8 and 4.11 will be described at the end of \S5.

Define the {\bit negative orbit complexity\/} of $f$ to be the sequence
$${\Cal N}(f)\=
 \Big(\Neg(f^{\circ 1})\,,\,\Neg(f^{\circ 2})\,,\,\Neg(f^{\circ 3})\,,\,
 \ldots\Big)$$
of non-negative integers, and define the relation ${\Cal N}(f)\gg{\Cal
N}(g)$
to mean that
$\Neg(f^{\circ p})\ge\Neg(g^{\circ p})$ for all $p$, then we can summarize
4.5, 4.9b, and 4.11 as follows.

{\QP{\bf Corollary 4.13.} \it The kneading data $\bold{K}(f)$
determines the negative orbit complexity ${\Cal N}(f)$,
which in turn determines the topological entropy $h(f)$, with
$$   \bold{K}(f)\gg\bold{K}(g)\quad \Rightarrow\quad {\Cal N}(f)\gg
{\Cal N}(g)\quad \Rightarrow\quad h(f)\ge h(g)~.\eqno(12) $$\smallskip}

\bigskip
%{\centerline {\bf $<$ To be continued $>$}}
\maybebreakhere
\centerline{\bf \S 5. The stunted sawtooth family.}\smallskip
\medskip

Closely related to kneading theory is a special family of
piecewise monotone maps which is rich enough to encompass in a canonical
way all possible kneading data and all possible itineraries.
(Compare [Gu], [BCMM], as well as Figure 9.)
In order to introduce this family, we first introduce the {\bit sawtooth
map\/} of specified shape $\boldsig=(\sigma_0,\ldots,\sigma_m)$.
This map, on an interval $J$,
can be characterized as the unique piecewise linear map $S:J\to\R$
which is boundary anchored of shape $\boldsig$, with slope $\pm s$
everywhere
where $s>m+1$ is some specified constant,\footnote
{In the preliminary publication [DGMT] (and also in 2.2)
we took $s=m+1$. However,
here we will need $s>m+1$ in order to get an actual embedding of the
set of all $\bold K$-admissible sequences in $\A^\N$ into the reals.} and
with folding points in arithmetic progression.

{\bf Caution:} This map $S$
does not carry the interval $J$ into itself. In fact
the folding values of $S$ all lie outside of $J$.

\midinsert
\centerline{\psfig{figure=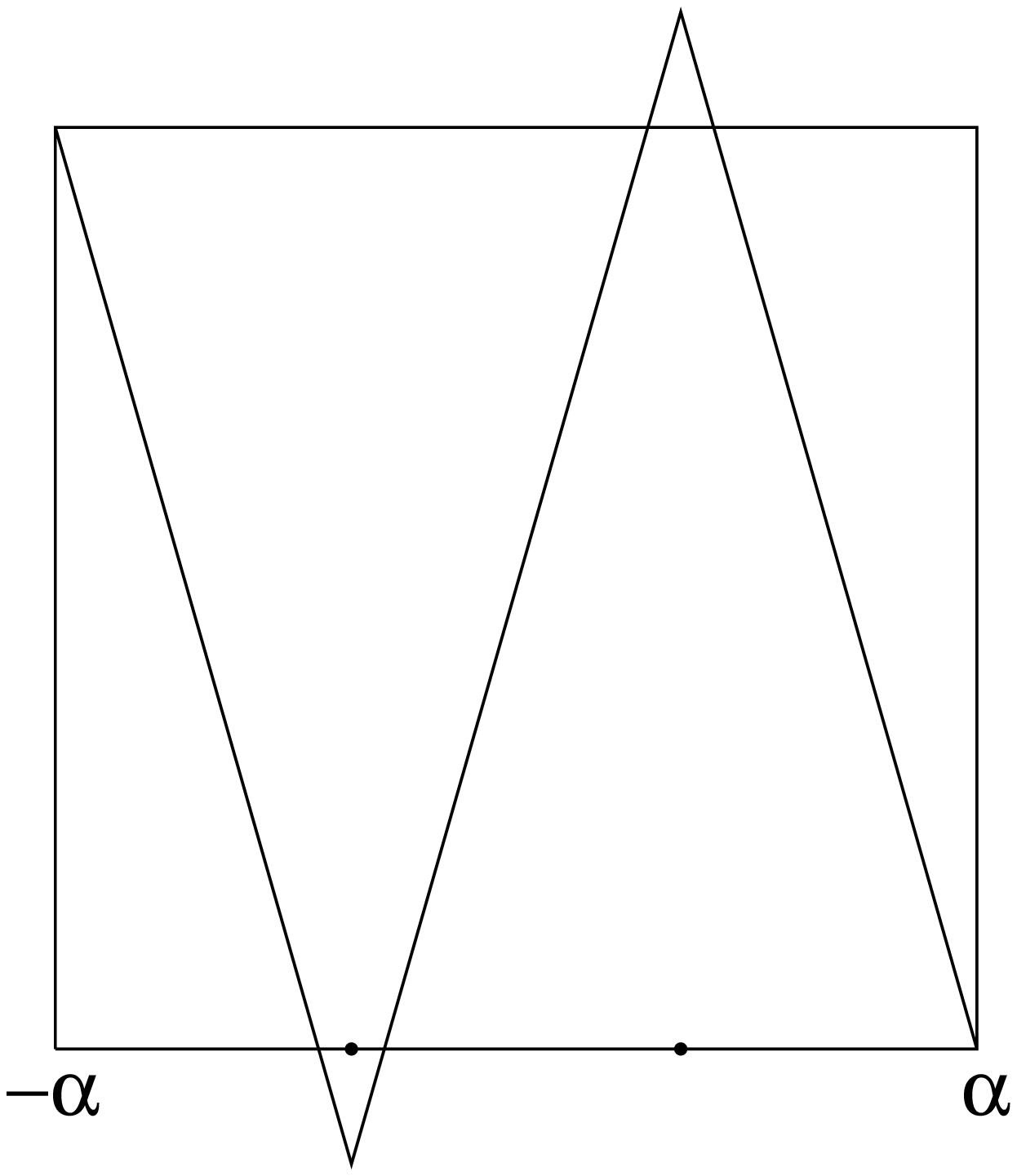,height=1.8in}\qquad\qquad
        \psfig{figure=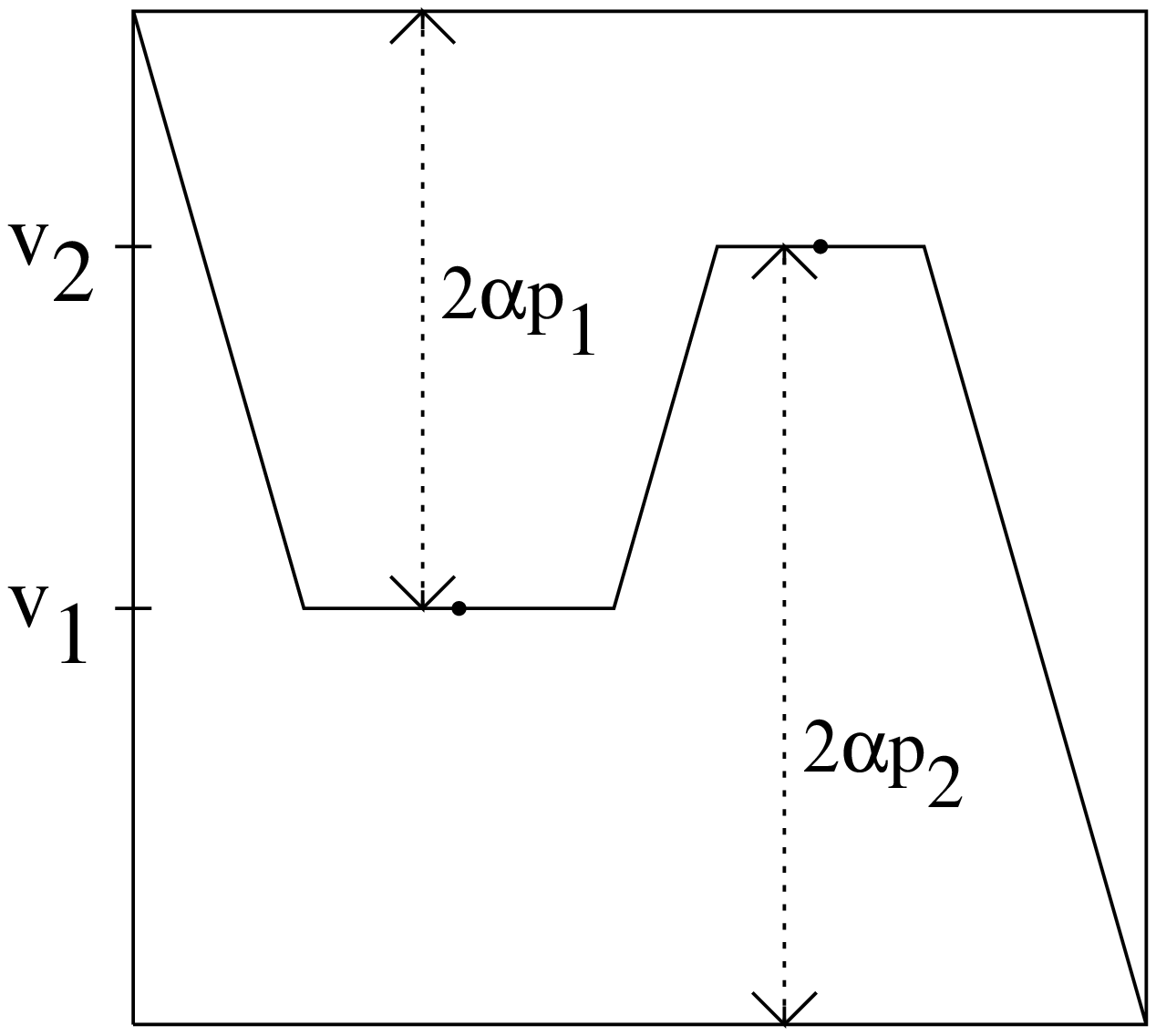,height=1.8in}}\smallskip
{\QP\bit Figure 9. The bimodal sawtooth map of shape $(-\,+\,-)$,
and a representative stunted sawtooth map of the same shape.\smallskip}
\endinsert

The precise choice of the
constant slope $s$ doesn't matter, but to fix our ideas, let us always
take $s=m+3/2$. One choice of the interval $J$ is particularly convenient
for kneading theory:
As folding points $\widehat C_1\,,\,\ldots\,,\,\widehat C_m$
choose the numbers $$-m+1\,,~-m+3\,,~\ldots\,,~m-1~,$$
with $\widehat C_j=2j-m-1$, and with $S(\hat C_j)=-\sigma_j s$ as
the corresponding folding values. Now choose a base
point $\widehat I_j=2j-m$ in the $j$-th lap, so that
$$   \widehat I_0~<~\widehat C_1~<~\widehat I_1~<~\cdots~<~\widehat
 C_m~<~\widehat I_m $$
are consecutive integers. Then $S$ is given by the formula
$$   S(x) \= \sigma_js\,(x-\hat I_j)\quad\text{for}\quad x\in I_j~.
     \eqno (13)$$
We must choose a domain of definition $J$ for $S$ so that $S(\partial
 J)\subset\partial J$. The appropriate choice is the interval $J=
[-\alpha\,,\,\alpha]$, where
$$   \alpha\={ms\over s-1}~<~m+1~<~s~. $$
Since $\alpha=s\,(\alpha-m)$, it follows
%$$ \alpha\=m\big(1+1/s+1/s^2+\cdots\big)\={m\over 1-1/s}~<~m+1~<~s~. $$
%$$  J\=[-\alpha\,,\,\alpha]\qquad\text{with}\qquad \alpha=m/(1-1/s)~. $$
%(Note that $m<\alpha<m+1<s$.)
that $S(\pm\alpha)\in\{\pm\alpha\}$,  as required.
%$S(-\alpha)=-\sigma_0\alpha$ and $S(+\alpha)=\sigma_m\alpha$,
%so that $S(\partial J)\subset\partial J$.

As in \S2, we work with the ordered alphabet
$\A=\{I_0\,,\,C_1\,,\,I_1\,,\,\ldots\,,\,C_m\,,\,I_m\}$.
The correspondence $I_j\mapsto\widehat I_j\,,\quad C_j\mapsto\widehat C_j$,
or briefly $A\mapsto\widehat A$ for $A\in\A$,
defines an order preserving embedding of this alphabet into the integers.
We can extend to a map\break $\Phi:\A^\N\to \R$ by setting
$$\eqalign{     \Phi(A_0\,,\,A_1\,,\,\ldots)~& =~\widehat
A_0\,+\,\epsilon_0
        \,\widehat A_1\,/s
        \,+\,\epsilon_0\,\epsilon_1\,\widehat A_2\,/s^2\,+\,\cdots\cr
        &=~\sum_{k=0}^\infty
        \;\epsilon_0\cdots\epsilon_{k-1}\,
\widehat A_k\,/s^k ~.}     \eqno (14) $$
(Recall that $\epsilon_i=\epsilon(A_i)$ with $\epsilon(I_j)=\sigma_j$ and
$\epsilon(C_j)=0$.) Order the space of sequences $\A^\N$ as in \S2.

{\QP{\bf Lemma 5.1.} \it The image $\Phi(\A^\N)$
is contained in the interval $J=[-\alpha\,,\,\alpha]$, and $\Phi$
is related to the sawtooth map on $J$ by the identity
$$   \Phi(A_1\,,\,A_2\,,\,A_3\,,\,\ldots)\=S\big(
 \Phi(A_0\,,\,A_1\,,\,A_2\,,\,\ldots)\big)$$
whenever $\epsilon_0\ne 0$, or in other words whenever
$ A_0\in\{I_0\,,\,I_1\,,\,\ldots\,, \,I_m\}$.
(Thus $\Phi$ semiconjugates the shift map on $\A^\N$ to the sawtooth
map on $J$, except when $A_0$ is a folding point symbol.)
This map $\Phi:\A^\N\to J$ is strictly order preserving,
provided that we consider only sequences which satisfy the first
compatibility
condition of \S2 with some fixed kneading data.
\ss}

{\bf Proof of 5.1.} Since $|\widehat A_j|\le m$, we have
$$   |\Phi(A_0\,,\,A_1\,,\,\ldots)|~\le~
m\,(1+1/s+1/s^2+\cdots)\=\alpha~,$$
hence $\Phi(\A^\N)\subset J$. We can write equation (14) as
$$   \Phi(A_0\,,\,A_1\,,\,\ldots)~=~\widehat A_0+\epsilon(A_0)\,
     \Phi(A_1\,,\,A_2\,,\,\ldots)\,/s~. \eqno(14')$$
If $\epsilon(A_0)=\pm 1$, then equation $(14')$ can be solved for
$$   \Phi(A_1\,,\,A_2\,,\,\ldots)~=~s\,\epsilon(A_0)\,\Big(
    \Phi(A_0\,,\,A_1\,,\,\ldots)-\widehat A_0\Big) ~.$$
Comparing (14), we can write the right hand side of this equation as\break
$S\big(\Phi(A_0\,,\,A_1\,,\,\ldots)\big)$, as required.
Finally, we must show that the map $\Phi:\A^\N\to J$ is strictly monotone,
in the sense that
$$(A_0\,,\,A_1\,,\,\ldots)~<~(B_0\,,\,B_1\,,\,\ldots)
\quad\Longleftrightarrow\quad
     \Phi(A_0\,,\,A_1\,,\,\ldots)~<~\Phi(B_0\,,\,B_1\,,\,\ldots)~, $$
provided that both sequences satisfy the First Compatibility Condition
of \S2 with the same kneading data. Let $k\ge 0$ be the smallest index
with $A_k\ne B_k$. It follows from $(15')$ that
$$   |\Phi(A_0\,,\,\ldots)-\widehat A_0)|~\le~ \alpha/s~<~1~ ,$$
with $\Phi(A_0\,,\,\ldots)=\widehat A_0$ when $A_0$ is a folding point
symbol.
If $k=0$ so that $A_0<B_0$, then it follows immediately from this
inequality
that
$$\Phi(A_0\,,\,A_1\,,\,\ldots)<\Phi(B_0\,,\,B_1\,,\,\ldots)~. $$
Now if $k\ge 1$, then $\epsilon(A_0)$ must be non-zero,
and an easy induction on $k$, using $(14')$, proves the corresponding
inequality.\QED

Now suppose that we are given an $m$-modal map $f:I\to I$ of shape
$\boldsig$.
For any non-folding point $x\in I\subset I$ the itinerary $\I(x)\in\A^\N$
maps to a real number
$$   \theta(x)\=\Phi(\I(x))~\in~J~, $$
which provides an invariantly defined coordinate for the point $x$
with respect to the map $f$. It follows from 5.1 together
with \S2 that $\theta:I\to J$ is a
monotone function, in the sense that
$$   x~<~y\quad\Rightarrow\quad \I(x)~\le~\I(y)\quad\Rightarrow\quad
 \theta(x)~\le~\theta(y)~.$$
(However, $\theta$ is certainly not continuous.)
Note that $\theta(c_j)=\widehat C_j$,
and that $\theta(x)=\theta(y)$
if and only if $x$ and $y$ (and hence all points between $x$ and $y$)
have the same itinerary. For any $x\in I_j$ it follows from 5.1 that
$$    \theta(f(x))~=~s\,\sigma_j\,\big(\theta(x)-\widehat I_j\big)\=
S\big(\theta(x)\big)~. \eqno (15) $$

Thus $\theta$ is a monotone but discontinuous semiconjugacy from the map
$f$
on $I$ to the sawtooth map $S$ on $J$, provided that we exclude folding
points. However, this identity $\theta(f(x))=S(\theta(x))$
definitely breaks down when $x=c_j$, since $\theta(f(c_j))\in J$ but
$S(\theta(c_j))=S(\widehat C_j)\not\in J$. To correct this problem,
we truncate the map $S$ as follows. Let us specify
some {\bit folding value vector\/} ${\v}=(v_1\,,\,\ldots
\,,\,v_m)\in J^m$, subject to the usual inequalities (2):
$$  \eqalign{ v_j~\le~ v_{j+1}\qquad &\text{when}\qquad \sigma_j=+1\cr
        v_j~\ge ~v_{j+1}\qquad &\text{when}\qquad \sigma_j=-1~~,\cr} $$
%\eqno (2) $$
for $0<j<m$.

{\bf Definition.}
The {\bit stunted sawtooth map\/} $S_\p:J\to J$ is obtained by chopping
off the successive peaks and pits of the map $S$
at the heights $v_1\,,
\,\ldots\,,\, v_m$, as shown in Figure 9.
Thus $S_\p$ is a continuous function which takes
the constant value $v_j$
throughout the largest connected neighborhood of $\widehat C_j=2j-m-1$ on
which
$$\eqalign{    S(y)~\ge~v_j\qquad&\text{if}\qquad
     S(\widehat C_j)\=+\alpha\cr
     S(y)~\le~v_j\qquad&\text{if}\qquad
     S(\widehat C_j)\=-\alpha~,} $$
but $S_\p(y)$ is equal to $S(y)$ otherwise.

Note: Here, as in \S3, we parametrize the admissible folding
value vectors $\v\in J^m$ by vectors $\p$ belonging to a standard
polyhedron $P^m\subset[0,1]^m$. Since $J=[-\alpha\,,\,\alpha]$, the
transformation $\p\leftrightarrow\v$ of (6) is given by
%parametrization (6) is given by
$$     v_j\= (2p_j-1)\,\sigma_j\,\alpha~\in~[-\alpha\,,\,\alpha]\qquad
%\Longleftrightarrow%
\text{or}\qquad   p_j~=~{1\over 2}+{\sigma_j\,v_j\over 2\,\alpha}~\in~[0,1]
 ~.
 \eqno (6')~.$$

It is not difficult to write down a more explicit formula as follows:
$$  S_\p(y)~~=~~\cases v_j &\text{if $\quad|y-\widehat C_j|\, \le\,
     1+\sigma_j\,v_j/s\quad$ for some $j$,}\\
    S(y) & \text{otherwise.}\endcases\eqno (16)
    $$
The interval of constancy
$$   \{\,y\in J~;~ |y-\widehat C_j|\, \le\, 1+\sigma_j\,v_j/s\,\} $$
will be called the $j$-th {\bit plateau\/} of $S_\p$. Since $s>\alpha\ge
|v_j|$, these plateaus always have strictly positive length. Note that the
heights of these plateaus are just the corresponding folding values.

Caution: If $v_j=v_{j+1}$ (or equivalently if $p_j+p_{j+1}=1$) then the
$j$-th and $j+1$-st plateaus have a common endpoint, and together form
a longer interval of constancy. More generally, it follows from (16)
that the distance between the $j$-th and $j+1$-st plateaus is equal to
$~|v_{j+1}-v_j|/s~\ge~0$.

By definition, the $j$-th designated folding point of $S_\p$ is
just the point $\widehat C_j$, which can also be described as the midpoint
of the $j$-th plateau.
Note that every $S_\p$ is boundary anchored, mapping
$\partial J=\{\pm\alpha\}$ into itself by the fixed map
$$ S_\p(-\alpha)\=-\sigma_0\alpha\,,\quad
S_\p(+\alpha)\=+\sigma_m\alpha~.$$

In order to relate this construction to kneading theory, we choose
the folding value vector $\v$ by setting $v_j=\Phi(\K_j)$, with $\Phi$
as in 5.1. In this way, we prove the following.

{\QP{\bf Theorem 5.2.} \it To any shape $\boldsig$ and any $m$-tuple
$\bold{K}$
of symbol sequences satisfying the Compatibility
Conditions 1, 2, 3 of \S2, there is associated a canonical
stunted sawtooth map $S_\p$ which has exactly these kneading data.
Furthermore a symbol sequence in $\A^\N$ actually occurs as the itinerary
of
some point under $S_\p$ if and only if it is admissible (i.e., satisfies
Compatibility Conditions 1 and 2).
\smallskip}

{\bf Proof.}  Let $S$ be the $m\!$-modal sawtooth
map with shape $\boldsig$, and let $S_\p$ be the associated truncated
map with folding
values $v_j=\Phi(\K_j)$. The Third Compatibility Condition guarantees
that these folding values satisfy the required inequalities $(2)$.
Notice that the orbits under $S$ and $S_\p$ are identical
as long as they do not enter a plateau of $S_\p$.
Next, using Compatibility Condition 2, it is not difficult
to show that the orbit of a folding value can enter the interior of
a plateau of $S_\p$ only at the folding point in this plateau.
It follows that the folding values $v_j$ have the same
itineraries, up to their first folding point if any occurs on the
itinerary,
under $S$ and under $S_\p$. We then use the First Compatibility Condition
to show that the full itinerary of $v_j$ under $S_\p$ is equal to the
given $\K_j$. Further details are straightforward, and
will be left to the reader. \QED

{\bf Remark.} The map $S_{\b1}$
corresponding to the vector $\p=\b1=(1,1,\dots ,1)$ is the unique
stunted sawtooth map of shape $\boldsig$ with maximal numbers of periodic
orbits and maximal topological entropy.
If $X_{\b1}$ is the set of points in $J$ whose orbit
under $S_{\b1}$ meets the various plateaus of
$S_{\b1}$ only at endpoints or midpoint, then we see
that each point of $X_\b1$ is uniquely characterized
by its itinerary, which can be either a completely arbitrary
infinite\break sequence in $\{I_0\,,\,I_1\,,\,\ldots\,,\,I_m\}^\N$, or
any finite sequence of symbols from\break $\{I_0\,,\,I_1\,,\,\ldots\,,\,I_m
\}$ followed by a folding point symbol $C_j$ (with $1\le j\le m$),
and then followed by either $\I_{\min}$ or $\I_{\max}$ according as
$\sigma_j$
is $+1$ or $-1$.
In particular, given any vector $\bold{K}=(\K_1\,,\,\ldots\,,\,\K_m)$
of symbol sequences, there
are unique points $v_j\in X_\b1$ such that the itinerary
of each $v_j$ under $S_{\b1}$ coincides with the given $\K_j$
up to the first folding point symbol (if any) in this sequence.
If the three Compatibility Conditions are satisfied, then
taking these points $v_j\in X_\b1$
to be the folding values for a stunted sawtooth map $S_{\p(\v)}$,
we easily obtain another proof of Theorem 5.2. (Compare [DGMT].)

\maybebreakhere
{\QP{\bf Corollary 5.3.} \it To any $m\!$-modal map $f:I\to I$
there is associated a canonical stunted
sawtooth map $S_\p:J\to J$ which has exactly the
same kneading data. Furthermore, the (monotone but discontinuous)
correspondence $\theta:I\to J$ semiconjugates $f$ to $S_\p$. That is
$$  \theta(f(x))~~=~~S_{\p}\big(\theta(x)\big) $$
for all $x$.\smallskip}

The proof is straightforward.\QED

To complete the discussion in Remark 4.12, we also need the following
observation.

{\QP{\bf Lemma 5.4.} \it An $m$-modal stunted sawtooth map can have
at most $m$ non-repelling periodic orbits.\smallskip}

For at most one periodic orbit can intersect any plateau. But
any period $q$ orbit which does not
meet any of the plateaus must be strictly repelling, with multiplier
$\pm s^q$.\QED

{\bf Remark 5.5.} Theorem 5.2 can be reformulated by saying that
for any shape $\boldsig$, the set of possible kneading data injects
canonically into the convex polyhedron $P^m$. We will write $\bold{K}\mapsto
\p(\bold{K})$. This allows us to replace the sometimes cumbersome
comparisons of symbol sequences by comparisons of numbers. Evidently
%If for some shape $\boldsig$, the kneading data
%$\bold{K}$ corresponds to $\p$ and $\bold{K} '$ to $\p '$, then
$$   \bold{K}~\gg~\bold{K}'\qquad\Longleftrightarrow\qquad
p_i(\bold{K})~\ge~p_i(\bold{K}')\quad
\text{for all}\quad i~,$$
with notation as in 4.5.
%$$  \K_i~>~\K'_i\quad\Longleftrightarrow\quad
%    \sigma_{i-1}p_i~>\sigma_{i-1}p'_i~. $$
%$$\sigma _i \K _i\geq \sigma _i\K '_i\Leftrightarrow p_i\geq p'_i\,.$$
\eject

{\bf Remark 5.6.} Conversely, suppose that $\p, \p'\in P^m$
satisfy $p_i\ge p'_i$ for
all $i$. Then it is not hard to show that $\bold{K}(S_\p)\gg\bold{K}(S_{\p'})$,
and it follows (or can be shown directly) that $\text{Adm}(S_\p\,,\,k)\ge
\text{Adm}(S_{\p'}\,,\,k)$ for all $k$, that ${\Cal N}(S_p)\gg{\Cal N}
 (S_{\p'})$, and that $h(S_\p)\ge h(S_{\p'})$. (Compare 4.4 and 4.13.)
Similarly, an easy argument shows that an
increase in $p_i$ can only increase the number of period $k$ orbits,
or the number of fixed points of $S_\p^{\circ k}$. (Compare the proof
of 6.1.)

{\QP{\bf Lemma 5.7.} \it The topological entropy of a stunted
sawtooth map depends continuously on its parameters, or in other
words on its vector of folding values.\ss}

{\bf Proof.} Lower semi-continuity follows immediately from [Mis2],
as noted in \S4. To prove upper semi-continuity, we note that the
topological entropy $h(f)$ depends only on the mapping $f$, and
not on which particular points
within the various plateaus are designated as folding points. However,
we can choose these folding points to be disjoint from the forward orbits
of all the folding values. Upper semi-continuity then follows from 4.6. \QED

Recall from 4.4 that topological entropy is uniquely determined by kneading
data. Combining 5.2 and 5.7 we have the following sharper form of 4.6.

{\QP{\bf Corollary 5.8.} \it Topological entropy
depends continuously on kneading data.\smallskip}

For the canonical model of 5.2 certainly
depends continuously on kneading data.
(Alternatively, a direct proof of 5.8 could be based on the methods used
in [MTh, Lemma 12.3].) \QED

{\bf Definition.} Recall that a polynomial map is called ``hyperbolic'' if
every critical point lies in the basin of some periodic attractor.
To simplify the analogous discussion for a stunted sawtooth map $S_\p$, we
consider only the case where $S_\p$ is
strictly $m$-modal; that is, we assume that consecutive folding values are
distinct. It is not
hard to see that a periodic orbit for such a map $S_\p$ is attracting,
and remains attracting under perturbation of the map, if and
only if it contains an
interior point of some plateau, and hence actually absorbs all orbits in a
neighborhood. Let us call a strictly $m$-modal $S_\p$ {\bit hyperbolic\/} if
the forward orbit of each folding value eventually lands in the interior
of some plateau. Clearly this is an open
condition. The generic hyperbolicity
property of stunted sawtooth maps can be stated as follows.

{\QP{\bf Lemma 5.9.} \it The stunted sawtooth maps which are strictly
$m$-modal and hyperbolic form a dense open set in the space of all stunted
sawtooth maps of specified shape.
\ss}

{\bf Proof of 5.9.} Openness is clear. Furthermore, it is clear that the
strictly $m$-modal maps are dense. Suppose inductively that every $S_\p$
can
be approximated by a strictly $m$-modal $S_\bold{q}$ for which the orbits
of the
first $k-1$ folding values eventually hit the interior of some plateau.
Then
we can choose some $\epsilon>0$ so that these conditions will remain true
as we change the $k$-th folding value throughout an
$\epsilon$-neighborhood.
Now choose an integer $n$ so that the product $s^n\epsilon$ is greater that
the length of the entire interval $J$. As we move the folding value $v_k$
with unit speed through its $\epsilon$-neighborhood,
its forward image $S_\bold{q}^{\,\circ n}$ will move with speed
$\pm s^n$, so long as its orbit does not pass through any plateau. But we
have chosen $n$ with $s^n\epsilon$ large enough so that this is impossible.
Hence some forward image must pass though a plateau; which completes the
inductive construction.\QED\ss

Combining 4.3 with the proof of 5.9, we have the following basic result.
%(Compare [BST], which is closely related.)

{\QP{\bf Theorem 5.10.} \it The topological entropy of a piecewise monotone
map can be effectively computed to any required degree of accuracy
from its kneading data.\ss}

{\bf Proof.} We must produce computable upper and lower bounds for $h(f)$,
arbitrarily close to each other. The proof will be by induction on the number
$m$ of folding points. To construct a lower bound,
first suppose that the kneading data $\bold{K}=\bold{K}(f)$
satisfies $\K_i\ne \K_{i+1}$ for every $i$, so that the associated stunted
sawtooth map $S_\p$, with $\p=\p(\bold{K})$ as in 5.5,
has distinct adjacent folding values. Then using the argument above we
see easily that there exists
a hyperbolic map $S_{\p'}$ so that $\p'$ is arbitrarily close to
$\p$, and with $p'_i\le p_i$ for all $i$. Here $h(S_{\p'})$ is effectively
computable by 4.3, and can be chosen arbitrarily close to $h(f)=h(S_\p)$
by 5.7. On the other hand, if $S_\p$ has two adjacent plateaus at the same
height, then by ignoring the two corresponding critical points it can be
considered as an $(m-2)$-modal map, and the conclusion follows by induction.
A completely analogous argument produces a computable upper bound
$$	h(f)~\le~h(S_{\p''}) $$
which is arbitrarily close to $h(f)$. One simply chooses
$\p''$ close to $\p$, with $p''_i\ge p_i$ for all $i$, so that the orbit
of each critical value is eventually periodic. In fact if $p_i=1$, then the
associated orbit is already eventually periodic, while if $p_i<1$ then
arguing as in 5.9 we can choose $p''_i$ so that its orbit eventually hits
a plateau.\QED

{\bf Remarks.} A closely related algorithm, more explicitly worked out,
is described in [BST]. For the special case of a bimodal map,
Block and Keesling [BK] have given a fast algorithm, based on comparison with
maps of $|\roman{slope}|=\roman{constant}$. (This
was used for the plots in Figure 8.) Note that the question of effective
computability of entropy for more general dynamical systems
remains open: compare [HKC].

As a corollary, we can give another proof of two basic results from \S4.

{\bf Proof of 4.11 and 4.8.} We first show that
$$	h(f)~=~\limsup_{k\to\infty}\log^+\big(\roman{Neg}
	(f^{\circ k})\big)/k $$
for every piecewise monotone map. This statement is certainly true
for the approximating maps $S_{\p'}$ and $S_{\p''}$ of 5.10, by Lemma
4.3, together with 4.10 and 5.4. Since
$$	{\Cal N}(S_{\p'})~\ll~{\Cal N}(S_p)={\Cal N}(f)~\ll~{\Cal N}(S_{\p''})
 $$
by 4.9 and 5.6, the corresponding statement for $f$ follows. The
inequality $h_\roman{per}(f)\ge h(f)$ is an immediate corollary.\QED
\bigskip

\maybebreakhere
\centerline{\bf \S 6. Contractibility of Isentropes for the
Stunted Sawtooth Family.}
\medskip

We consider the family of stunted sawtooth maps of some specified
shape\break
$\boldsig= (\sigma_0\,,\,\ldots\,,\,\sigma_m)$, which remains fixed
throughout this section. Parameterize this family by the standard
polyhedron
$P^m$, consisting of all vectors
$$  \p\=(p_1\,,\,\ldots\,,\,p_m)~\in~ [0,1]^
m\quad\text{satisfying}\quad
p_i+ p_{i+1}\ge 1\quad\text{for}\quad 1\le i<m~~,$$
as described in \S3.
Let $h^\roman{saw}_{\boldsig}:P^m\to[0\,,\,\log( m+1)]$ be the {\bit
topological
entropy function}
$$  h^\roman{saw}_{\boldsig}(\p)~~=~~h_{\text{top}}(S_{\p}\,)~~. $$
These functions are continuous by Lemma 5.7.

{\bf Caution:} %Most of the details which follow depend only on $m$.
%However, the
This function $h^\roman{saw}_{\boldsig}:P^m\to[0\,,\,\log(m+1)]$
definitely
depends on the
choice of shape $\boldsig$, and also on the fact that we are working with
the family of stunted sawtooth  maps with specified slope $s>m+1$, rather
than
some other family. For $m$ odd, the two
choices of $\boldsig$ are related by a canonical involution of $P^m$.
However for $m$ even it is important to realize that the two possible
choices
of shape yield families which are essentially
different, and have quite different topological entropy functions.

We will use the notation $\{\p\in P^m\,;\,h^\roman{saw}_{\boldsig}(\p)=
h_0\}$, or briefly
$\{h^\roman{saw}_{\boldsig}=h_0\}$, for the {\bit isentrope\/}
% $h^\roman{saw}_{\boldsig}^{-1}(h_0)\subset P^m$,
consisting of all
$\p\in P^m$ with topological entropy $h^\roman{saw}_{\boldsig}(\p)$ equal
to $h_0$.
Similarly we sometimes write $\{h^\roman{saw}_{\boldsig}\le h_0\}$ for the
 compact subset
$\{\p\in P^m~;~h^\roman{saw}_{\boldsig}(\p)\le h_0\}$. The object of
this section is to prove the following.

{\QP{\bf Theorem 6.1.} \it For each $h_0\in[0\,,\,\log (m+1)]$ the
isentrope $$~\{\p\in P^m~;~h^\roman{saw}_{\boldsig}(\p)= h_0\}~$$ is
contractible.\medskip}

\comment ((((((((((((((((((
[An example of a non-contractible set with this property is provided
by the closure of the set of points $(x\,,\,\sin(1/x))$
with $0<x\le 1$.
We conjecture that the set $\{h=h_0\}$ is actually
contractible, but have not been able to prove it. For comparison,
we mention that the zero isentrope in the sine family of lifts of
circle maps $f_{a,b}(\theta)=\theta +a +{b\over 2\pi}
\sin (2\pi\theta)$, $(a,b) \in {\R}\times {\R}^+$, is not
contractable. More precisely the
boundary of the zero isentrope in the parameter space
has a comb structure (compare [FT]).]
\endcomment %%%%%%%%%%%%%%%%%%%%

The proof will be based on two lemmas.
We first construct a partial ordering $\>$
of the polyhedron $P^m$
so that if $\p\>\q$ then the corresponding stunted sawtooth map
$S_{\p}$ has at least as many periodic
orbits as does $S_{\q}$. By 4.10, this will imply that $h^{\roman
saw}_{\boldsig}(\p)\ge h^\roman{saw}_{\boldsig}(\q)$.
Within the interior of $P^m$ we can simply say
that $\p\>\q$ if and only if
$p_j\ge q_j$ for all $j$. (As any coordinate $p_j$ increases, any
periodic orbit which intersects the $j$-th plateau deforms continuously,
while any other periodic orbit remains unchanged.)
However, to take care of implications
on the boundary we must give a more complicated definition. It will be
convenient to say that $\p$ contains a {\bit level block\/} of length
$\ell$ if there are indices $1\le i_0\le j_0\le m$ with $\ell=j_0-i_0+1$
so that
$$  p_i+p_{i+1}~~=~~1\qquad\text{for}\qquad i_0\le i<j_0~~.$$
This means that the corresponding folding value
vector $\v=\v(\p)$ satisfies
$$   v_{i_0}=v_{{i_0}+1}=\cdots=v_{j_0}~, $$
so that the corresponding map $S_{\p}$
has $\ell$ consecutive
plateaus at the same level. (Here it is convenient to
allow the uninteresting case $\ell=1$.) We now define
$\>$ to be the smallest transitive relation on $P^m$ which satisfies the
following condition: {\it If $\p$ coincides with $\q$ except that
$\p$ contains a level block
$$      p~,~1-p~,~p~,~\ldots~,~1-p~,~p$$
of odd length $\ell\ge 1$ which is replaced
in $\q$ by a corresponding block
$$      q~,~1-q~,~q~,~\ldots~,~1-q~,~q$$
of the same length, where $p\ge q$, then $\p\>\q$.}
%\eject

{\QP{\bf Lemma 6.2.} \it For each $0\le h_0\le \log(m+1)$, the isentrope
$\{h^\roman{saw}_{\boldsig}=h_0\}$ is a deformation
retract of the region $\{h\le h_0\}$.\ss}

{\bf Proof of 6.2.} Starting with any
point $\hat\p\in P^m$ we construct a topological entropy
increasing path $t\mapsto\p(t)
\in P^m$ for $0\le t\le 1$ with $\p(0)=\hat \p$ and $h^{\roman
saw}_{\boldsig}(\p(1))=
\log(m+1)$. In fact, let
$$   p_j(t)~~=~~\min(\hat p_j+t\,,\,1)~~. $$
Alternatively, this deformation can be described
by the differential equation
$$  {dp_j\over dt}~~=~~\cases +1 & \text{if $~~p_j<1$}\\
    ~0 & \text{if $~~p_j=1$~~.}\endcases
$$
Clearly the resulting path depends continuously on $\hat\p$ and
$t$, takes values in $P^m$, and satisfies $~\p(t)\<\p(t')~$
whenever $t\le t'$. In particular,
it follows that\break $h^\roman{saw}_{\boldsig}(\p(t))\le h^{\roman
saw}_{\boldsig}(\p(t'))$ whenever $t\le t'$.

It will be convenient to use the norm
$$   \|\p-\q\|~~=~~\max_j\,|p_j-q_j|$$
for $\p\,,\,\q\in P^m$.
If $\|\hat\p-\hat\q\|<\epsilon$, then
clearly $\|\p(t)-\q(t)\|<\epsilon$ for
$t\in [0,1]$, and furthermore
$$   \p(t+\epsilon)~\>~\q(t)~,\qquad \q(t+\epsilon)~\>\p(t)~~.$$

Now suppose that $h^\roman{saw}_{\boldsig}(\hat p)\le h_0$.
Let $t(\hat\p)$ be the smallest value of
$t\in[0,1]$ with $h^\roman{saw}_{\boldsig}(\p(t))= h_0$. If
$\|\hat\p-\hat\q\|<\epsilon$, then it
follows that $|t(\hat\p)-t(\hat\q)|<\epsilon$.
Therefore, the homotopy which is defined by
$$  h_t(\hat\p)~=~\cases \p(t)& \text{if $~~t\le t(\hat\p)$}\\
    \p(t(\hat\p))&\text{if $~~t\ge t(\hat\p)$ }\endcases
$$
yields a deformation retraction from $\{h^\roman{saw}_{\boldsig}\le h_0\}$
 onto the
isentrope $\{h^\roman{saw}_{\boldsig}= h_0\}$.\QED

{\QP{\bf Lemma 6.3.} \it The region $\{h^\roman{saw}_{\boldsig}\le
h_0\}\subset P^m$
is contractible.\ss}

{\bf Proof of 6.3.} The first step is to construct a topological
entropy decreasing
deformation\break $(\hat\p\,,\,t)\mapsto \p(t)$ which continuously
``flattens out the bumps'' on
the corresponding stunted sawtooth maps $S_{\p(t)}$.
This deformation is defined by
a differential equation as follows. Let $p\,,\,1-p\,,\, p\,,\,\ldots$
be any maximal level block contained in $\p$, and let $\ell\ge 1$ be
its length. We set
$$  {dp\over dt}~~=~~\cases -1& \text{if $\ell$ is odd
 and $p>0$, and}\\
    ~0& \text{if $\ell$ is even or $p=0$~.}\endcases
$$
(Thus the majority of the entries in a block of odd length move down,
but a minority of the entries move up.)
Evidently topological entropy decreases monotonically along each path.
Furthermore the quantity $p_1+\cdots+ p_n$, which is linearly
related to the total variation of the map
$S_{\p}$, decreases with derivative $\le -1$
until we reach a stationary state, corresponding to a monotone map
with topological entropy equal to zero. This proves that
the sub-polyhedron consisting
of such stationary states is a deformation retract of the region
$\{h^\roman{saw}_{\boldsig}\le h_0\}$. In the case $m$ odd, this
sub-polyhedron consists
of a single point, corresponding to the constant map. This completes
the proof of 6.3 for $m$ odd.

For $m=2k$, this sub-polyhedron is a $k$-simplex, consisting of all
$$ (p_1\,,\,1-p_1\,,\,p_3\,,\,1-p_3\,,\,\ldots\,
,\,p_{m-1}\,,\,1-p_{m-1}) $$
with $p_1\le p_3\le\cdots\le p_{m-1}$. Since this simplex is itself
contractible, this completes the proof of 6.3. \QED

Theorem 6.1 now follows, since a retract of a
contractible space is clearly contractible
(if $F_\lambda$ is a contraction of the total space, $r$ a retraction,
and $e$ the embedding of the retract in the total space,
then $r\circ F_\lambda \circ e$ is a contraction of the retract).\QED
\bigskip

\maybebreakhere
\centerline{\bf \S 7. Bones in $P^2$.}
\medskip

The stunted sawtooth families are well understood, but the polynomial
families of the same shape are poorly understood beyond the unimodal case.
Our analysis of the relationship between these two families in the cubic
case will rely
on the study of parameter points for which at least one of the two critical
orbits is periodic.

\comment ((((((((((((((((((

By definition, a polynomial map is {\bit post-critically finite\/}
if the orbit of every critical point is periodic or eventually periodic.
Similarly, we will say that a
piecewise monotone map has {\bit finite folding point
orbits\/} if the orbit of every folding point is periodic or eventually
periodic.

{\QP{\bf Thurston Uniqueness Theorem for Real Polynomial Maps.} \it
A post-critically finite real polynomial map with distinct real critical
points, is uniquely determined by its kneading data, up to positive
affine conjugation.\smallskip}

For a full existence and uniqueness statement, together with some
indication
of the proof, see the discussion in Appendix B. There is an
analogous statement for the stunted sawtooth family, as follows. As usual,
we assume that the slope $s>m+1$ has been fixed.

{\QP{\bf Lemma 7.1.} \it Given any
piecewise strictly monotone $m\!$-modal map which is
post-critically finite, there exists one and
only one stunted sawtooth map of the same shape
with the same kneading data.\smallskip}

{\bf Caution.} Here it is essential that the kneading
data comes from a piecewise {\bit strictly\/} monotone map. Here is an
example. Recall that the ``critical points'' of a stunted sawtooth map
are defined to be the midpoints $\widehat C_i$ of the plateaus.
Let $S_\p$ be a stunted sawtooth map of shape $(+-+)$ with
critical values $v_1=\widehat C_2+\epsilon>v_2=\widehat C_2$. If
$\epsilon>0$
is sufficiently small, then this map is post-critically finite,
with kneading data
$$   \K_1\=(I_2\,,\,C_2\,,\,C_2\,,\,\ldots)\,,\qquad \K_2\=(C_2\,,\,
C_2\,,\,\ldots)~. $$
Thus there is an entire 1-parameter family of stunted sawtooth maps with
the same post-critically finite kneading data. Yet no polynomial map, and
more generally no piecewise strictly monotone map, can have this kneading
data, which requires that an entire interval in $I_2$ must map to the
second
critical point. There is a simple criterion to establish that the kneading
data for some given post-critically finite
map $f$ can be realized by a map which is
piecewise strictly monotone. Let $\{y_1\,,\,\ldots\,,\,y_n\}$ be the union
of the forward orbits of the critical points, where $y_1<\cdots<y_n$. If
$f(y_i)\ne f(y_{i+1})$ for $1\le i<n$, then interpolating linearly between
the
$y_i$, and extending suitably over the rest of the interval $[a,b]$,
we easily obtain a piecewise {\it strictly\/} monotone map with
the same kneading data.\smallskip

{\bf Proof of 7.1.} Existence follows from 5.2. To prove uniqueness of
the resulting stunted sawtooth map $S_\p$, note
that no critical orbit of $S_\p$
can land in a plateau, except at its midpoint $\widehat
C_j$. For otherwise, any map with the same kneading data would have
an interval of points eventually mapping to the same point,
which would contradict the hypothesis. Therefore, orbit points are uniquely
determined by their itineraries, hence $S_\p$ itself
is uniquely determined.\QED
\endcomment %%%%%%%%%%%%%%%%%%%%%%%%%%

We will need a precise terminology for describing periodic orbits.
As in \S1, let $\O$ be a cyclic permutation of the integers
$\{1,2,\ldots,q\}$.
By definition, a periodic
orbit $\Cal{O}=\{x_1\,,\,\ldots\,,\,x_q\}$ of an interval map $f$ is said
to
have {\bit order type\/} $\O$ if $f$ maps each $x_i$ to $x_{\O(i)}$,
where $x_1<x_2<\cdots<x_q$.
We will sometimes use the notation $\O=(i_1\,i_2\,\ldots i_p)$ for the
permutation which satisfies
$\O(i_j)=i_{j+1}$, so that
$$   f(x_{i_j})\=x_{i_{j+1}}~.$$
Here the subscripts $j$ are to be taken modulo $p$.

{\bf Definition.} A cyclic permutation
$\O$ will be called {\bit $m$-modal of shape $\boldsig$\/}
if there exists an $m$-modal map of shape $\boldsig$ which has a periodic
orbit with order type $\O$. It will be called {\bit strictly
$m$-modal of shape $\boldsig$\/} if there exists such a map with all
$m$ of its critical points on this orbit. Equivalently, this means that
there are integers
$$1~\le~\gamma_1~<~\gamma_2~<~\cdots~<~\gamma_m~\le~q~, $$
where $x_{\gamma_1}<\cdots<x_{\gamma_m}$ are to be the critical points,
with the following property. Setting $\gamma_0=1\,,\,\gamma_{m+1}=q$,
the restriction of the permutation $\O$ to the integers in each
interval $[\gamma_i\,,\,\gamma_{i+1}]$ must be either monotone increasing
or monotone decreasing according as $\sigma_i$ equals $+1$ or $-1$.
Of course this condition implies that the period $q$ must satisfy
$q\ge m$. Given any periodic orbit of this order type, for a map of shape
$\boldsig$, note that the address of the orbit
point $x_i$ is necessarily given by
$$   A(x_i)\=\cases I_0\quad if\quad i<\gamma_1\\
     I_j\quad if\quad \gamma_j<i<\gamma_{j+1}\\
     I_m\quad if\quad \gamma_m<i~.\endcases $$
%$$  A(x_i)\=I_j\qquad\text{whenever}\quad \gamma_j<i<\gamma_{j+1}~,$$
%with $A(x_i)=I_0$ for $i<\gamma_1$ and $A(x_i)=I_m$ for $i>\gamma_m$.

However, when $i$ is precisely equal to $\gamma_j$ the address $A(x_i)$
is not uniquely determined: it can be $C_j$ but it can also be either
of the adjacent intervals $I_{j-1}$ or $I_j$. Thus altogether there are
$3^m$
distinct possibilities. It follows that there are $3^m$ different
possibilities for the {\bit kneading type\/} of the periodic orbit (that
is, the
itinerary of a representative point, which is well defined up to a shift).
Compare the proof of 7.1 below.

It is not difficult to check
that a bimodal order type of shape $(+-+)$ is {\it strictly\/}
$(+-+)$-bimodal
if and only if its period satisfies
$q\ge 2$. On the other hand, a bimodal order type of shape $(-+-)$
is strictly $(-+-)$-bimodal if and only if $q\ge 3$.

\midinsert
\centerline{\psfig{figure=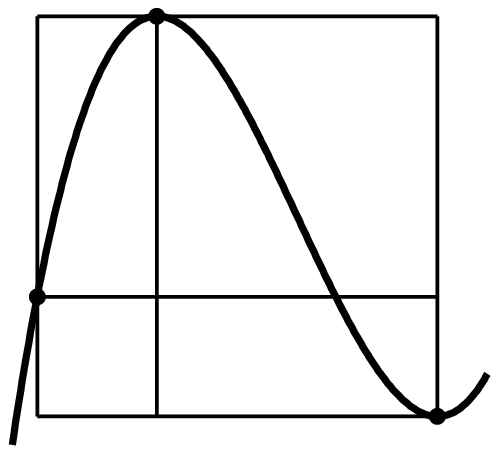,height=1in}\qquad
        \psfig{figure=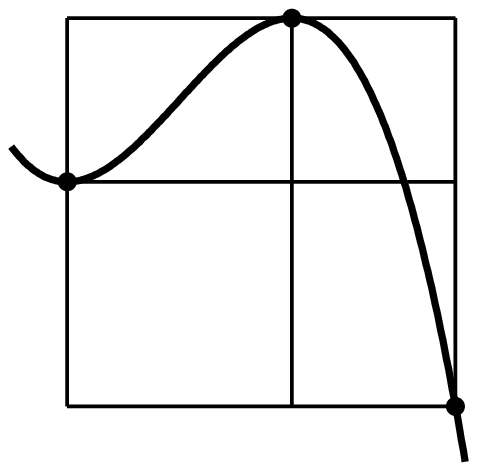,height=1in}\qquad
        \psfig{figure=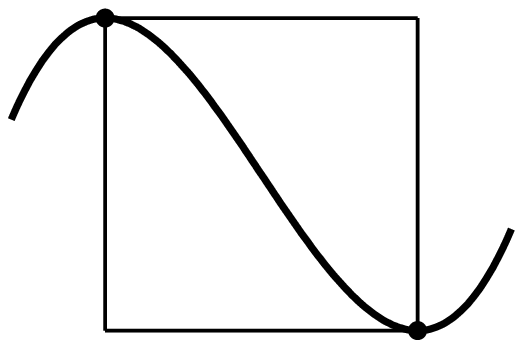,height=1in}\qquad
        \psfig{figure=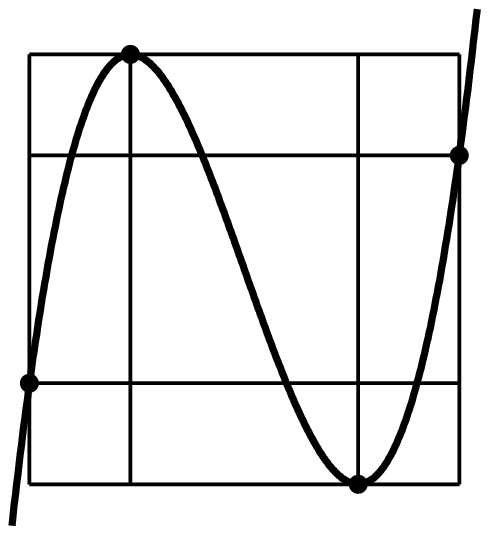,height=1in}}
\smallskip
{\QP\bit Figure 10. The period 3 order type ${\bold o}=(123)$ is strictly
bimodal of
shape either $(+-+)$ or $(-+-)$, but the period 2 order type is strictly
bimodal only of shape $(+-+)$, and the period 4 order type ${\bold o}
=(1243)$ is bimodal only of shape $(+-+)$.\medskip}
\endinsert

In order to relate two families of the same shape in the bimodal
case, we introduce some terminology from MacKay
and Tresser [MaT1]. As in \S3, consider a family of bimodal maps of
some fixed shape $\boldsig$, parameterized by the triangle $P^2$.

{\bf Definition.} By a {\bit bone} in the parameter
space $P^2$ we mean the compact set consisting of all parameter
values for which a specified critical point has periodic orbit with
specified order type. More precisely, the {\bit left bone\/} $B_-(\O)$ is
the set of
parameter values for which the left hand critical point is periodic with
order type $\O$. The dual {\bit right bone\/}
$B_+(\O)$ is the set of parameter values
for which the right critical point is periodic with this same order type.
Note that two left bones, or two right bones, are disjoint,
almost by definition.
These definitions make sense either for the stunted sawtooth family or
for the cubic family.
We will insert the superscript ``saw'' respectively
``cub'' in order to distinguish these two cases. Similarly, we will
use the notation $P^{\thinspace\roman{saw}}$ or $P^{\,\roman{cub}}$ for the
 parameter
triangle $P^2$ when we want to emphasize that it is being considered as the
parameter space for stunted sawtooth maps or for cubic maps.

In discussing the parameter triangle $P^2$, we will refer to the vertex
$\p=(1,1)$, corresponding to a map of entropy $\log 3$, as the {\bit top
vertex\/}. The opposite edge, corresponding to monotone maps, with
entropy zero, will be called the {\bit bottom edge\/} of $P^2$.

In either the stunted sawtooth or the cubic family of shape $\boldsig$,
it is not difficult to check that a bone $B_{\pm}(\O)$ is non-vacuous
if and only if its order type $\O$ is bimodal of shape $\boldsig$.
We will concentrate on order types which are {\it strictly\/}
$\boldsig$-bimodal. In the three exceptional cases where $\O$ is not
strictly bimodal of this shape
(that is for period  two with shape $(-+-)$ or
for period one with either shape), the bones behave rather differently.
For example, in these cases only, the
corresponding bones intersect the bottom edge of the triangle $P^2$
(see Figure 14). These exceptional cases will play only a minor
role in our argument. (Compare \S8.)

\midinsert
\centerline{\psfig{figure=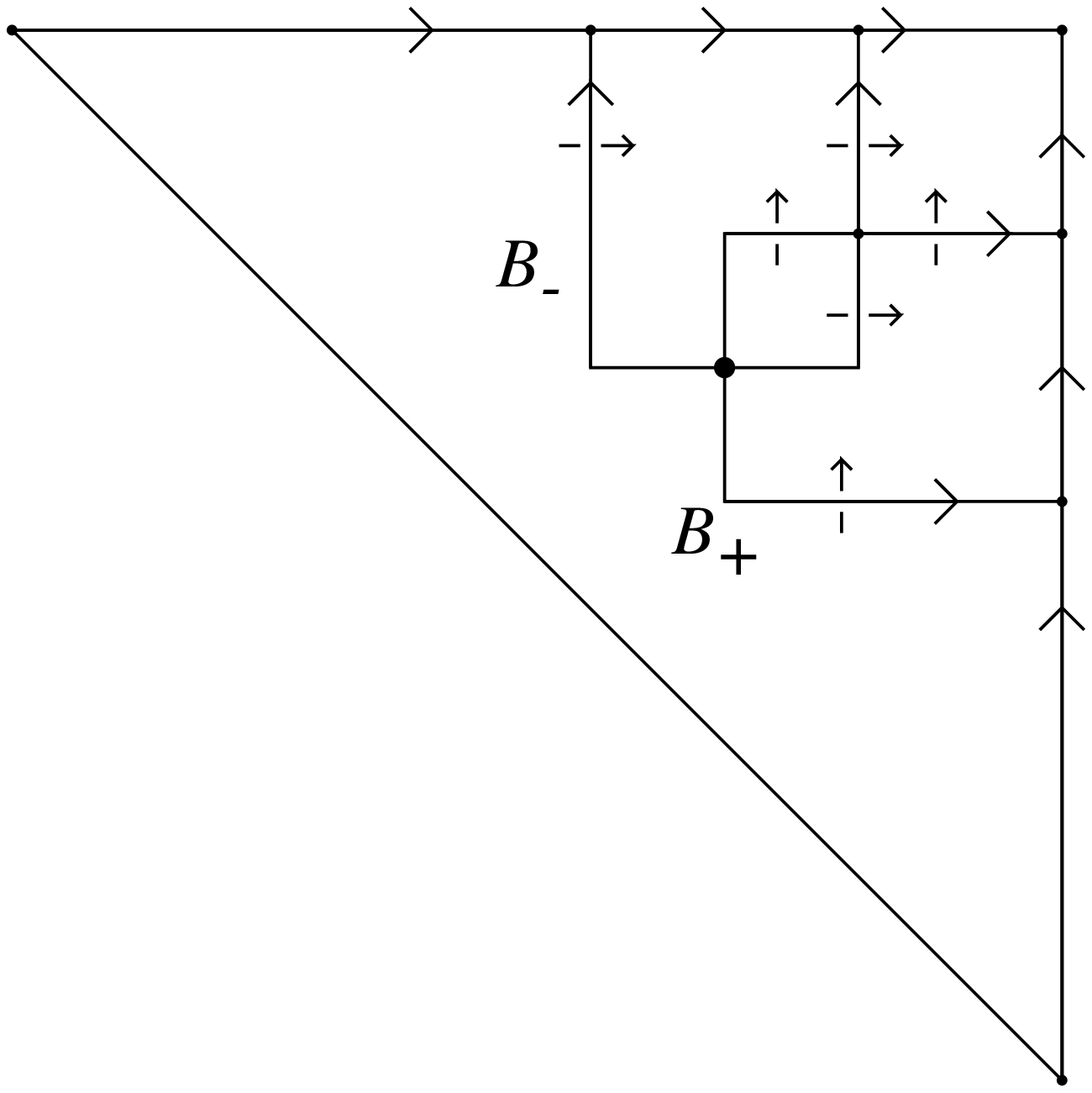,height=2in}\qquad
        \psfig{figure=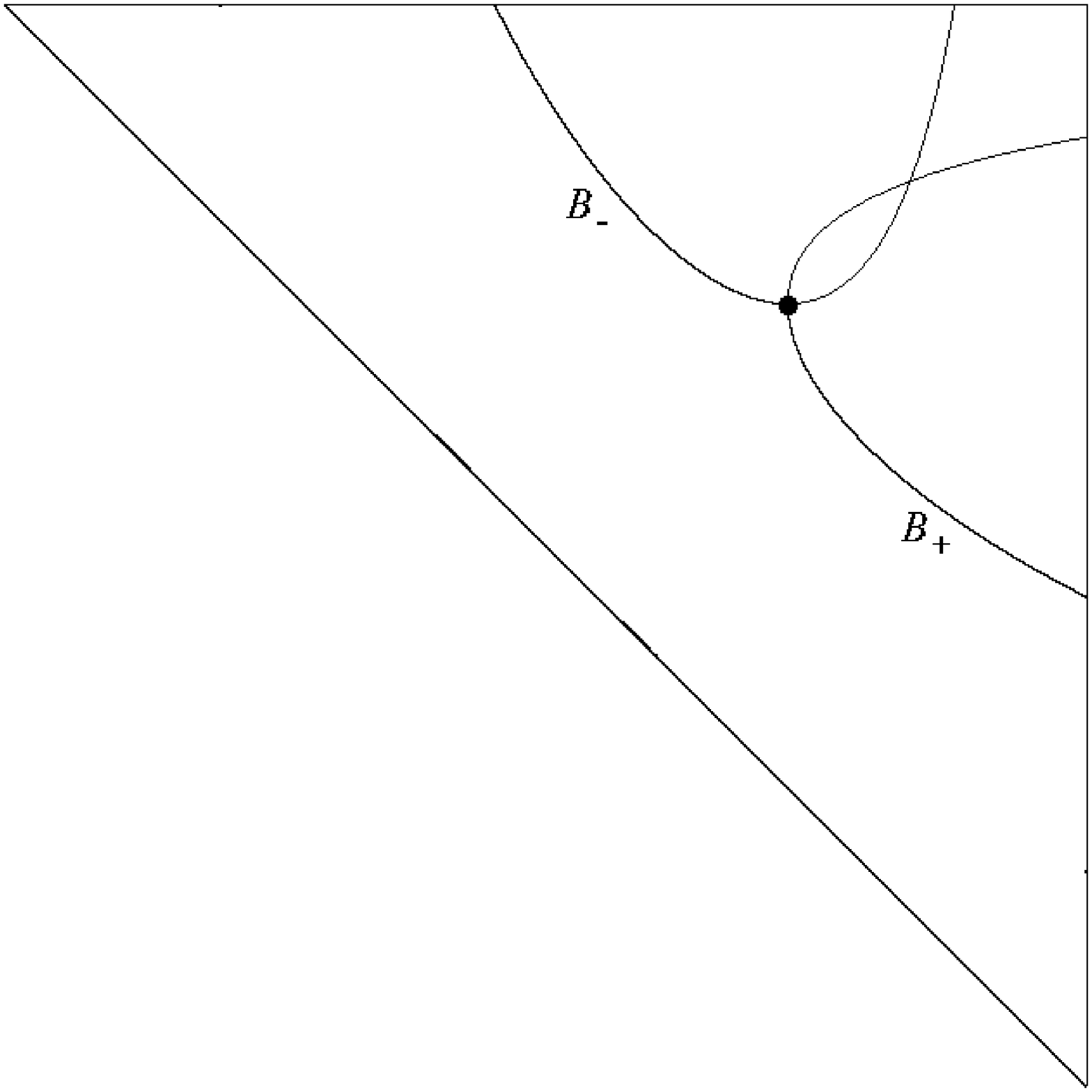,height=1.791in}}
\smallskip
{\wideQP\bit Figure 11. The period 2 bones for the stunted sawtooth family
of shape\break $(+-+)$ on the left, and for the corresponding
cubic family on the right,
with the primary intersection points emphasized.
Arrows point in the direction of increasing dynamic complexity.
(The pictures for either shape $\boldsig$ and for any strictly
bimodal order type of this shape would look qualitatively the same.)
\medskip}
\endinsert

Suppose then that $\O$ is strictly bimodal of shape $\boldsig$. We
consider first the stunted sawtooth case, which is much easier to analyze.

{\QP{\bf Lemma 7.1.} \it For either bimodal shape $\boldsig$ and for each
strictly $\boldsig$-bimodal order type $\O$:

{\rm (i)} Each bone $B_\pm^{\text{saw}}(\O)$
is a simple arc with both endpoints on a common vertical or horizontal
edge of the triangle $P^2=
P^{\,\roman{ saw}}$, and is made up out of three straight
line segments which are alternately horizontal and vertical.

{\rm(ii)} The dual bones $B_-^{\roman{saw}}(\O)$ and
$B_+^{\roman{saw}}(\O)$
intersect transversally
in exactly two points. The intersection at the midpoints of
their middle line segments
will be called their {\bit primary intersection point\/}. It corresponds
to the unique map $S_\p$ for which both critical points lie on a common
periodic orbit of order type $\O$. There is also a {\bit secondary
intersection point}, higher in $P^{\,\roman{saw}}$,
corresponding to the unique map $S_\p$ which has two
disjoint critical orbits, both periodic with the same order type $\O$.

{\rm(iii)} This pair of dual bones cuts the triangle $P^{\,\roman{saw}}$
into five regions. There is a large region which contains the entire
bottom edge of $P^{\,\roman{saw}}$ and contains the primary intersection
point in its boundary.
Adjoining this are two side regions, each touching just one edge of
$P^{\,\roman{saw}}$.
Next there is a central region, which is disjoint from all
other bones so that the topological entropy takes a constant value
(equal to the logarithm of an algebraic unit) throughout.
Finally, there is a top region whose boundary
contains both the secondary intersection point
and the top vertex of $P^{\,\roman{saw}}$. % in its boundary.

{\rm(iv)}
For $\O\ne \O'$, the bones $B_-^{\roman{saw}}(\O)$ and $B_+^{\text saw}
(\O')$ intersect trans\-versally in either 0, 2, or 4 points. Each of these
intersection points corresponds to a map with two disjoint periodic
critical orbits. Again the entropy is the logarithm of an algebraic unit.
\bigskip}

{\bf Proof of (i), (ii) and (iv).}
Since the order type $\O$ is strictly bimodal of shape
$\boldsig$, it follows easily from 5.3 that there exists a point
$\p^0=(p_1^0\,,\,p_2^0)\in P^{\,\roman{saw}}$ so that
both critical points of the map
$S_{\p^0}$ lie in a common period $q$ orbit with order type $\O$.
(The argument will show that this point $\p^0$ is unique. Compare 8.1
below.)
Evidently there is a unique $0<n<q$ so that the iterate
% $$ S_{\p^0}^{\circ n}(\widehat C_1)\=\widehat C_2\,,\quad
%    S_{\p^0}^{\circ (q-n)}(\widehat C_2)\=\widehat C_1~. $$
$S_{\p^0}^{\circ n}$ maps the critical point
$\widehat C_1$ to $\widehat C_2$ while $S_{\p^0}^{\circ (q-n)}$ maps
$\widehat C_2$ to $\widehat C_1$.
The {\bit kneading type\/} $\I^0$ of this periodic orbit, that is the
itinerary of a
representative point (well defined up to a shift), is periodic of period
$q$, and contains each of the
critical point symbols $C_1$ and $C_2$ exactly once in each period. For
an arbitrary periodic orbit of the same order type, as noted above,
there are $3^m=9$ different possibilities for the kneading type:
If we replace $(p_1^0\,,\,p_1^0)$ by some nearby point $(p_1\,,\,p_2)$,
then there will still be a unique periodic orbit which intersects both
plateaus. However, depending on the sign of $p_2-p_2^0$ and of $p_1-p_1^0$,
we can replace the symbol $C_1$ in the kneading type $\I^0$ by either of
the two
adjacent symbols $I_0$ or $I_1$, and we can replace $C_2$ by either
$I_1$ or $I_2$. However, no other replacements are possible. In particular,
for any orbit of order type $\O$ for a bimodal map of shape $\boldsig$,
the remaining $q-2$ symbols must remain unchanged.

To fix our ideas, consider a left bone $B_-(\O)$.
For $\p\in B_-(\O)$ the left critical point $\widehat
C_1$ is periodic. Evidently each plateau
of a stunted sawtooth map can contain at most one periodic point, and for
$\p\in B_-(\O)$ it is easy to check that
only the point $S_{\p}^{\circ n}(\widehat C_1)$
in the orbit of $\widehat C_1$ can lie in the right hand plateau.
We divide the discussion into three cases,
corresponding to the three line segments which make up the bone $B_-(\O)$.
% according as this orbit does or does not lie in the right plateau.

{\bf Case 1.} If $S_{\p}^{\circ n}(\widehat C_1)=S_\p^{\circ(n-1)}(v_1)$
does lie in the right hand plateau, then the iterate $S_\p^{\circ(n-q-1)}$
must
map the critical value $v_2=S_\p(\widehat C_2)$ to $\widehat C_1$. Since
the
partial itinerary from $v_2$ to $\widehat C_1$ is uniquely determined by
the
order type, this yields a linear equation which we can solve for $v_2$,
and hence for $p_2=p_2^0$.

{\bf Cases 2, 3.}
On the other hand, if $S_{\p}^{\circ n}(\widehat C_1)$ lies either to the
left or to the right of the right hand plateau, then the entire
partial itinerary from $v_1$ back to $\widehat C_1$ is uniquely determined.
In either of these two cases,
this yields a linear equation which we can solve for $v_1$.

Alternatively, we can describe these three cases as follows. Suppose that
we
fix $p_2=p_2^0$ and vary $p_1$ throughout a neighborhood of $p_1^0$.
Then the size of the right hand plateau remains fixed, and
there will be a largest
interval $[p_1^-\,,\,p_1^+]$, symmetric about $p_1^0$, so that,
for $p_1$ in this interval, the image $S_\p^{\circ(n-1)}(v_1)$ lies in the
right hand plateau. For the extreme values $p_1=p_1^\pm$, this image will
lie at either end of the plateau. Now fix $p_1=p_1^\pm$ and let $p_2$
vary over the interval $[p_2^0\,,\,1]$. As $p_2$ increases, the right
hand plateau will move up or down (depending on $\boldsig$), shrinking in
length, but the orbit of $\widehat C_1$ will remain fixed, missing the
right
hand plateau completely when $p_2>p_2^0$. It follows that
$B^-(\O)$ is the union of three line segments:
$$   B_-(\O)~\= ~p_1^-\times[p_2^0\,,\,1]~~\cup~~[p_1^-\,,\,p_1^+]\times
      p_2^0    ~~\cup~~ p_1^+\times[p_2^0\,,\,1]~. $$
The discussion of right bones in completely analogous. This proves Part
(i);
and Parts (ii), (iv) follow easily.

{\bf Proof of (iii).} It is easy to see that  a pair of dual bones
partitions
$P^{\,\roman{saw}}$ into five regions.
To prove that the central region is disjoint from all
other bones, let us follow the left bone $B_-(\O)$ from the primary
intersection point $\p^0$ to the secondary intersection point. For each
point
$\p$ in this path, the image $C'(\p)=S_\p^{\circ n}(\widehat C_1)$ is
periodic, belonging to the orbit of $\widehat C_1$. Let $J_\p$ be the
interval with endpoints $\widehat C_2$ and $C'(\p)$. For $\p=\p^0$ this
interval degenerates to a point, and for $\p$ in a large neighborhood of
$\p^0$ the iterate $S_\p^{\circ q}$ maps all of $J_\p$ to the endpoint
$C'(\p)$. However, as $\p$ moves towards the secondary intersection point,
eventually its image under $S_\p^{\circ q}$ will become bigger. However,
the
restriction of $S_p^{\circ q}$ to $J_\p$ will remain monotone, and its
image
will remain a proper subset of $J_\p$ until $\p$ reaches the secondary
intersection point, at which time $S_\p^{\circ q}(\widehat C_2)$ will
equal $\widehat C_2$, so that $J_\p$ maps onto itself with both endpoints
fixed. Until $\p$ reaches this point, clearly $\widehat C_2$ will not be
periodic, so no other bone can cross through to the central region.

The statement that topological entropy is constant in any region without
bones follows easily from 4.11. (Compare 9.1 below.) Such a region
necessarily contains hyperbolic points by 5.9. Hence this constant value
must be the logarithm of an algebraic integer by 4.3.
This proves (iii), and completes the proof of 7.1.\QED

The statement corresponding to Lemma 7.1 for the cubic family
is much more difficult. Here is some preliminary information.
Again let $\O$ be strictly bimodal of shape $\boldsig$.

{\QP{\bf Lemma 7.2.} \it Both bones $B_\pm^{\text {cub}}(\O)$ are smooth
1-dimensional manifolds with boundary.
% equal to $B_\pm^{\text {cub}}(\O)\cap\partial P^{\,\roman{cub}}$.
Further, the boundary of $B_-^{\text {cub}}(\O)$ is equal to the
transversal intersection of this bone with the horizontal (top) edge
of the triangle $P^{\,\roman{cub}}$, while the boundary of $B_+^{\text
{cub}}
(\O)$ is equal to the transversal intersection of this bone with the
vertical
(right hand) edge
of $P^{\,\roman{cub}}$. Any intersection $B_-^{\text {cub}}(\O)\cap
 B_+^{\text {cub}}(\O')$ between two such bones is also transverse.\ss}

(If we make use of Heckman's Theorem that there are no bone-loops, then
we have the much sharper statement that each bone is a connected simple
arc.)

\comment ((((((((((((((((((
, and these boundary points belong to a common vertical
or horizontal edge of the
triangle $P^2$. Furthermore, if $\O'$ is also a strictly bimodal order
type,
then the intersection between $B_-^\roman{cub}(\O)$ and
$B_+^\roman{cub}(\O')$ is transverse, and the number of intersection
points is precisely equal to the number of intersections between
the corresponding bones
$B_-^\roman{saw}(\O)$ and $B_+^\roman{saw}(\O')$. In particular,
in the case of dual bones with $\O=\O'$, there are exactly two intersection
points.\smallskip}
\endcomment %%%%%%%%%%%%%%%%

The proof of 7.2 begins as follows. First consider the corresponding
statement for the family
of complex maps $z\mapsto z^3-3a^2z+b$, with critical points $\pm a$.
It is proved in [M3] that the locus ${\Cal S}_\pm(p)$ of points for which
$\pm a$ has period $p$ is a smooth complex curve. Furthermore, for each
$p$ and $q$ the curves ${\Cal S}_+(p)$ and ${\Cal S}_-(q)$ intersect
transversally. (See [M2].)
In fact, ${\Cal S}_+(p)$ has transverse intersection with
any curve consisting
of points for which the other critical point $-a$ is
preperiodic. The proofs make essential
use of quasi-conformal surgery. Compare [St], where analogous
results for quadratic rational maps are proved by similar methods.
Alternatively, these statements have been proved by quite different
methods by Epstein [E].

Restricting to the real $(a,b)\!$-plane, we obtain a corresponding
statement for real cubic maps of shape $(+-+)$.
{\it In the family of real maps\vskip -.2in
$$x~\mapsto~ x^3-3a^2x+b~,\eqno (17^+)$$\vskip -.1in
\noindent
the locus of pairs $(a,b)$ for which $a$ (or $-a$) is periodic of period
$p$ forms a smooth $1\!$-dimensional manifold
without boundary.} A complex linear change of coordinates yields the
corresponding statement for the family of maps
$$   x~\mapsto~-(x^3-3a^2x+b)~, \eqno (17^-) $$
of shape $(-+-)$.

In order to relate this to the family of cubic maps $f_\p$,
we note that each $f_\p$ is {\bit positively
affinely conjugate\/} to a unique
map in the normal form $(17^\pm)$ with
$a\ge 0$. That is, for each $\p$ there
is a unique affine map
$L(x)=c\,x+d$ with $c>0$ so that $L\circ f_\p\circ L^{-1}$ has the
required form. It follows that there is a well defined continuous mapping
$\phi:P^{\,\roman{cub}}\to\R^2$ which associates to each $\p\in
P^{\,\roman{cub}}$
the associated pair $\phi(\p)=(a,b)$ in the half-plane $a\ge 0$. Evidently
the pre-image of the curve ${\Cal S}_\pm(p)$ under $\phi$
is the union of all bones
$B_\pm^{\text {cub}}(\O)$ of period $p$. If $\phi$ were a diffeomorphism,
then it would follow immediately that each bone is a smooth 1-manifold,
with
boundary points precisely at the intersections with the
horizontal or vertical part of the boundary of $P^{\,\roman{cub}}$.
In fact, the
situation is somewhat more complicated, and can be described as follows.
(Compare [DGMT, Figures 8,9.)

{\QP{\bf Lemma 7.3.} \it Let $P^\roman{cub}$ be the parameter triangle for
cubic maps of specified shape $\boldsig$, and let $U\subset P^\roman{cub}$
be the open subset consisting of those parameter values
$\p$ such that the boundary $\partial I$
is strictly repelling for the associated cubic map $f_\p:I\to I$. Then
every bone of strictly bimodal order type in $P^\roman{cub}$ is
contained in this open set $U$. Furthermore, the
map $\phi:P^{\,\roman{cub}}\to\R^2$ embeds $U$ diffeomorphically into
$\R^2$.\ss}

Clearly Lemma 7.2 will follow from 7.3, together with the discussion above.

{\bf Proof of 7.3.} If $\p$ belongs to the complement $P^\roman{cub}\ssm
U$,
then some boundary orbit of period one or two for $f_\p$ has multiplier
in the interval $[0,1]$. Since $f_\p$ has negative Schwarzian derivative,
this implies that there is at least one critical point in the immediate
basin,
and this clearly implies
that the interval between the two critical points is also
contained in this immediate basin. Thus there cannot be any periodic
critical
orbit, unless one and hence both critical points are actually contained
in the boundary
$\partial I$. This occurs only if $\p$ is one of the two lower
corner points of the parameter triangle. Evidently the corresponding
order type $\O$ is not strictly bimodal. (These zero entropy corner points
are
quite anomalous, since each one
represents an isolated point of the associated bone of
period one or two.)

For any $\p\in U$ the boundary periodic orbit (in the $(-+-)$ case) or
orbits (in the $(+-+)$ case) are strictly repelling, and hence vary
smoothly
under smooth deformation of the polynomial. It follows that there is a
smooth local inverse function from the image $\phi(U)\subset\R^2$ back to
$U$.
This inverse function is single valued, since otherwise we would have
two different intervals, say $I$ and $I'$, both containing the two critical
points, so that the same polynomial function $f_\p$ maps each of these
intervals into itself with strictly repelling boundary. Then any component
of
$I\ssm I'$ would map diffeomorphically
into itself under $f_\p$ or $f_\p\circ f_\p$, with
both boundary points repelling, which is impossible for a map of negative
Schwarzian. This completes the proof of 7.3 and 7.2.\QED

\comment
\midinsert
\centerline{\psfig{figure=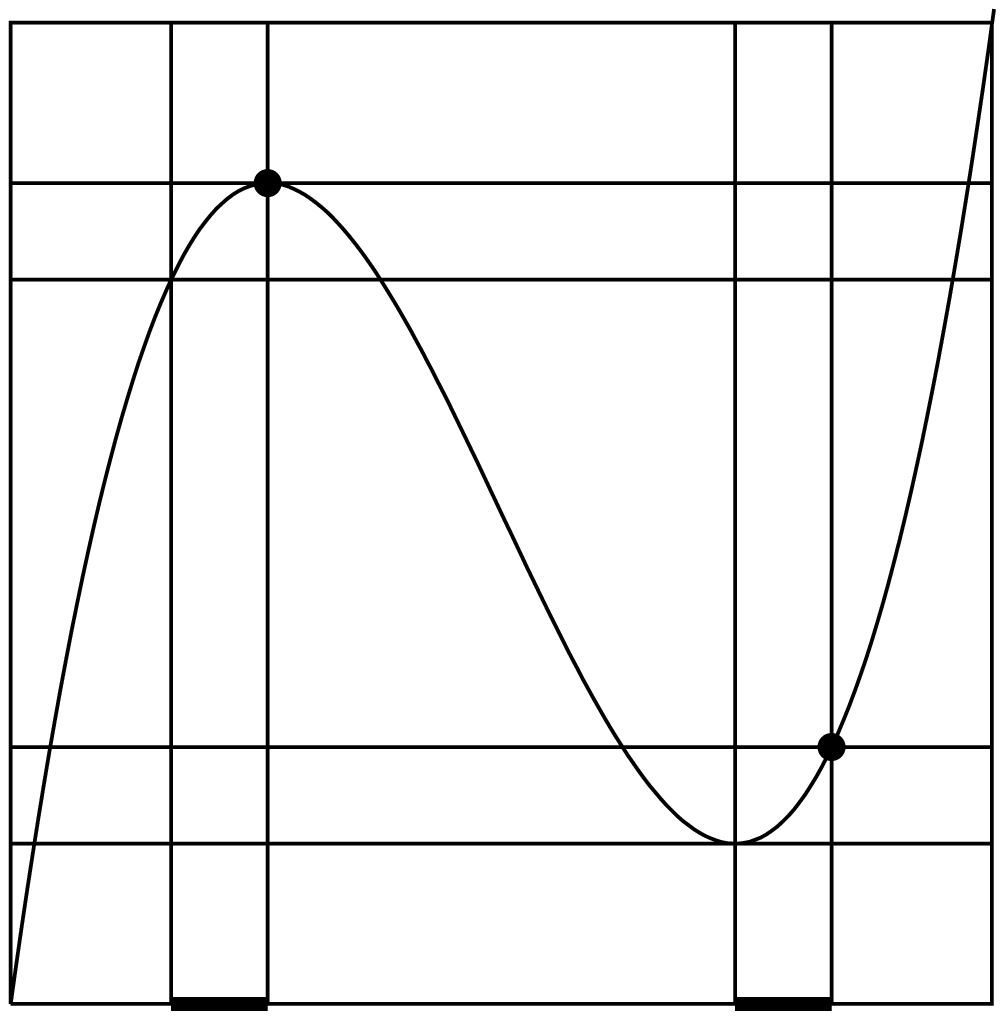,height=2.4in}}\smallskip
{\QP\bit Figure 12. Graph of the map representing the secondary
intersection point of the period two bones in the family of
cubic maps of shape $(+-+)$. One of the two critical orbits has been
emphasized. Note the emphasized intervals, which join the two periodic
orbits and are mapped homeomorphically by all iterates of the map.
\medskip}
\endinsert
\endcomment

\midinsert
\centerline{
\psfig{figure=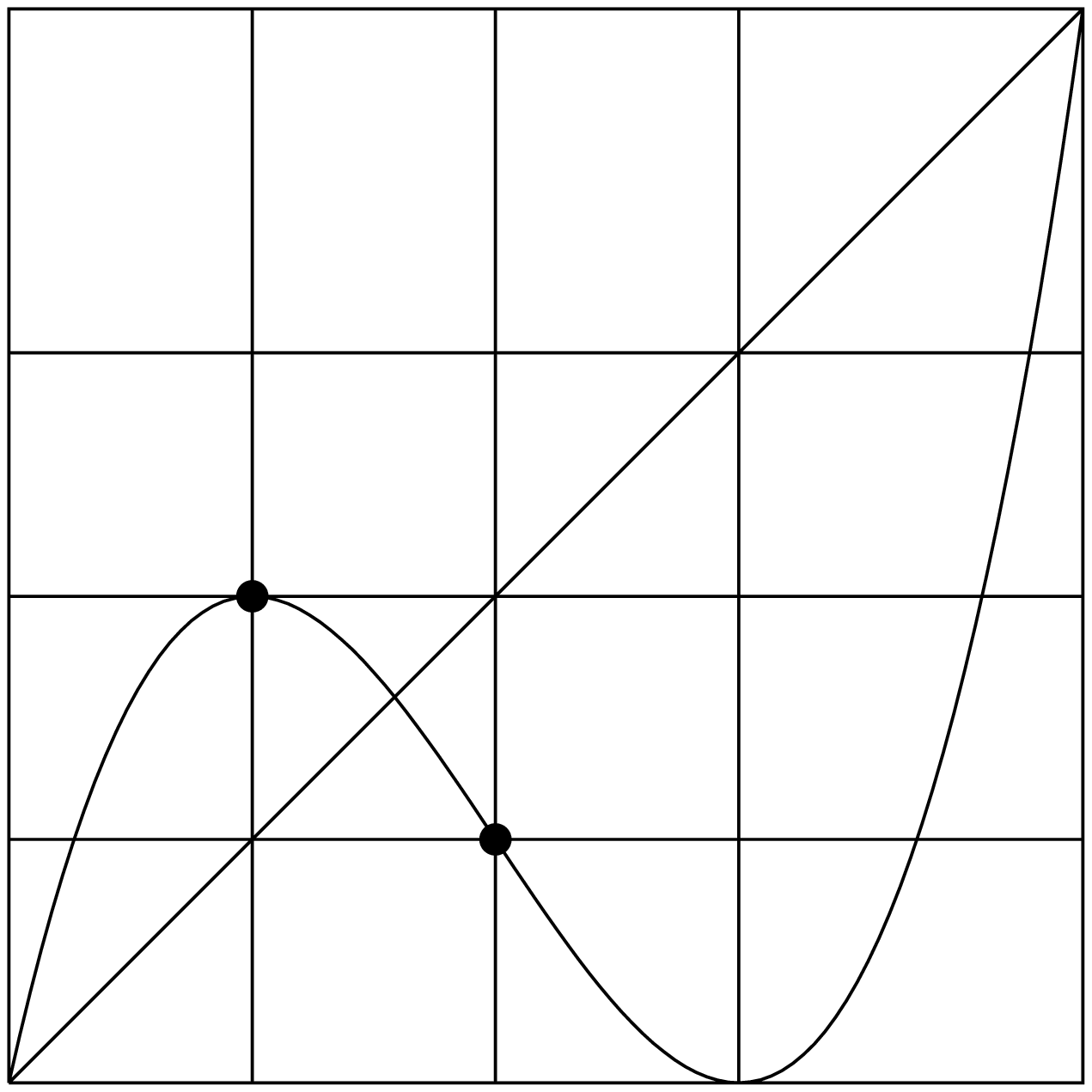,height=1.2in}\qquad
\psfig{figure=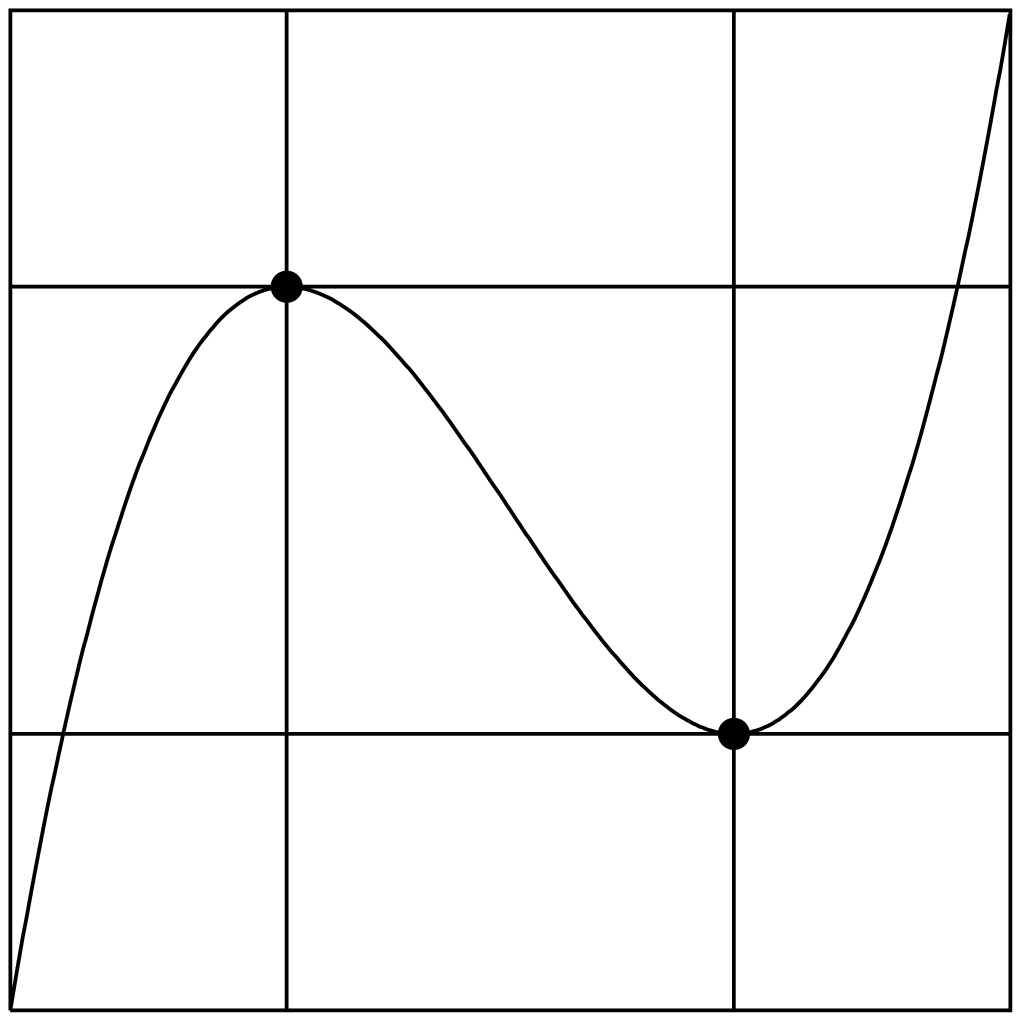,height=1.2in}\qquad
\psfig{figure=per2X2.ps,height=1.2in}\qquad
\psfig{figure=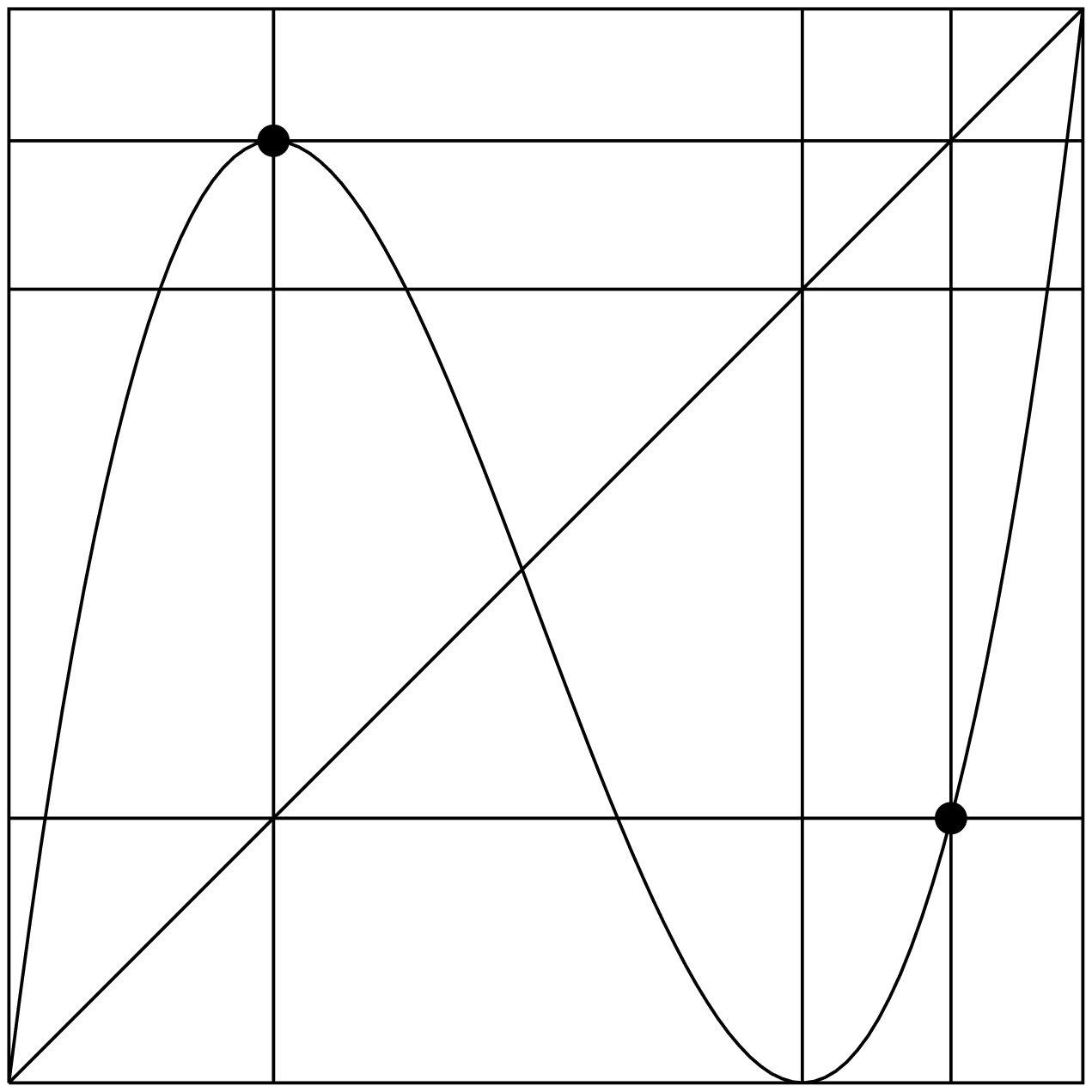,height=1.2in}}\ss
{\QP\bit Figure 12. Graphs of cubic maps of shape $(+-+)$
representing successive points along the left bone of period two.
Illustrated
are the ``lower'' endpoint of this bone, the primary intersection point
where
dynamic complexity is minimized,
the secondary intersection point, and the ``upper'' endpoint. In the third
graph there are three different period two orbits, all contained in the
emphasized intervals, and these period two orbits all persist to the right
hand graph, which has the largest negative orbit complexity.\ss}
\endinsert

{\bf Remark 7.4.} There is a fundamental relationship between bones and
negative orbit complexity. (Compare 9.1.) In fact, as we traverse a path
in the parameter triangle, the number of fixed points of negative type
for the $n$-th iterate can clearly change only as we cross a bone of period
dividing $n$. More precisely, as we cross a period $p$
bone in the `positive' direction,
as indicated in Figure 11, it turns out that exactly one
new periodic orbit of negative type is created, with period equal to $p$.
This fact can be used to count the number of period $p$ bones. In the
$(+-+)$-case, note that the number of fixed points of negative type for
the $n$-th iterate increases from zero, along the bottom edge of the
parameter triangle, to a maximum of $(3^n-1)/2$ for the map, say $S_{(1,1)}$,
corresponding to the top (high entropy) vertex of the triangle. {\it It follows
easily that $\#\roman{Bones}(p)$, the number of period $p$ bones in the
parameter triangle $($suitably interpreted when $p=1)$,
is precisely equal to the number
of period $p$ orbits of negative type for the map $S_{(1,1)}$.\/}
From this, we derive the following formula, which can be used to compute
$\#\roman{Bones}(p)$.
$$	(3^n-1)/2~=~\sum_{p|n~\text{with}~ n/p~\text{odd}} p\cdot
	\#\roman{Bones}(p) $$
For example:
$$      \matrix p=~~ & 1 & 2 & 3 & 4 & 5 & 6 & 7 & 8& 9\cr
        (3^p-1)/2=~~ & 1 & 4 & 13 & 40 & 121 & 364 & 1093 & 3280 & 9841~\cr
	\#\roman{Bones}(p)= & 1 & 2 & 4 & 10 & 24 & 60 & 156 & 410 & 1092.
  \endmatrix $$
Note that $\#\roman{Bones}(p)$ is even for $p>1$, since the bones occur in
dual pairs. The discussion for the $(-+-)$-case is slightly different, but
the number $\#\roman{Bones}(p)$ remains the same.

{\bf Remark 7.5.} Evidently the intersection points of bones are precisely
those points in parameter space for which {\it both\/} critical
points are periodic. We will see in \S8 that there is a one-to-one
correspondence, preserving kneading data, between these intersection
points for the stunted sawtooth case and for the cubic case. Alternatively,
this statement would follow easily from a form of
Thurston's Theorem, as stated in B.3 and B.4 of Appendix B.
In general, for the cubic map corresponding to any intersection of bones,
the two critical points will belong to
disjoint periodic orbits. The only exception is for the primary
intersection point between two dual bones. In this exceptional case, the
two critical points belong to a common orbit. (See Figure 11 for the
picture
in $P^2$, and see Figure 10 for graphs of some cubic maps corresponding
to primary intersections between dual bones.)

Similarly, the endpoints of strictly $\boldsig$-bimodal bones correspond to
 maps for
which one critical point is periodic while the other is preperiodic,
mapping directly to a boundary fixed point or period two orbit.
Again there is a canonical
one-to-one correspondence between such endpoints for the stunted sawtooth
and cubic families. In particular, it follows that
each cubic bone has exactly two endpoints. (See Figure 12 for
graphs of several maps corresponding to parameter values along
a left bone, including the two endpoints.)

\comment %%%%%%%%%%%%%%%%%%%%%
In order to analyze these possible boundary points
$B_\pm^{\text {cub}}\cap \partial P^2$, we will need to make use of
the following basic principle. By definition, a piecewise monotone
map is {\bit post-critically finite\/}
if the orbit of every critical point is periodic or eventually periodic.

{\QP{\bf Thurston's Theorem for Real Polynomial Maps.} {\it Given any
post-critically finite $m\!$-modal map, there exists one and
up to positive affine conjugation only one polynomial of degree $m+1$
with the same kneading data.
\footnotemark\footnotetext
{Although this statement is well known to experts, it is difficult
to find in the literature. The proof makes essential use of
complex methods. In fact this statement is an easy corollary of a
much more complicated statement for complex polynomials (or for
complex rational maps). For further information, the reader is referred
to [DH1], [Po1] and [MvS2], as well as the discussion in [MTh].}
}\smallskip}

Assuming Thurston's Theorem, the proof of Lemma 7.3 concludes
as follows.

Consider a map $f_\p$ belonging to the intersection
$B_-^{\text {cub}}(\O)\cap
\partial P^2$. Thus the left critical point of $f_\p$
is periodic, with some period $p\ge 2$. It is easy to
check that the two critical points cannot coincide. (Here we
make use of the fact that our maps have shape $+-+$.
For maps of shape $-+-$ we would rather have to assume
that $p\ge 3$ in order to avoid the case $c_1=c_2$.) Since
$f_\p\in\partial P^2$, it follows
that either $p_1=1$ or $p_2=1$. In other words, at least
one of the two critical points must map to a boundary fixed point.
But $c_1$ is periodic, so it follows that we must be on the edge $p_2=1$.
If we cut the interval $I$ at the points of the period $p$ orbit
$\{c_1,f_\p(c_1),f_\p^\circ(c_1),\ldots\}$, then since our map
must have shape $+-+$ a little work shows that there are exactly
two of the $p+1$ complementary intervals
where the right hand critical point can be placed. Using
Thurston's Theorem, it follows that the intersection
$B^{\text {cub}}_-(\O)\cap\partial P^2$
consists of exactly two points. Further, this intersection is
transverse, so these two intersection points must belong to the
boundary $\partial B^{\text {cub}}_-(\O)$.
\QED
\endcomment %%%%%%%%%%%%%%%%

%From tresser@us.ibm.com Tue Feb 23 10:01:00 1999
%Date: Mon, 22 Feb 1999 19:03:32 -0500
%From: tresser@us.ibm.com

%\documentstyle{amsppt}
%\magnification=1200
%\vsize=9 true in
%\hsize=6.5 true in
%\voffset=0pt
%\hoffset=0pt
%\overfullrule=0pt
%\endcomment %%%%%%%%%%%%%
%%%%
%\parskip=4pt
%\mathsurround=2pt
%\input amssym.def
%\input amssym
%%\pageno=39
%%\adjustfootnotemark 5
%\input psfig
%\def\boldsig{\sigma\!\!\!\!\sigma}
%\def\[{$\,}
%\def\]{\,$}
%\def\Neg{\roman{Neg}}
%\def\cl{\centerline}
%\def\min{\roman{min}}
%\def\max{\roman{max}}
%\def\C{{\Bbb {C}}}
%\def\N{{\Bbb {N}}}
%\def\Q{{\Bbb {Q}}}
%\def\R{{\Bbb {R}}}
%\def\Z{{\Bbb {Z}}}
%\def\A{{\frak {A}}}
%\def\cf{{\frak {C}}}
%\def\I{{\Cal {I}}}
%\def\K{{\Cal {K}}}
%\def\O{{\bold o}}
%\def\p{{\bold p}}
%\def\q{{\bold q}}
%\def\v{{\bold v}}
%\def\w{{\bold w}}
%\def\>{\gg} %{>\!\!>}
%\def\<{\ll} %{<\!\!<}
%\def\={~=~}
%\def\sgn{\roman{ sgn}}
%\def\QED{  \rlap{$\sqcup$}$\sqcap$ \smallskip}
%\def\QP{\smallskip\leftskip=.4in\rightskip=.4in\noindent}
%\def\wideQP{\smallskip\leftskip=.25in\rightskip=.25in\noindent}
%\def\b1{{\pmb{1}}}
%\def\ss{\smallskip}
%\font\bit=cmssi12 at 12truept
%\font\tenmsy=msym10
%\font\sevenmsy=msym7
%\font\sma=cmr10 at 10 truept
%\textfont8=\tenmsy
%\def\ssm{\smallsetminus}
%\mathsurround = 1pt
%\abovedisplayskip=4pt
%\belowdisplayskip=4pt
%\parskip=5pt
%\def\ref{\hangindent=1pc \hangafter=1 \noindent}
%\def\s{\thinspace}
%\def\vs{\vskip .05in}
%
%\centerline{\bf On Entropy and Monotonicity\qquad $\cdots$
% \qquad[ \S8-end ]}\ss

%\centerline{6-16-98 last hand = jm}\ss

\bigskip

\bigskip

\maybebreakhere
\centerline{\bf\S8. Monotonicity, Intersections of Bones, and
the $n\!$-Skeleton.}
\medskip
\bigskip

In the quadratic case, it is
known that the number of period $q$ points for the interval map
$Q_v(x)=4vx(1-x)$ increases monotonically as the folding value
parameter $v\in [0,1]$ increases. Similarly the kneading sequence
$\K_1(Q_v)$ increases monotonically with $v$. Proofs of this result have
been given by
Sullivan,\footnote
{As far as we know, Sullivan did not actually publish a proof.
However he communicated one orally at an early date, and
conversations with Sullivan led to the Milnor-Thurston proof.}
Douady and Hubbard, and by Milnor and Thurston.
Compare [DH2, n$^{\!\text{o}}$VI],
[MTh], [D], % [McS ??],
as well as [MvS2].

We can adapt similar techniques to
study special one parameter families of cubic maps. (Compare [NN].) %,
[KN].)
However, instead of working with the full kneading data, it will suffice
for our purposes to work with the ``negative orbit complexity''
$${\Cal N}(f)
=\left(\Neg(f^{\circ 1})\,,\,\Neg(f^{\circ 2})\,,\,\ldots\right)$$
of 4.13. In fact, we will establish monotonicity of ${\Cal N}(f)$ in a
collection of one parameter families
large enough to give us a good picture of the full two parameter family,
making use of Heckman's Theorem that there are no bone-loops.
% help us capture essentially all of the full picture
%if we assume that there are no bone-loops, or
%if we assume that the Generic Hyperbolicity Conjecture holds true.

By definition, a polynomial map is {\bit post-critically finite\/}
if the orbit of every critical point is periodic or eventually periodic.
Similarly, we will say that a
piecewise monotone map has {\bit finite folding point
orbits\/} if the orbit of every folding point is periodic or eventually
periodic. (If $f:I\to I$ is piecewise strictly monotone and boundary
anchored, then a completely equivalent condition is that $f$ be a ``Markov
map'' in the sense that there exist a subdivision of $I$ into finitely many
subintervals, each of which maps homeomorphically onto some union of these
subintervals.)

%Recall that a polynomial map is said to be
%{\bit post-critically finite\/}
%if the orbit of every critical point is periodic or eventually periodic.
%Our proofs of monotonicity will be based on the following basic principle.

{\QP{\bf Thurston Uniqueness Theorem for Real Polynomial Maps.} \it
A post-critically finite real polynomial map of degree $m+1$ with
$m$ distinct real critical points is uniquely determined, up to positive
affine conjugation, by its kneading data.\smallskip}

See Appendix B, which augments this to a full existence and uniqueness
statement, and gives some indication of the proof. There is an
analogous statement for the stunted sawtooth family, as follows. As usual,
we assume that the slope $s>m+1$ has been fixed.

{\QP{\bf Lemma 8.1.} \it Given any
piecewise strictly monotone $m\!$-modal map such that each folding point
has finite forward orbit, there exists one and
only one stunted sawtooth map of the same shape
with the same kneading data.\medskip}

\midinsert
\cl{\psfig{figure=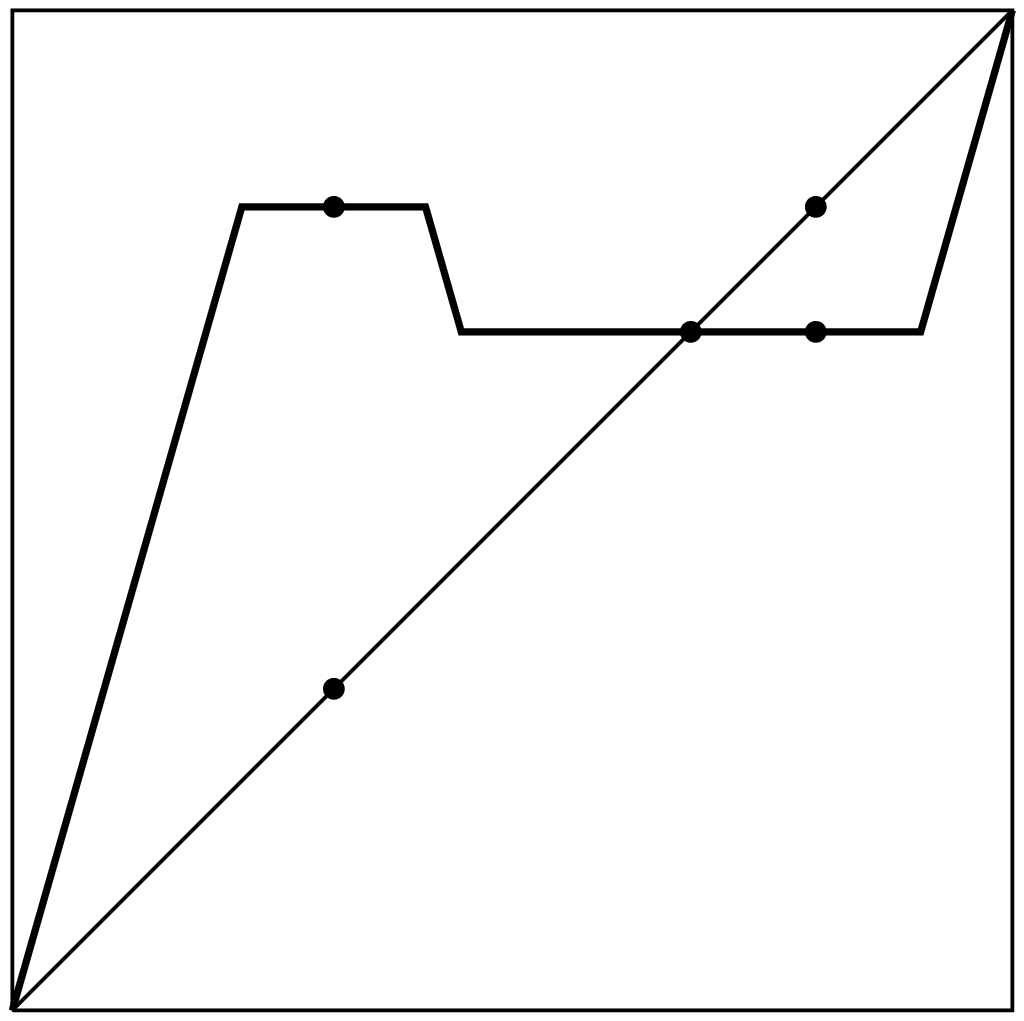,height=1.3in}}
\endinsert

{\bf Caution.} Here it is essential that the kneading
data comes from a piecewise {\bit strictly\/} monotone map. The following
figure shows an example to illustrate this.
Recall that the ``folding points'' of a stunted sawtooth map
are defined to be the midpoints $\widehat C_i$ of the plateaus.
Let $S_\p$ be a stunted sawtooth map as illustrated, of shape $(+-+)$ with
folding values $v_1=\widehat C_2+\epsilon>v_2=\widehat C_2$. If
$\epsilon>0$
is sufficiently small, then this map is post-critically finite,
with kneading data
$    \K_1=(I_2\,,\,\overline C_2)\,,\; \K_2=\overline C_2$,
where the overline stands for infinite repetition.
%$$  \K_1\=(I_2\,,\,C_2\,,\,C_2\,,\,\ldots)\,,\qquad \K_2\=(C_2\,,\,
%C_2\,,\,\ldots)~. $$
Thus there is an entire 1-parameter family of stunted sawtooth maps with
the same post-critically finite kneading data. Yet no polynomial map, and
more generally no piecewise strictly monotone map, can have this kneading
data, which requires that an entire interval in $I_2$ must map to the
second
folding point. (Compare Appendix B.)

In fact there is a simple criterion which guarantees that the kneading
data for some given post-critically finite
map $f$ can be realized by a map which is
piecewise strictly monotone. Let $\{y_1\,,\,\ldots\,,\,y_n\}$ be the union
of the forward orbits of the folding points, where\break
$y_1<\cdots<y_n$. {\it If
$f(y_i)\ne f(y_{i+1})$ for $1\le i<n$, then interpolating linearly between
the
$y_i$, and extending suitably over the rest of the interval $[a,b]$,
we easily obtain a piecewise {\it strictly\/} monotone map with
the same kneading data.} \smallskip

{\bf Proof of 8.1.} Existence follows from 5.3. To prove uniqueness of
the resulting stunted sawtooth map $S_\p$, note
that the orbit of a folding point of $S_\p$ can never
land in a plateau, except at its midpoint $\widehat
C_j$. For otherwise, any map with the same kneading data would have
an interval of points eventually mapping to the same point,
which would contradict the hypothesis. Therefore, orbit points are uniquely
determined by their itineraries, hence $S_\p$ itself
is uniquely determined.\QED\smallskip

%We now proceed to study the monotonicity properties along the one
%parameter families which make up the essential skeletons.

We now proceed to apply these results to the study of monotonicity
along special curves. As a first example, consider the
top edge $p_2=1$ of the parameter triangle for either cubic
or stunted sawtooth maps of fixed shape
$\boldsig$. This edge corresponds to the set of maps $f:I\to I$
(in either family) such that the second critical value $v_2$
is an endpoint of the interval $I$, and hence is either a fixed
point, in the $(+-+)$ case, or a period two point, in the $(-+-)$ case.
Thus the kneading sequence ${\Cal K}_2$ remains constant as we traverse
this
edge, and only $\K_1$ varies.

{\QP{\bf Lemma 8.2} \it In either the cubic or the stunted sawtooth family
of specified shape $\boldsig$,
as we traverse the top edge $p_2=1$ of the parameter triangle in the
direction of increasing $p_1$, the negative orbit complexity $\Cal N$
increases
monotonically. More precisely, $\Cal N$ increases
whenever we cross the endpoint of a left bone, and remains
constant otherwise. Furthermore these crossing points occur in the same
order
in the two families: To each crossing point in one family,
there corresponds one and only one crossing point in the other family with
the same kneading data; where
%     between crossing points which preserves the kneading data,
the crossing point $f$ (in either family) occurs before $f'$ (in the same
family) if and only if ${\bold{K}}(f)\ll{\bold{K}}(f')$.\medskip}

\midinsert
\centerline{\psfig{figure=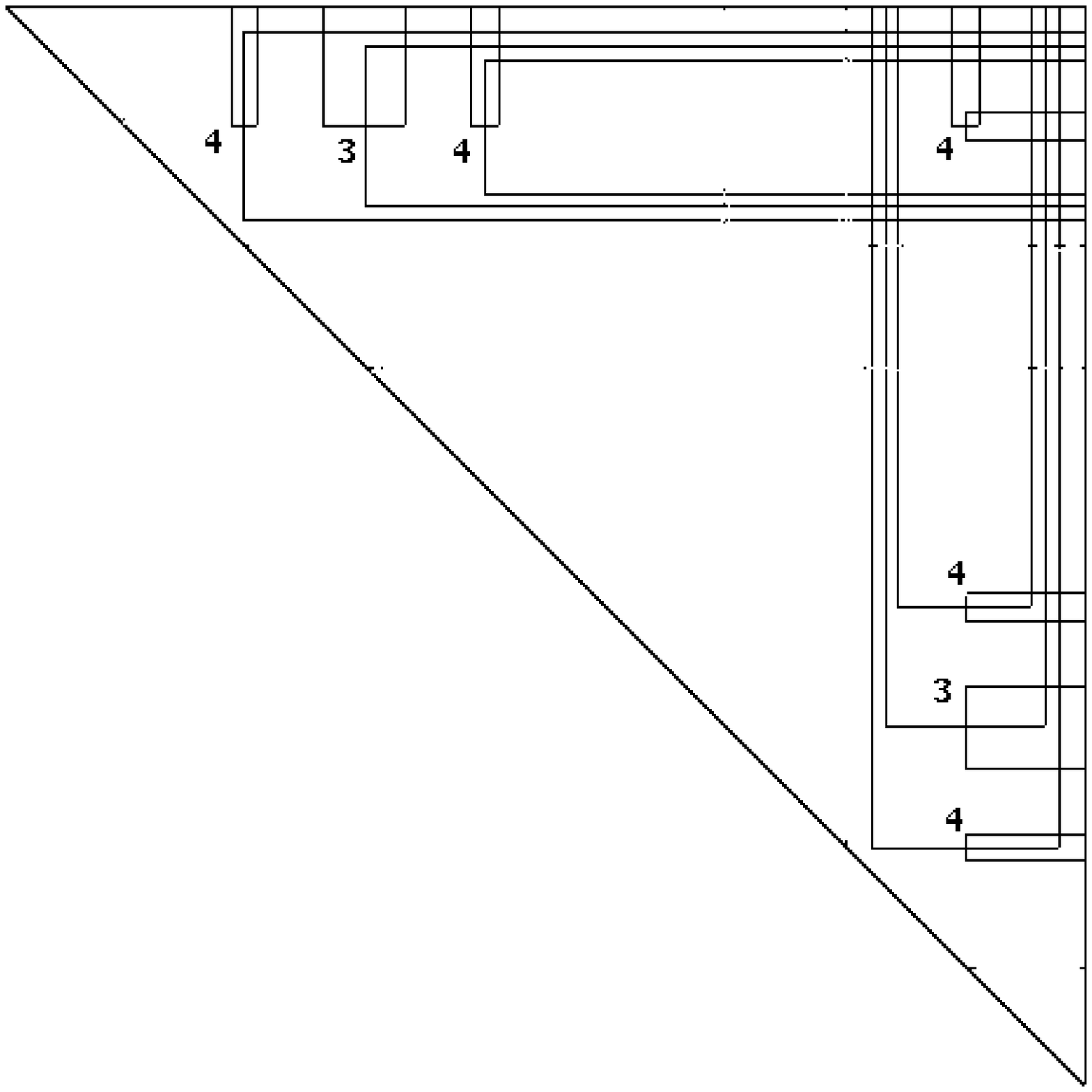,height=2.3in}\qquad
        \psfig{figure=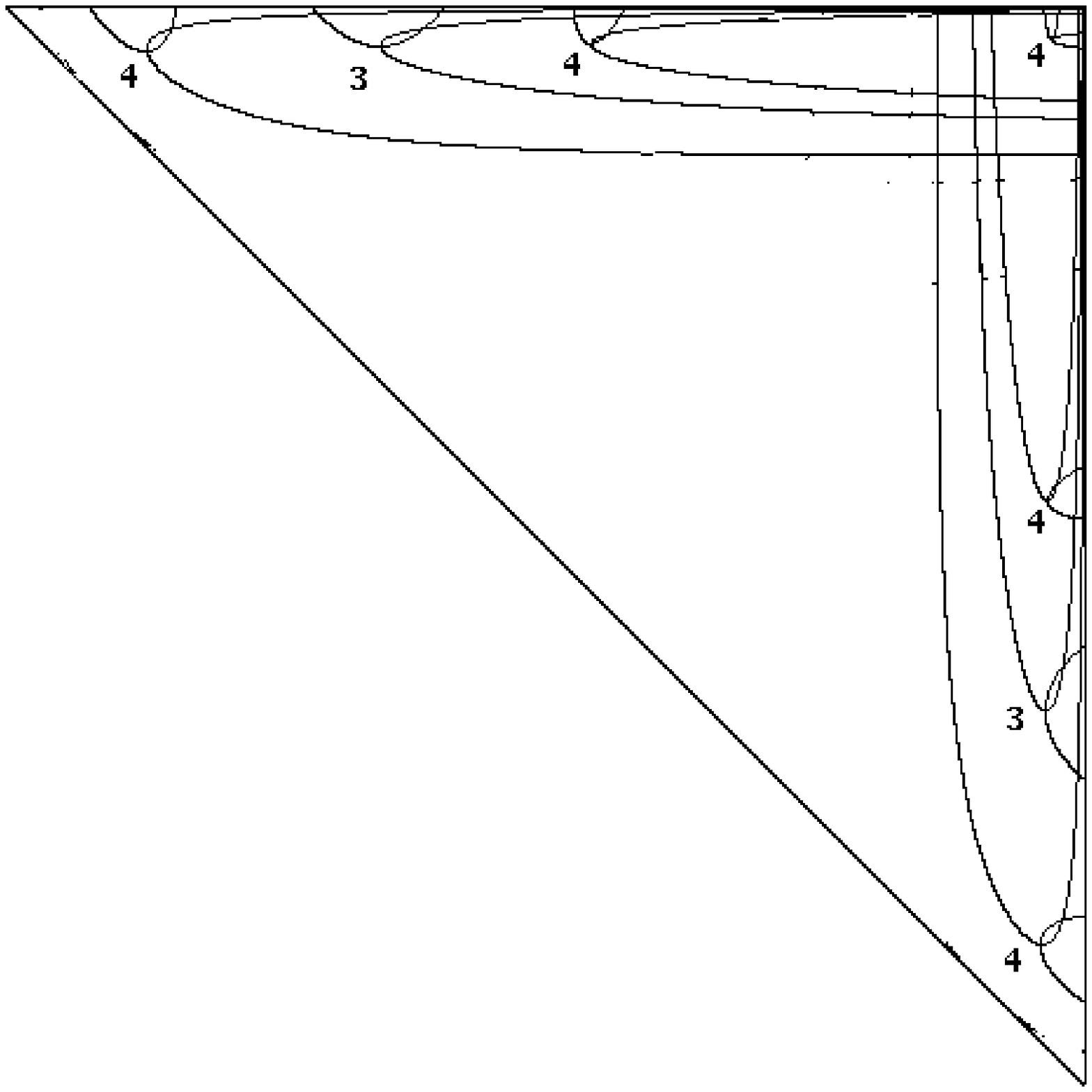,height=2.3in}}\smallskip
{\QP\bit Figure 13. Bones of period 3 and 4 for the stunted sawtooth
family of shape $(-+-)$ on the left and for the cubic family of shape
$(-+-)$ on the right. Periods are indicated near the primary intersection
points.\bigskip}
\endinsert

(Compare 4.5 for this last notation.) It follows that topological entropy
also increases monotonically as we traverse this edge.
Of course there is a completely
analogous statement as we traverse the right hand edge
$p_1=1$ in the direction of increasing $p_2$. It must be emphasized that
this lemma compares cubic and stunted sawtooth maps {\it for some fixed
shape
$\boldsig$.\/} If we switch shape, then the combinatorics becomes
completely
different.

{\bf Proof of 8.2.} At the left endpoint $\p=(0,1)$ of this edge, either in
the cubic or the stunted sawtooth family,  the associated map is monotone,
with no periodic orbit of negative type other than a single
fixed point of negative type in the $(-+-)$ case. On the other hand, at the
right endpoint $\p=(1,1)$, we have a bimodal map of maximal entropy, with
three non-overlapping subintervals each of which maps onto the entire
interval $I$. For such a
maximal entropy map
$f$ there are $3^n$ fixed points of $f^{\circ n}$ with distinct
itineraries,
and approximately half of these have negative type. Thus, as $p_1$
increases,
many periodic orbits of negative type must be created. However, an orbit of
negative type is a relatively robust object which can only be created or
destroyed at a parameter value for which one of the orbit points becomes a
folding point. Here only the left folding point $c_1$
can be the one in question,
since the condition $p_2=1$ guarantees that the orbit of the right folding
point $c_2$ never returns to $c_2$.

More precisely, we claim the following: {\it In either family, as
$p_1$ increases
through a value such that the orbit of $c_1$ is periodic, a periodic orbit
of negative type is created. This orbit persists until we reach the right
hand endpoint with $\p=(1,1)$, and all periodic orbits of negative type
$($with the exception of a single fixed point of negative type in the
$(-+-)$
case$)$ arise in this way.\/}

There is only one step we have to be careful about: We must be sure that an
orbit of negative type, once created, cannot disappear again as $p_1$
increases. But this would imply that there were two different parameter
values
along this edge for which the left folding point is periodic, {\it with the
same itinerary in both cases.\/} (It is not difficult to check that there
is
only one point of such an orbit which can merge with the left critical
point.)
This is impossible by Thurston's Theorem in the cubic case,\footnote
{As in the proof of 7.2, we must be careful since the map from parameter
triangle to positive affine conjugacy classes of cubic polynomials is not
one-to-one, but rather folds over two corners of the triangle. However
by 7.3, there
are no periodic critical orbits for parameters in the folded over region.}
or by 8.1 in the stunted sawtooth case.

Thus the invariant $\Cal N$ increases monotonically as we follow the edge,
and
increases as we cross each bone. But according to 4.9, $\Cal N$ is
determined
by the kneading data. Hence bone-ends with the same kneading data, in the
cubic and stunted sawtooth families, must occur in the same order, as
asserted.\QED

Now let us restrict attention to bones which are strictly
$\boldsig$-bimodal.
Thus we exclude only a few low period bones which are atypical
(as illustrated in Figure 14).

{\QP{\bf Lemma 8.3.} \it In the cubic family, just as in the stunted
sawtooth
family, each strictly $\boldsig$-bimodal left bone $B_-(\O)$ intersects the
corresponding right bone $B_+(\O)$ in exactly two points, and just one of
these (the ``primary intersection point'')
has the property that both folding points belong to a common periodic
orbit.
\ss}

{\bf Proof.} If there were no intersection point, then we could find a path
from the maximal entropy vertex $\p=(1,1)$ of the parameter triangle
which avoids both $B_-(\O)$ and $B_+(\O)$ and yet leads to the
opposite edge. This is impossible, since any negative orbit with strictly
bimodal order type $\O$ must disappear before we reach this lower edge,
which parametrizes only monotone maps. Hence there must be at least
one crossing point.

(Alternative Proof.
Compare Figure 12. As we follow the bone $B_-(\O)$ from one endpoint to the
other, the kneading invariant $\K_1$ remains periodic, of order type $\O$.
However, this kneading invariant cannot be the same at the two endpoints,
since
this would contradict Thurston's Theorem, and the only way for it to change
is to pass a point, where {\bit both} critical points belong to this same
periodic orbit.)

It follows from simple plane topology, using the Jordan curve
theorem, that the number of crossing points must be even, so there are
at least two. Finally, there cannot be more than two crossing points,
since each crossing point corresponds
to a postcritically finite cubic map, which is uniquely determined by its
kneading data $\bold{K}$. There are only two crossing points in the
stunted sawtooth case, hence there are only two possible choices for $\bold
 K$,
and the conclusion follows.\QED

{\bf Definition.} It will be convenient to divide
each strictly $\boldsig$-bimodal bone $B_\pm(\O)$ into two halves by
cutting
at the primary intersection point. Each of these two will be called a
{\bit half-bone\/}.

{\QP{\bf Theorem 8.4.} \it In either the cubic or the stunted sawtooth
family,
as we follow any left half-bone from its
primary intersection endpoint to the edge $p_2=1$, the negative orbit
complexity $\Cal N$ increases monotonically, with an actual increase every
time
we cross a right bone. Again there is a one-to-one correspondence between
crossing points in the two families which preserves the order along the
half-bone, and which preserves kneading data.\medskip}

The proof is completely analogous to the proof of 8.2. As we follow the
half-bone, the kneading sequence $\K_1$ remains constant. Hence an orbit
of negative type can appear or disappear only by passing through a
parameter
value for which the right folding point is periodic. Again there can be
only
one such parameter value for each periodic itinerary of negative type.
Details will be left to the reader.\QED

It follows that topological entropy increases monotonically as we traverse
either left half-bone, starting at the primary intersection point. Again
there
is a completely analogous statement for right half-bones.

{\bf Definition.}
By the $n\!$-{\bit skeleton\/} $S^{\roman{saw}}_n$ for the stunted
sawtooth family of some specified shape, either $\boldsig=(+-+)$ or
$\boldsig=(-+-)$, we will mean the union of all bones
$B^{\roman{saw}}_\pm(\O)\subset P^{\,\roman{saw}}$ of period $p\le n$,
together
with the boundary $\partial P^{\,\roman{saw}}$. Similarly, the $n$-skeleton
$S^{\roman{cub}}_n$ for the cubic family of the same shape
is the union of all cubic bones $B^{\roman{cub}}_\pm(\O)$ of
period $p\le n$, together with $\partial P^{\,\roman{cub}}$.
For elements of
the analysis of the structure of these skeletons, see [RS], [RT],
[MaT1], [MaT2] as well as [Mu].

\midinsert
\centerline{\psfig{figure=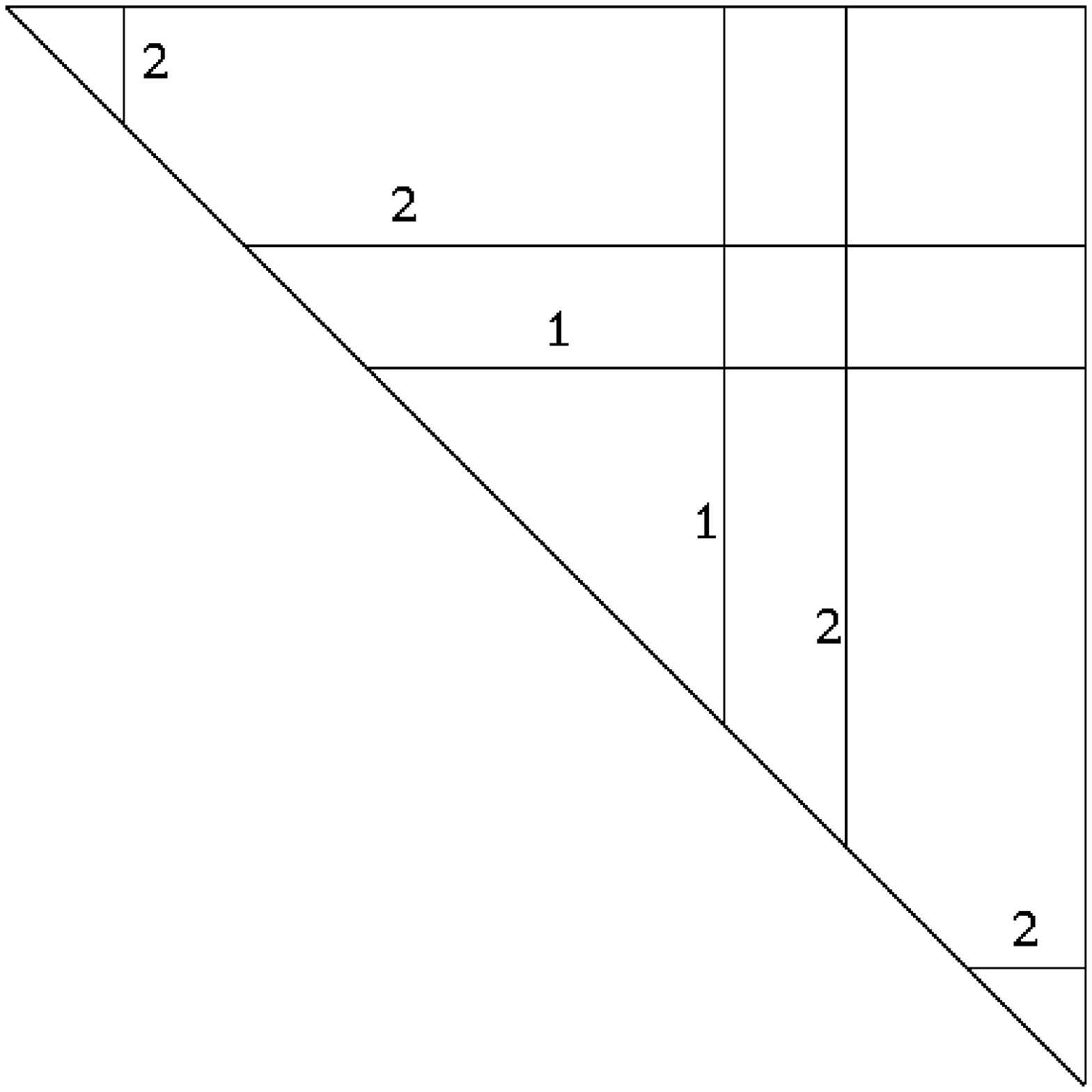,height=1.6in}\qquad
        \psfig{figure=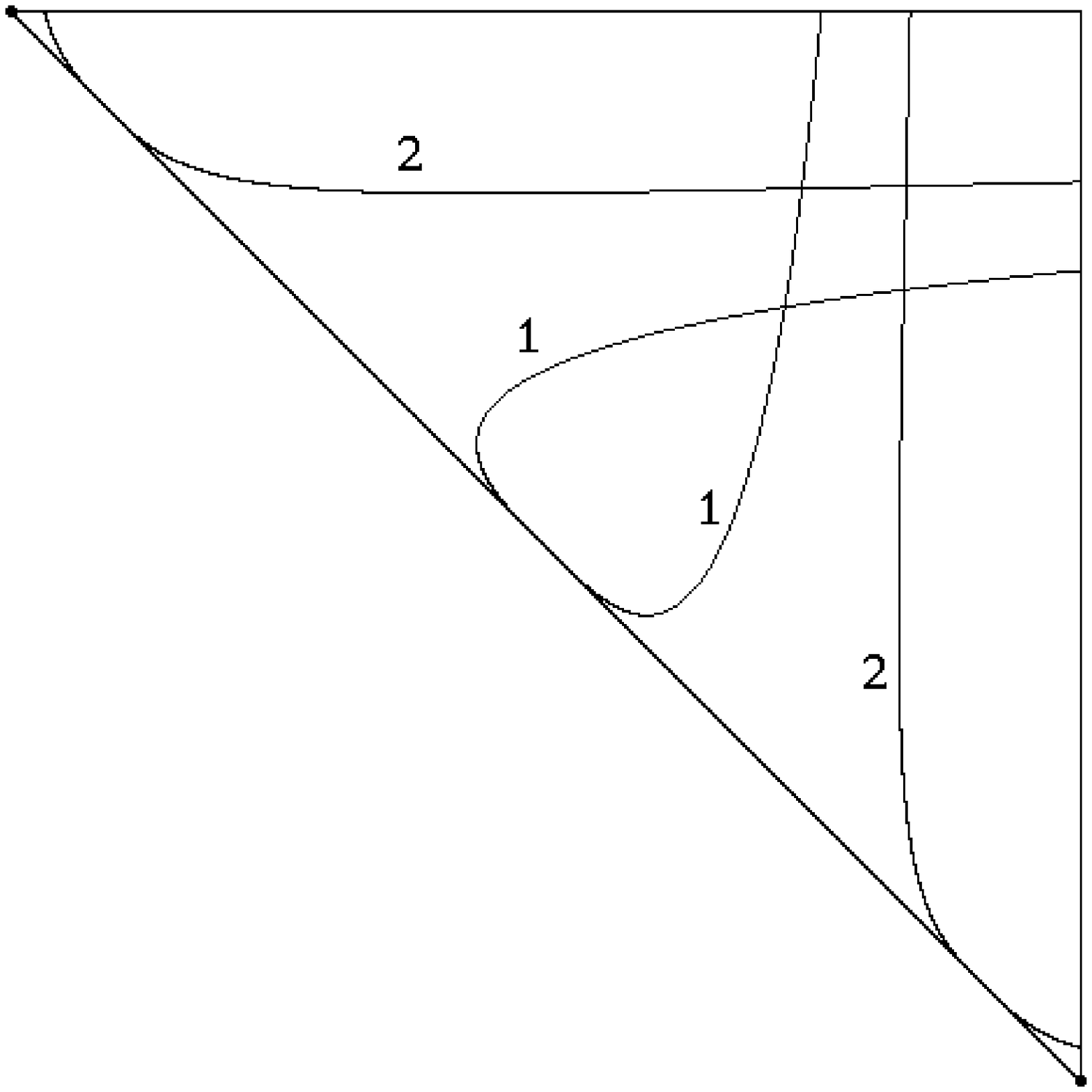,height=1.6in}}\smallskip
{\QP\bit Figure 14. The bones of period $\le 2$ for the
shape $(-+-)$ are not strictly bimodal, and hence are quite atypical.
In particular, their structure  in the stunted sawtooth case (shown on
the left) is qualitatively different from that
in the cubic case (on the right). The period one bones
for the $(+-+)$ case are even simpler, with no intersection point at all.
In order to get a homeomorphism
from one parameter triangle to the other preserving these low period bones,
we would have to first remove the bottom edge
of the triangle or collapse it to a point.
\smallskip}
\endinsert

It will be convenient to exclude the low period bones (those whose orbit
types are not strictly bimodal), since they behave rather awkwardly.
(Compare Figure 14.) Thus, in either family, we define
the {\bit essential $n$-skeleton\/}
$E^{\roman saw}_n$ or $E^{\roman cub}_n$ to be the boundary of the
parameter
triangle together with the union of all bones
associated with {\it strictly\/} $\boldsig$-bimodal order types $\O$ of
period
$p\le n$. In other words, we require
that $2\le p\le n$ in the $(+-+)$ case, and that $3\le p\le n$ in the
$(-+-)$ case. By a {\bit vertex\/} of this essential skeleton, we mean
either
an endpoint of one of its bones or an intersection point of two of its
bones.

As an example, Figure 13 shows the essential 4-skeletons for the stunted
sawtooth and cubic families of shape $(-+-)$. In comparing the pictures
for these two families, it is important that we consider
only the {\it essential\/} skeletons, since low period bones do not behave
in
quite the same way.

{\QP{\bf Theorem 8.5} \it For either bimodal shape $\boldsig$, and for any
$n>2$ there is a homeomorphism from the parameter triangle
$P^{\,\roman{saw}}$
to the parameter triangle $P^{\,\roman{cub}}$ which maps the essential
skeleton $E_n^{\roman{saw}}$ homeomorphically onto the essential skeleton
$E_n^{\roman{cub}}$, carrying each bone of order type $\O$ to a bone of
the same order type  $\O$, carrying each vertex
to a vertex with the same kneading data, and carring each edge of the
parameter triangle to the corresponding edge.\medskip}

The proof, making use of 8.2 and 8.4 and the fact that all intersections
are
transverse, is a straightforward exercise in plane topology. For example
one can proceed inductively, starting with the identity map of $P^2$
and then adjusting it to behave correctly on one bone at a time. Details
will be left to the reader. \QED\bigskip

\maybebreakhere
\centerline{\bf\S9. From Connected Bones to Connected Isentropes.}
\smallskip

As usual. we fix one of the two possible shapes $\boldsig$ for bimodal
maps.
Working either in the cubic family or the stunted sawtooth family,
we have the following.

{\QP{\bf Lemma 9.1.} \it
Let $\p$ and $\p'$ be two points in the parameter triangle $P^2$
such that the associated maps have topological
entropy $h(f_\p)\ne h(f_{\p'})$.
Then any path from $\p$ to $\p'$ in the parameter triangle
$P^{\,\roman{cub}}$ must cross infinitely many bones.\smallskip}

{\bf Proof of 9.1} According to Theorem 4.11
the difference $|\roman{Neg}(f_\p^{\circ k})-
\roman{Neg}(f_{\p'}^{\circ k})|$ must be unbounded as $k\to\infty$.
As we deform $\p$ along some path in $P^2$, the number
${\roman{Neg}}(f_\p^{\circ k})$, which measures the number of decreasing
laps of
$f_\p^{\circ k}$ whose graph crosses the diagonal, will remain constant
except as we pass through a map which has a folding point with periodic
orbit
of period $q$ dividing $k$, with $k/q$ odd, in which case this number jumps
 by
$\pm q$. (Compare 7.4.) Thus, ${\roman{Neg}}(f^{\circ k})$ remains
constant unless we pass through a bone of period dividing $k$.
Further details are easily supplied.\QED

Associated with the essential skeleton $E_n^{\roman cub}$ or
$E_n^{\roman{saw}}\subset P^2$ is
a {\bit topological cell structure\/} on $P^2$. That is, we can
partition $P^2$ into subsets, each of which is homeomorphic
either to a point, an open interval, or a 2-dimensional
open unit disk. Furthermore, these subsets
fit together nicely so that the closure of each one is topologically
a point, a closed interval, or a closed 2-disk.
By definition, the open 2-cells in this cell structure are the
connected components of the complement $P^2 \smallsetminus
E_n$, %^{\roman{saw}}$,
the 0-cells are the vertices as described in
Theorem 8.5, and the open 1-cells are the connected
components of $E_n \smallsetminus\{{\roman{vertices}}\}$.
The resulting cell complex will be denoted by $P_n^2$, or more
precisely by $P_n^{2,\,\roman{cub}}$
or $P_n^{2,\,\roman{saw}}$. %There is a completely analogous cell
%complex $P_n^{2,\,\roman{cub}}$ for the cubic family.
These two cell complexes $P_n^{2,\,\roman{cub}}$
and $P_n^{2,\,\roman{saw}}$ are homeomorphic as a consequence of
Theorem 8.5, in a homeomorphism $\eta _n$ which
takes each vertex to a vertex of the same topological entropy and
each edge to an edge with the same interval of entropies.

{\QP{\bf Lemma 9.2.} \it For each $\epsilon>0$ there exists an integer
$n$ such that if two points $\q$ and $\q'$ belong to a common
closed cell of the complex $P_n^{2,\,\roman{cub}}$ or
$P_n^{2,\,\roman{saw}}$, then the corresponding maps $f_\q$ and $f_{\q'}$
satisfy $|h(f_\q)-h(f_{\q'})|<\epsilon$.\ss}

\comment
$$  |h(S_{\q'})-h(S_\q)|~<~\epsilon~. $$
If the Connected Bone Conjecture is true, then there is an analogous
statement for the complex $P_n^{2,\,\roman{cub}}$ and the
family of cubic maps $f_\p$.\smallskip}
\endcomment

{\bf Proof in either family.} Otherwise we could find $\epsilon>0$ so that
for each $n$ there existed points $\q_n$ and $\q'_n$ in a common cell of
$P_n^2$ with $|h(f_{\q'_n})-h(f_{\q_n})|\ge
\epsilon$. After passing
to infinite subsequences, we could assume that both sequences converge,
say $\q_k\to\q$ and $\q'_k\to\q'$,
and furthermore that all of the $\q_k$ and $\q'_k$
belong to a common closed cell of $P_n$
whenever $k\ge n$.
Thus the limit points $\q$ and $\q'$ would belong to a common closed
cell of $P_n^{2}$ for every $n$, but by continuity
the associated entropies would differ by at least $\epsilon$. This is
impossible by Lemma 9.1.\QED

{\QP{\bf Lemma 9.3.} \it The topological
entropy function $\q\mapsto h(f_\q)$, either for the cubic family or
the stunted sawtooth family, restricted to any closed cell of the cell
complex $P_n^{2}\!$, takes its maximum
and minimum values on the boundary (and in fact on the set
of boundary vertices). Hence the interval of entropy values which are
realized
by any cell in the cubic family is the same as the interval of entropy
values
for the corresponding cell in the stunted sawtooth family.
% If the Connected Bone Conjecture is true, then
%the analogous statement holds for the topological entropy function
%$\p\mapsto h(f_\p)$ for cubic maps, restricted to any
%closed cell of the cell complex $P_n^{2,\,\roman{cub}}\!$.
\smallskip}

{\bf Proof.} For the stunted sawtooth family, the first statement
follows easily from the fact that entropy is a monotone function of either
coordinate $p_1$ or $p_2$. To prove the analogous statement
for the cubic family, suppose for example
that for some point $\p_0$ of a closed cell $C$, the value
$h(f_{\p_0})$ were strictly larger than the maximum
value $h_{\max}$
on the boundary of $C$. Let $2\epsilon = h(f_{\p_0})-h_{\max}$.
According to Lemma 9.2 we can choose $m>n$ so that $h$ varies
by less than $\epsilon$ on each cell of $P_m^{2,\,\roman{cub}}$. Let
$C'\subset C$ be a cell of $P_m^{2,\,\roman{cub}}$
which contains this point $\p_0$,
and let  $\p'$ be any vertex of $C'$. Then it follows that
$h(f_{\p'}) > h_{\max}$.
Since the homeomorphism $\eta_m$ carries
vertices to vertices with the same topological
entropy, this would yield a vertex in
the complex $P_m^{2,\,\roman{saw}}$ satisfying a corresponding
inequality. But this is impossible, since we already proved
Lemma 9.3 to be true for the stunted sawtooth family. The second statement
follows immediately.\QED

{\QP{\bf Main Theorem 9.4.} \it Every isentrope $\{\p\in
P_\pm^{2,\roman{cub}}
 ~;~ h(f_\p)=h_0\}$ for either cubic family  is connected.\medskip}

{\bf Proof.} For each $n$, the union of all
closed cells of $P_n^{2,\,\roman{saw}}$
which touch the $h_0\!$-isentrope forms a compact set, which is connected
by Theorem 6.1. The
union of corresponding cells in $P_n^{2,\,\roman{cub}}$
forms a compact connected set by Theorem 8.5, and
contains the $h_0$-isentrope for the cubic family by Lemma 9.3.
The intersection of these sets, as $n\to\infty$, will be precisely equal
to the required isentrope by Lemma 9.2. Since an intersection of compact
connected sets is compact and connected, the conclusion follows.\QED

Recall that a compact subset of $\R^n$ is called {\bit cellular\/} if
it is the intersection of a nested sequence of closed topological
$n$-cells,
each contained in the interior of its predecessor. (Compare [Br].)

{\QP{\bf Corollary 9.5.} \it Each isentrope $\{h=h_0\}$ for either cubic
family is a cellular set. Further, if $0<h_0<\log(3)$ then this isentrope
separates the parameter triangle into two connected subsets.\ss}

{\bf Proof.} Since we are in dimension two, it is only nessary to check
that a
set is compact and connected, with connected complement in $\R^2$, in order
 to
prove that it is cellular. But it follows from 9.4 that no isentrope can
separate the plane. For otherwise, if $\p_0$ belonged to the bounded
component
of the complement of some isentrope $\{h=h_0\}\subset P^2$ within $\R^2$,
then choosing a point on
the boundary of the parameter triangle with the same topological entropy as
$p_0$, we would see that the isentrope $\{h=h(\p_0)\}$ could not be
connected, contradicting 9.4.
Since we know by 8.2 that each isentrope intersects the two top edges of
the parameter triangle in connected sets, it follows easily that the
complement
$P^2\ssm\{h=h_0\}$ has at most two connected components.\QED

\centerline{------------------------------}\ss\bigskip
%\vfil\eject

\maybebreakhere
\centerline{\bf Appendix A :
Characterization of a Polynomial by its Critical Values.}\medskip

\centerline{ A letter from Adrien Douady (slightly edited).}\medskip

%\centerline{--------------------------}\medskip

\rightline{ Tuesday, 13 July 1993}

%Cc: sentenac@matups.matups.fr (Pierrete Sentenac)\bigskip\bigskip

\noindent Dear Milnor,\ss

Putting order in my papers, I found some notes
I had written with Pierrette
Sentenac after a talk by Arnold in ENS, where he asked for
a proof in $\R$. Let:

${\Cal P}_n = \{$ monic centered real polynomials of degree  $n \}
     = \{x^n + c_{n-2}x^{n-2} +\cdots+c_0 \}~, $\ss

$\Omega = \{ f \in {\Cal P}_n~ |~ f$ has distinct real roots $ a_1 ,
\ldots, a_n\}~.$\ss

\noindent Then

$f \mapsto (a_1 ,\ldots, a_n )$ is a diffeo of  $\Omega$  onto
  $\{(a_1,\ldots,a_n)~ |~ a_1<\cdots<a_n ,~ \sum a_i = 0 \} ,$\ss

\noindent or equivalently\ss

$f \mapsto (\ell_1,\ldots,\ell_{n-1} )$ , where
$\ell_i = a_{i+1} - a_i$ ,
is a diffeo of $\Omega$ onto $(\R_+)^{n-1}$ .\ss

Let $\Phi(f)=(s_1,\ldots,s_{n-1})$ where
 $$s_i = \int_{a_i}^{a_{i+1}} |f(x)| dx = (-1)^{n-i}
\int_{a_i}^{a_{i+1}} f(x) dx~.$$\ss

{\bf Proposition.} $\Phi: f \mapsto (s_1,\ldots,s_{n-1} )$  is a diffeo
of  $\Omega$ onto  $(\R_+)^{n-1}$ .\medskip

{\bf Proof.}  a) $\Phi$  is a local diffeo. In fact,
for  $g$  of degree  $n-2$ , we have
$$   (-1)^{n-i}
 {d \Phi_i(f+\epsilon g)\over d\epsilon}~=~\int_{a_i}^{a_{i+1}}
 g(x)\, dx ~+~{da_{i+1}\over d\epsilon} f(a_{i+1})~-~
 {da_{i}\over d\epsilon} f(a_{i}) $$
at $\epsilon=0$,
where the last two summands vanish since the $a_i$ are zeros of $f$. Thus
%$$\Phi_i ( f + \epsilon g) = \Phi_i (f) +  (-1)^{n-i} \epsilon
%\int_{a_i}^{a_{i+1}} g(x)\, dx  + O(\epsilon^2) .$$
%The linear term
$d\Phi/d\epsilon|_{\epsilon=0}$ could vanish only
if  $g$  had zero average
on each of the $(n-1)$ intervals  $[a_i,a_{i+1}]$ , and hence had at least
$n-1$  roots; which is impossible for  $g \ne 0$ of degree $n-2$.\medskip

b)  $\Phi$ is proper. In fact  $s_i\to \infty$ as the length
$\ell_i=a_{i+1}-a_i$ tends to
$\infty$, and $s_i\to 0$ as $\ell_i$ tends to $0$
with $\sum\ell_j$ bounded. (Use the fact that $f(x)=\prod (x-a_i)$.)
 \medskip

{\bf Corollary.} Let $\Omega'\subset {\Cal P}_{n+1}$ be
the set of monic centered
polynomials $F$ of degree $n+1$ with $n$ distinct real critical points.
Taking $(n+1)f$ to be the derivative of $F$, it follows that
the map  $F \mapsto (v_1,\ldots,v_n)$  from  $\Omega'$  to  $\R^n$,
where $v_i$ are the critical values (in the order of the critical points)
is a diffeo of  $\Omega'$ onto the set of $(v_1,\ldots,v_n)$ with
${\sgn} (v_i - v_{i-1}) = (-1)^{n-i}$.\medskip

\hskip .5in  Best wishes,\ss\ss

\hskip 1in  Adrien\ss
%\eject
\centerline{------------------------------}\ss

{\bf Addendum.} This argument can easily be extended to the case where
the roots $\{a_j\}$ of $f$ (alias the critical points of $F$) need not be
distinct, but must only satisfy\break $a_1\le a_2\le\cdots\le a_n$ with
$\sum a_i=0$. Let $\hat a_1<\cdots<\hat a_r$ be the distinct roots
of $f$, with multiplicities $m_i\ge 1$ so that
$$   f(x)\= (x-\hat a_1)^{m_1}\,\cdots\,(x-\hat a_r)^{m_r}~,$$
with
$$   m_1+\cdots+m_r\=n\qquad\text{and}\qquad
 m_1\hat a_1+\cdots+m_r\hat a_r\=0~.$$
To deform
$f$ within polynomials of this same form, choose direction parameters
$w_i\in\R$ with $m_1w_1+\cdots +m_rw_r=0$ and set
$$   f_t(x)\=\big(x-\hat a_1(t)\big)^{m_1}\,\cdots\,\big(x-a_r(t)
 \big)^{m_r}\qquad \text{where}\qquad \hat a_i(t)\=\hat a_i+tw_i~. $$
As before, let
$$   s_i(t)\=\int_{\hat a_i(t)}^{\hat a_{i+1}(t)}|f_t(x)|\,dx
 \=\pm\int_{\hat a_i(t)}^{\hat a_{i+1}(t)}f_t(x)\,dx~. $$
A brief computation shows that the
derivative $ds_i/dt$ at $t=0$ can be expressed as the
integral from $\hat a_i$
to $\hat a_{i+1}$ of the product of a fixed polynomial
$$ \mp(x-\hat a_1)^{m_1-1}\,\cdots\,(x-\hat a_r)^{m_r-1} $$
which has no zeros within this interval, %$\hat a_i<x<\hat a_{i+1}$
and a polynomial
$$   g(x)\=\sum_i m_iw_i\prod_{\{j\,;\,j\ne i\}} (x-\hat a_j) $$
of degree $\le r-2$. This polynomial $g(x)$ is non-zero
unless all of the $w_i$
are zero, since %its value at $\hat a_i$ is $
$$   g(\hat a_i)\=m_iw_i\prod_{\{j\,;\,j\ne i\}}(\hat a_i-\hat a_j)~, $$
which is non-zero unless $w_i=0$.
Hence $g$ has at most $r-2$ zeros. Arguing as above,
it follows that the linear first derivative map
$$   (w_1\,,\,\ldots\,,\,w_r)~\mapsto~(ds_1/dt\,,\,\ldots\,,\,ds_{r-1}/dt)
 |_{t=0} $$
is bijective, and it follows as before that $F$ is uniquely determined by
its critical value vector $(v_1,\ldots, v_n)$, where now the case
$v_i=v_{i+1}$
is allowed.

\comment
a2222222222222222222
\endcomment
%\vskip .5in

\bigskip
\maybebreakhere
\centerline{\bf Appendix B: Tight Symbol Sequences and Thurston's Theorem.}
\bigskip

We first make a rather formal definition, to be
justified by the properties
which follow. Let $\I(x)=(A_0\,,\,A_1\,,\,A_2\,,\,\ldots)\in\A^\N$ be the
itinerary of a point under some $m$-modal map of shape $\boldsig$
with kneading data $\bold{K}=(\K_1\,,\,\ldots\,,\,\K_m)$.\smallskip

{\bf Definition.} We will say that this symbol sequence $\I(x)=
(A_0\,,\,A_1\,,\,\ldots)$ is {\bit flabby\/} (flasque)
if some {\bit terminal segment}
$$\roman{shift}^{\circ k}\I(x)\=(A_k\,,\,A_{k+1}
\,,\,\ldots)$$
has either the form $(I_j\,,\,\K_j)$ or the form
$(I_{j-1}\,,\,\K_j)$, for any $1\le j\le m$. In other words, some point of
the associated
orbit which is not a folding point must have the same itinerary as an
immediately adjacent folding point. This symbol sequence will be called
{\bit tight\/} if it is not flabby.
Similarly, we will say that
the kneading data $\bold{K}$ is {\bit tight\/}
if each of the sequences $\K_1\,,\,\ldots\,,\,\K_m$ is tight.\smallskip

As a first example, consider the family of stunted sawtooth maps
with fixed shape $\boldsig$ and fixed slope $s>m+1$.

{\QP{\bf Lemma B.1.} \it Let $S_\p:J\to J$ be a stunted sawtooth map in
this family with kneading data $\bold{K}$,
and let $\I^0$ be a symbol sequence in $\A^\N$. There exists
one and only one $x\in J$ with itinerary $\I^0$ if and only if this symbol
sequence is tight and $\bold{K}$-admissible.

Similarly, for any kneading data $\bold{K}$, there exists one
and only one map $S_\p$ in the family with this kneading data if and only
if
$\bold{K}$ is tight and $\boldsig$-admissible.\medskip}

(Existence was shown in \S5, so only the uniqueness part is new and
requires
tightness.) Proofs will be given at the end of this appendix. We can also
give
a more geometric criterion for tightness in the stunted sawtooth case.

{\QP{\bf Lemma B.2.} \it The itinerary of $x$ under $S_\p$ is tight if and
only if the orbit of $x$ never hits any plateau except at its central
folding point. Similarly the kneading data for $S_\p$ is tight if and
only if the orbit of each folding point satisfies this condition.\medskip}

Recall from \S5 that any admissible kneading data is represented by
a canonical stunted sawtooth map $S_\p$, which has the property that the
orbit
of a folding point cannot hit the {\it interior\/} of a plateau except
at a folding point. Thus, if we use this canonical map, the only
extra requirement for tightness of $\bold{K}$ is that the
folding point orbits cannot even hit the endpoints of the
plateaus.\smallskip

The case of polynomial maps is much more difficult, and not well
understood.
However, for the {\bit postcritically finite\/} case (where every
critical
orbit is periodic or eventually periodic) we have the following
result,\footnote
{{\bf Erratum.} This result was stated incorrectly, without the hypothesis
of tightness, in our preliminary manuscript [DGMT, \S5]. A counter-example
to the version stated there is provided in Example 3 below.
(Figure 15-c.)}
which is an adaptation of a basic theorem of Thurston in a
form due to Poirier.

{\QP{\bf Theorem B.3.} \it Suppose that the $m$-modal
kneading data $\bold{K}$ is\break
$\boldsig$-admissible, with $\K_i\ne\K_{i+1}$ for all $i$.
There exists a postcritically finite polynomial map of degree $m+1$
and shape $\boldsig$ with kneading data $\bold{K}$
if and only if each $\K_i$ is periodic or eventually periodic,
and also tight. This polynomial is always unique when it
exists, up to a positive affine change of coordinates, or as a boundary
anchored map of the interval.\medskip}

To illustrate this result, here are four examples of non-tight data.
The notation $\overline A$ will stand for an infinite sequence $(A\,,\,A\,,
\,A\,,\,\ldots)$ of identical copies of the symbol $A\in\A$.\smallskip

\midinsert
\centerline{\psfig{figure=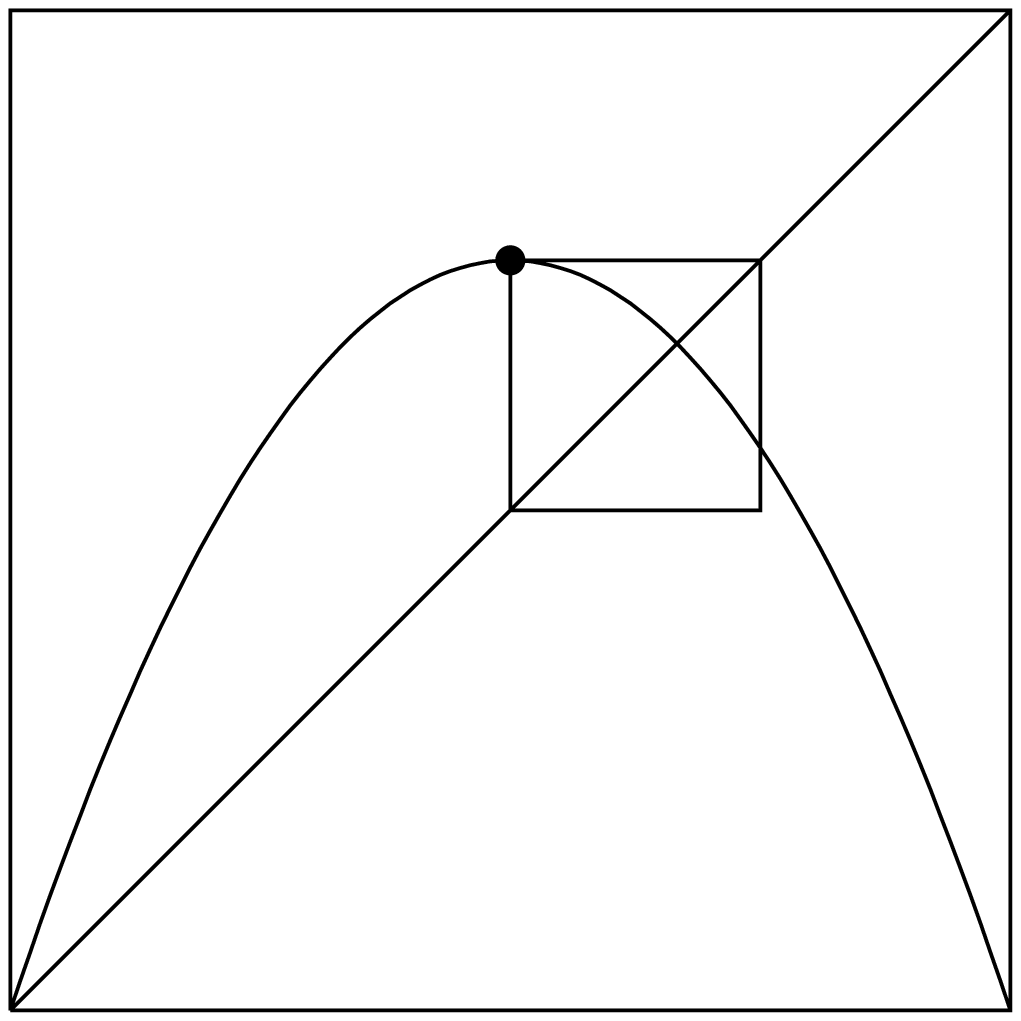,height=1.4in}\qquad
        \psfig{figure=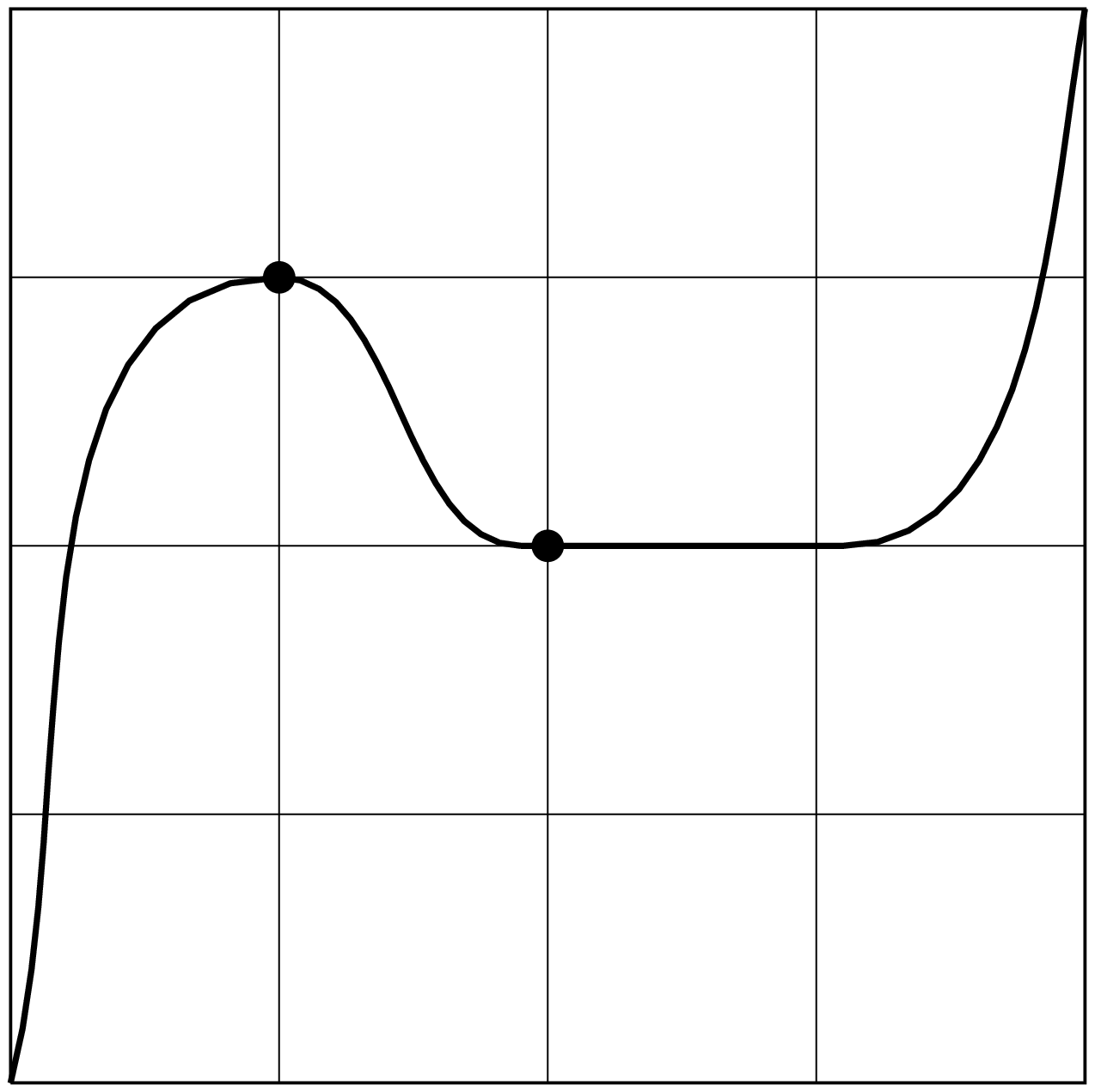,height=1.4in}\qquad
        \psfig{figure=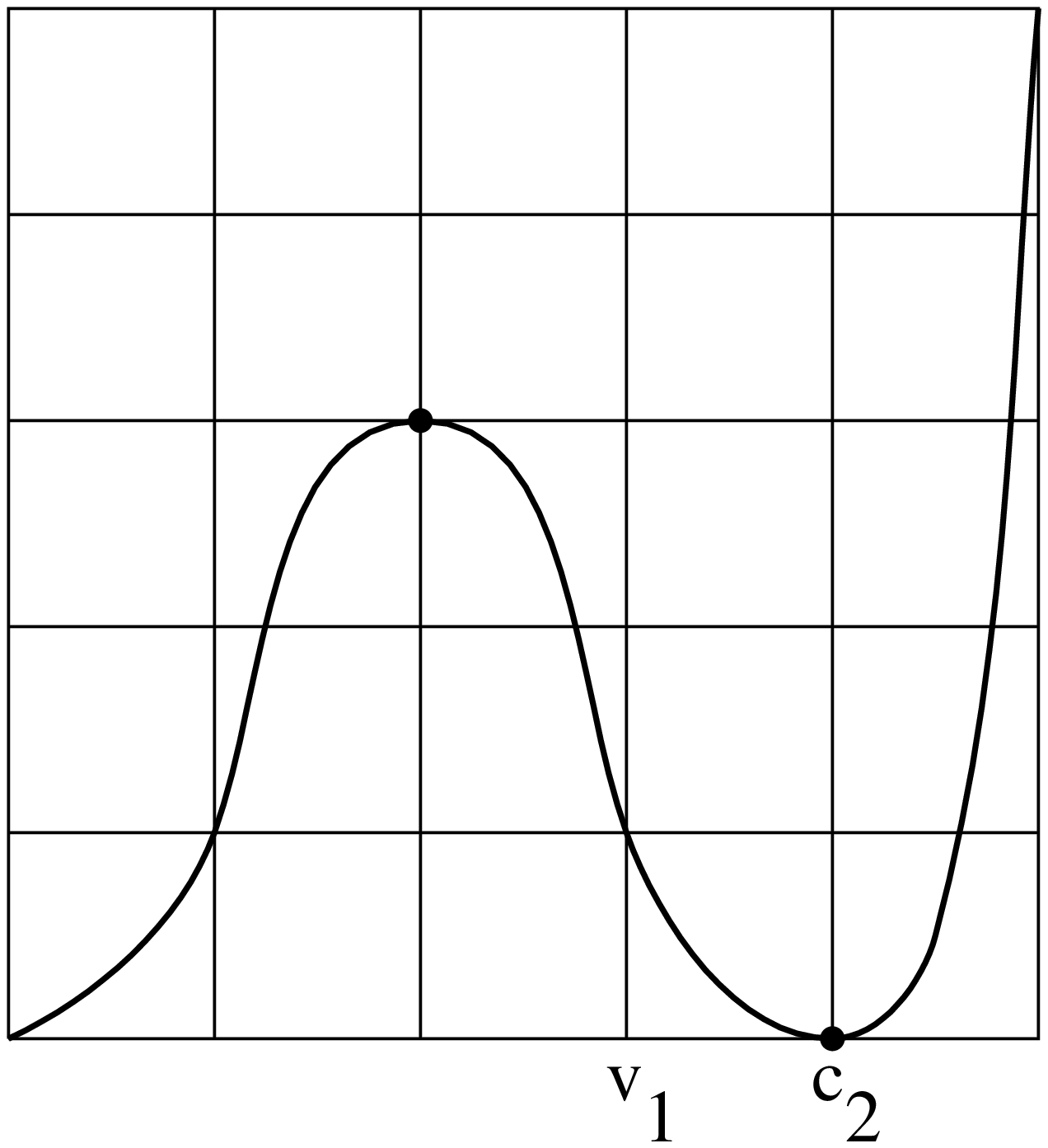,height=1.4in}}
{\QP\bit Figure 15. Three maps with eventually periodic but non-tight
kneading
data. (See Examples 1, 2, 3.) The folding points are marked.\smallskip}
\endinsert

{\bf Example 1.} If $2<4v<1+\sqrt 5$, the unimodal map $Q_v(x)=4vx(1-x)$
has periodic kneading sequence $~\K_1=\overline I_1=(I_1\,,\,\K_1)~$.
This kneading
sequence is not tight. It is realised by an entire one-parameter family of
distinct quadratic maps (Figure 15-a), but not by any postcritically
finite quadratic map.

{\bf Example 2.} The $(+-+)$-bimodal kneading data $~\K_1=(I_2\,,\,
\overline C_2)\;,\;\K_2=\overline C_2~$
% $$ \K_1\=(I_2\,,\,C_2\,,\,C_2\,,\,\ldots)\;,\quad \K_2\=(C_2\,,\,C_2\,,\,
%    \ldots) $$
is admissible, but cannot be realised by any piecewise strictly monotone
map.
(Figure 15-b. Compare the discussion of 8.1.)

{\bf Example 3.} The $(+-+)$-bimodal data $~\K_1=(I_1\,,\,\overline I_0)\;,
\;\K_2=\overline I_0~$
%$$  \K_1\=(I_2\,,\,I_0\,,\,I_o\,,\,\ldots)\;,\quad
%    \K_2\=(I_0\,,\,I_0\,,\,\ldots) $$
can be realised by a postcritically finite
4-th degree polynomial, but not by any cubic polynomial.
(Figure 15-c.)

{\bf Example 4.} (See [MaT1, p. 179].) More generally, whenever
the kneading data for $f$
satisfies $\K_1=(I_1\,,\,\K_2)$, it follows that all
points in the open interval $(v_1\,,\,c_2)\subset I_1$ have
the same itinerary.
If $\K_2$ is not eventually periodic, then $(v_1\,,\,c_2)$
is a wandering interval, hence $f$ cannot be real
analytic.\footnote{Compare [L1], [BL], [MvS1], [MMS]. It would be
interesting
to know whether any such kneading data can be realized by a map which is
$C^\infty$ and piecewise strictly monotone.  Examples of $C^\infty$
maps with wandering intervals have been given in in [Ha] and [SI].}
There are uncountably many distinct choices for
$\K_2$. In fact, using the
stunted sawtooth model with folding points $\pm 1$, we can take $v_1=1-
\epsilon$ and let $v_2$ range from $-\alpha$ to $v_1$. The corresponding
topological entropy ranges between
$\log 2$ and zero. (For example, the right picture in Figure 15
illustrates the $\log 2$ case.) Thus we can certainly choose infinitely
many $\K_2$ which are not eventually periodic.
\smallskip

For the application, we need a special case where the
hypothesis of tightness is automatically satisfied:

{\QP{\bf Corollary B.4.} \it Let $\bold{K}$ be admissible kneading
data such that every kneading sequence $\K_j$ either

1) contains the symbol $C_j$, so that the folding point itinerary
$(C_j\,,\,\K_j)$
is periodic, or

2) satisfies the condition that $\K_j$ equals either $\I_\min$ or $\I_\max$
according as the sign $\sigma_j$ is $+1$ or $-1$.

\noindent Then this data is tight, and hence is realized
by one and only one postcritically finite polynomial in the parametrized
family of \S3, and by
one and only one stunted sawtooth map in the family of \S5.\medskip}

(However, compare Figure 16.)
Here the various periodic critical orbits need not be disjoint.
It follows that there is a canonical one-to-one correspondence
between intersection of bones in the cubic family and in the corresponding
stunted sawtooth family, and also between endpoints of bones in these
two families. This would yield an alternate proof of 8.5.

\midinsert
\centerline{\psfig{figure=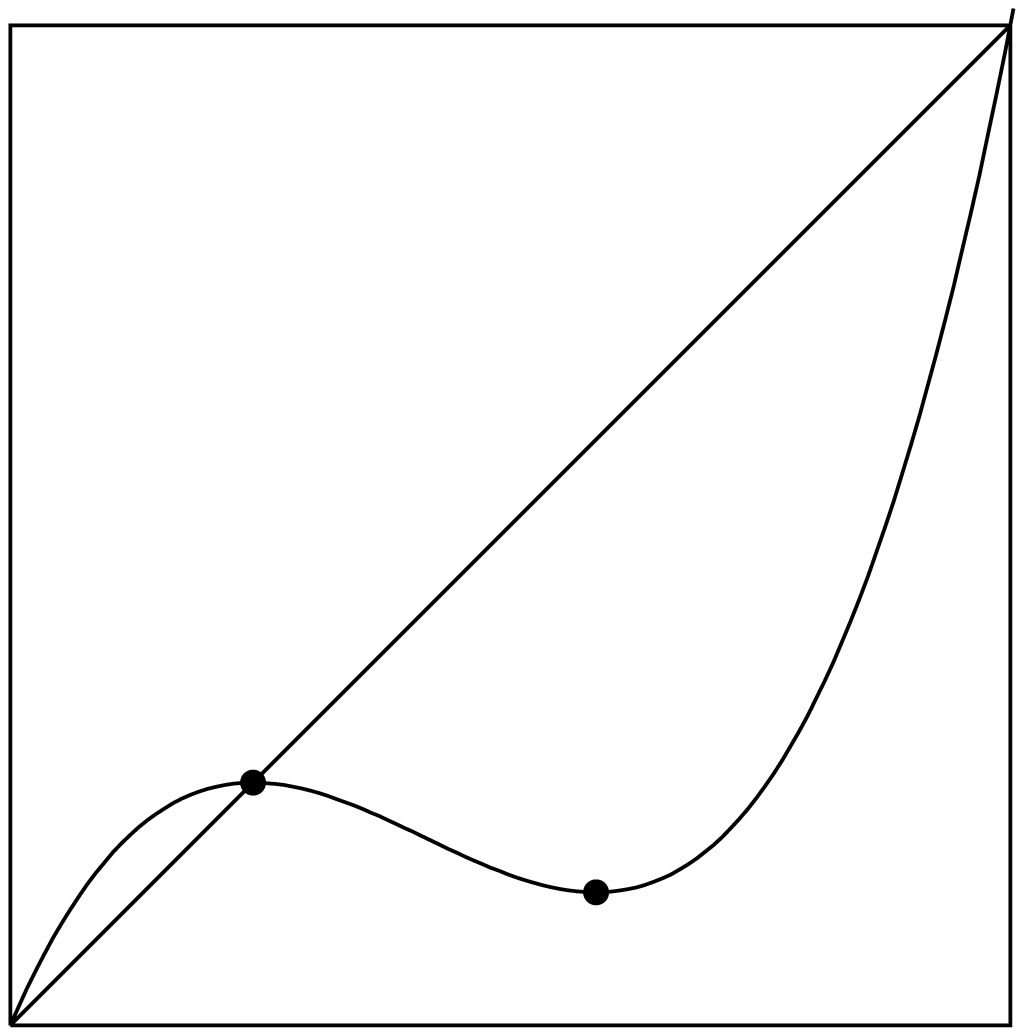,height=1.2in}\qquad\qquad
        \psfig{figure=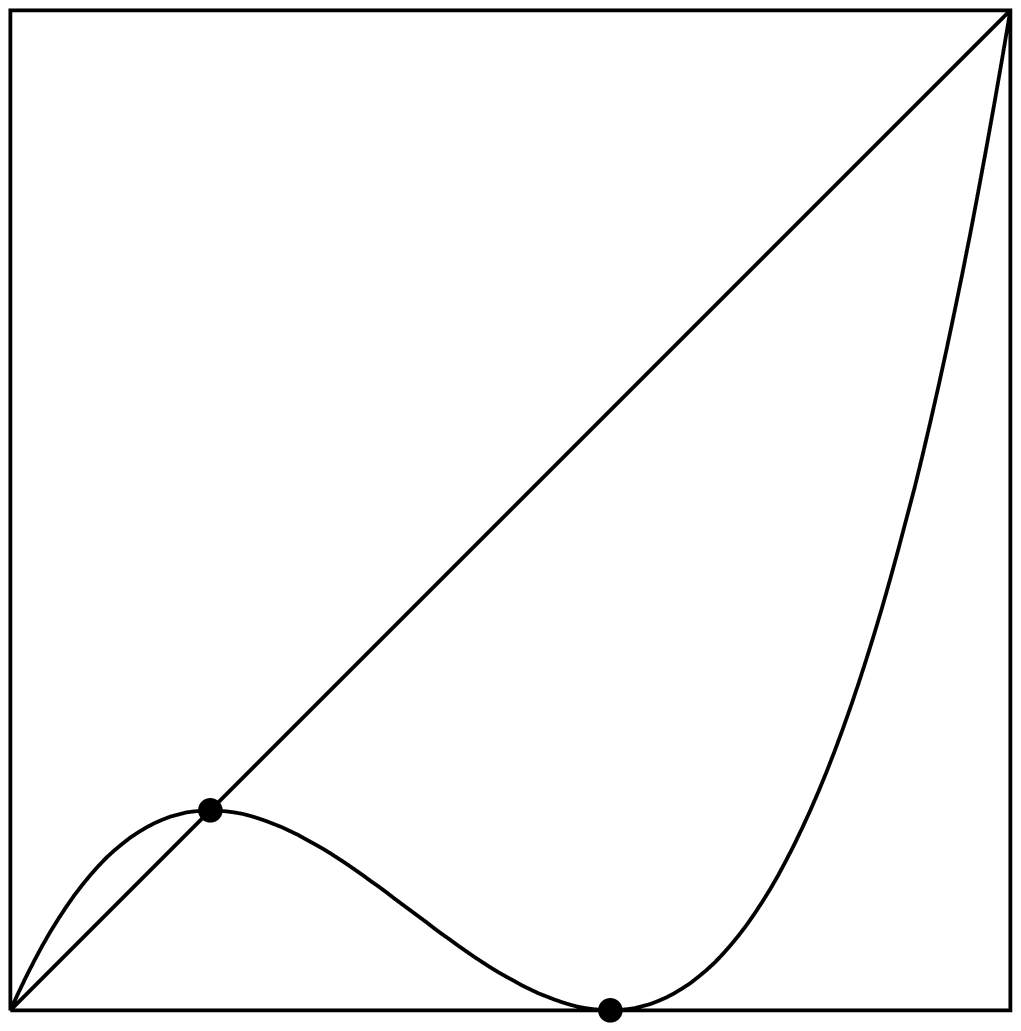,height=1.2in}}
{\wideQP\bit Figure 16. The periodic kneading data $~\K_1=\overline
C_1\,,\;
\K_2=\overline I_0~$ is tight, yet is
represented by an entire one-parameter family of
cubic maps. However, only one of these is postcritically finite.\par}
\endinsert

The proofs follow.\smallskip

{\bf Proof of B.1 and B.2.} First note that the two statements B.1 and B.2
are
equivalent to each other: The itinerary of any non folding
point on the $j$-th plateau is either $(I_{j-1}\,,\,\K_j)$ or
$(I_j\,,\,\K_j)$.
If the orbit of $x$ hits each plateau at most at its center point,
then as in \S5 one can compute $x$ from its itinerary, and it is not
difficult to see that this itinerary
cannot contain any $(I_{j-1}\,,\,\K_j)$ or $(I_j\,,\,\K_j)$.

If the orbit of $x$ hits some plateau off-center,
then clearly we can move $x$ slightly without changing its itinerary.
Similarly, if the $i$-th folding point orbit hits the {\it interior\/} of
some
plateau at a non folding point, then varying the $i$-th folding value
but keeping all other folding values fixed we obtain a one-parameter
family of stunted sawtooth maps with the same kneading data. If,
for one or more values of $i$, the $k(i)$-th forward image of the folding
value $v_i$ lies at an {\it endpoint\/} of the $j(i)$-th plateau, then the
argument is slightly more delicate: Vary all such
folding values at unit speed in such a way that this $k(i)$-th
forward image moves with speed $s^{k(i)}\ge 1$
towards the interior of the corresponding plateau. Of course the $j(i)$-th
folding value may itself be varying at unit speed, but in this
case the endpoints of the
plateau will then move at speed $1/s<1$, so that the image will still
move into the interior of the plateau, and the kneading data will not
change.
Thus $S_\p$ is not tight.

Conversely, if the orbit of $x$ hits every plateau at most at its center
point,
then it follows from the remarks above that $x$ is tight, and the
corresponding statement for kneading data follows easily.\QED

{\bf Proof of B.3.} We will make use of a very special form of {\bit
Hubbard
tree\/}, as described by Douady and Hubbard [DH1], and in a more precise
form by Poirier [Po1, Theorem A]. Since we are interested only in real maps
 with
non-degenerate critical points, our
Hubbard tree $T$ can be identified simply with an interval $[1,n]$ of real
numbers, with the integer points $1\,,\,2\,,\,\ldots
\,,\,n$ as vertices. Here $m$ of these vertices
are to be designated as {\bit critical points\/} (hence $n\ge m$). As
final element of structure, there is to be a mapping $F$ from this finite
set
$\{1\,,\,2\,,\,\ldots\,,\,n\}$ to itself which satisfies three conditions:

(i) $F$ must take consecutive vertices to distinct points.

(ii) Extending linearly over the intermediate intervals $[i\,,\,i+1]$,
we require that two consecutive intervals should map with opposite
orientations (so that their images overlap)
if and only their common vertex is a critical point.

The third condition can be stated in several equivalent forms
(compare [Po1, II.3.11]), but the following will be convenient for our
purposes.

(iii) {\bf Poirier Expansion Condition.} For any interval
$[i\,,\,i+1]\subset T$ there must be an integer $k\ge 0$ so that
the $k$-th forward image of this interval contains a critical point.

{\QP{\bf Assertion:} \it These conditions are satisfied if and only if
there is a polynomial map $f$ of degree $m+1$, necessarily
unique up to positive affine conjugacy, and an order preserving embedding
of the vertices $\{1\,,\,2\,,\,\ldots\,,\,n\}$ of $T$
into the real line which sends the designated
critical points in $T$ to the critical points of $f$, and which
conjugates $F$ (on this set of vertices) to $f$.\smallskip}

Poirier's proof (of a more general theorem which implies this Assertion)
is based on Thurston's
characterization of postcritically finite rational maps
(see [DH2]), and relies on earlier work by Douady, Hubbard, Bielefeld,
Y. Fisher, and others (see [Po2] and references therein).

To apply this result to the proof of B.3, we start with the finite
subset $\Sigma\subset
\A^\N$ consisting of the critical orbits $(C_j\,,\,\K_j)$ together with
all of their terminal segments. Number the elements of $\Sigma$
in increasing order from $1$ to $n$, and let the shift map on $\Sigma$
correspond to the mapping $F$ on $\{1\,,\,\ldots\,,\,n\}$. Then
Condition (i) above follows from
tightness, together with the hypothesis that
$\K_i\ne \K_{i+1}$, while Condition (ii) follows from admissibility,
and Condition (iii) follows trivially
from the definition of kneading sequences
as itineraries. Thus the Assertion above yields a polynomial map
$f:\R\to\R$. Let $I\subset\R$ be the set of points with bounded orbit.
Since $f$ is postcritically finite, all critical orbits are bounded.
It follows easily that $I$ is a closed
interval, and that $f|_I$ is boundary
anchored. There cannot be any smaller closed interval $J\subset I$ which
contains all critical points and is $f$-invariant and boundary anchored.
For then a component of $\overline{I\ssm J}$ would contain
an attracting point (of period $\le 2$), which is impossible
in the postcritically finite case.
\QED

Closely related is an alternative characterization of tightness,
in the eventually periodic case.

{\QP{\bf Corollary B.5.} \it Let $\bold K$ be
admissible eventually periodic
kneading data with $\K_i\ne\K_{i+1}$ for all $i$.
Then $\bold K$ is tight if and
only if the collection of all terminal segments of
the folding point itineraries $(C_i\,,\,\K_i)$ has the
property that consecutive elements (in the ordering of \S2) map to
distinct points under the shift map.\medskip}

{\bf Proof.} This follows immediately from the discussion
above.\QED\smallskip

{\bf Proof of B.4.} Although it is not difficult to prove B.4 directly
from the definition of tightness, the argument would be completely formal.
It is more intuitive to make use of B.5. In other words, we
look at the associated Hubbard tree, whose set of
vertices can be identified with the union of the critical orbits
for a representative map, and check directly that
consecutive vertices map to distinct points. For the union
of the periodic critical orbits, there is nothing to prove since distinct
points clearly map to
distinct points. Any two critical points which map to the right endpoint,
corresponding to $\I_\max$, have the same parity, and hence are separated
by another critical point. We must also check that the rightmost critical
point and the $\I_\max$ point have distinct images,
$~\K_m\ne\roman{shift}(\I_\max)$. But if $\sigma_m=+1$
then $~\roman{shift}(\I_\max)=\I_\max\ne\K_m$, while if
$\sigma_m=-1$ then $\roman{shift}(\I_\max)=\I_\min\ne\K_m$,
using the hypothesis of B.4 in both cases.
The discussion of critical points mapping to $\I_\min$ is completely
analogous.\QED\medskip

\comment
\centerline{$***$ ALTERNATIVE ARGUMENTS FOLLOW. Which (or both) is best ?
$***$}

{\bf Alternative Proof.} We will first show
that the sequences $\I_\min$ and $\I_\max$
are tight. Since every terminal segment of $\I_\min$
or $\I_\max$ is itself equal to either $\I_\min$ or $\I_\max$,
it suffices to show that $\I_\max\ne (I_m\,,\,
\K_m)$ and that $\I_\min\ne(I_0\,,\,\K_1)$. Note more precisely that

$$   \I_\max\=\cases(I_m\,,\,\I_\max)\qquad\text{if}\qquad \sigma_m\=+1~,\\
     (I_m\,,\,\I_\min)\qquad\text{if}\qquad \sigma_m\=-1~.
\endcases $$
On the other hand, by the hypothesis of B.4, the kneading
sequence $\K_m$ can only equal $\I_\max$ when $\sigma_m=-1$, and can
only equal $\I_\min$ when $\sigma_m=+1$. Together with
a completely analogous argument
for $\I_\min$, this proves that both $\I_\max$ and $\I_\min$ are tight.

Suppose now that the critical itinerary $(C_i\,,\,\K_i)$ is periodic but
not tight. Then it must contain some
periodic terminal segment of the form
$(I_{j-\epsilon}\,,\,\K_j)$ with $\epsilon=0$
or $\epsilon=1$. This $\K_j$ cannot be equal to $\I_\min$ or
$\I_\max$ since it contains the symbol $C_i$,
so it follows from the hypothesis of B.4 that $(C_j\,,\,\K_j)$ must
be periodic. But $(I_{j-\epsilon}\,,\,\K_j)$ is also periodic, which is
impossible since two periodic sequence cannot differ at only one place.
\QED
a33333333333333333333
\endcomment
\bigskip
\maybebreakhere
\centerline{\bf Appendix C:
Monotonicity vs Antimonotonicity.}\medskip

As we mentioned in \S 8, one nice monotonicity property in parametrized
families would be to have a line through each point in the parameter
space such that the dynamical complexity increases monotonically
along the line.
This augmentation of complexity could be precisely associated to
the fact that no periodic orbit disappears as we follow the line: in one
dimension this implies (but is not required for) monotonicity of the
topological entropy. For the quadratic family, it is true that no
periodic orbit disappears. (If one does not like to follow deformation of
orbits
from their birth or to their destruction, one can think of the equivalent
statement that the number of orbits of each order type does not
decrease). Now think of the quadratic maps as singular
maps of $\R ^2$, writing
$\tilde Q_{v,0}(x,y)= (4vx(1-x)+y,0)$.
Such
maps are limits of diffeomorphisms of the plane, the H{\'e}non maps,
obtained for each $b\not = 0$ by writing
$\tilde Q_{v,b}(x,y)= (4vx(1-x)+y,by)$. Numerical computations have
suggested for a long time that periodic orbits both appear
and disappear along typical one-parameter families extracted from
the two-parameter H{\'e}non family. The deep reason for these
observations was finally given by the following result (Compare [KKY]).

{\QP{\bf Antimonotonicity Theorem.} \it
In any neighborhood of a nondegenerate, homoclinic-tangency parameter
value of a one-parameter $C^3$-smooth family
of dissipative diffeomorphisms
of the plane, there must be both infinitely many orbit-creation and
infinitely many orbit-annihilation parameter values.\smallskip}

It has been conjectured ([DG], [DGYKK], see also [DGK]) that
a similar statement should work for $m$-modal maps as soon as
$m>1$. More precisely, calling {\bit antimonotone} a parameter value
approached both by infinitely many orbit-creation and
infinitely many orbit-annihilation parameter values,
the following can be extracted from [DGYKK]:

{\QP{\bf Antimonotonicity Conjecture.} \it
A smooth one-dimensional map depending on one parameter has antimonotone
parameter values whenever two critical points have disjoint orbits and
are contained in the interior of a chaotic attractor.\smallskip}

This conjecture would not overtly contradict 9.4. If we follow a (not
necessarily smooth) path within some isentrope, note that any significant
amount of periodic orbit creation must be offset by a roughly equivalent
amount of orbit annihilation, so that topological entropy will remain
constant.
Even if we could find a smooth curve which is transverse to the family of
isentropes, so that topological entropy increases monotonically along it,
it
seems possible that orbit annihilations could be dense, but outweighed
by even more orbit creations. Similarly, the Antimonicity Theorem above
does not contradict the possibility that isentropes for the family of
real H\'enon maps may be connected, although we have no reason to
conjecture
this.

The analogue of the Antimonoticity
Conjecture for the stunted sawtooth families is certainly
false, since by 5.8, it is very easy to find smooth
curves along which there are only
orbit creations. Thus, if the conjecture is true for the cubic family,
then any complexity preserving correspondence
between the stunteed sawtooth and cubic parameter triangles must be very
wild
indeed.

\bigskip\bigskip\bigskip
\maybebreakhere

\centerline{\bf References.}\smallskip

\ref [AKM] R.L. Adler, A.G. Konheim, and M.H. McAndrews, ``Topological
entropy'', {\sl Trans. Amer. Math. Soc. \bf 114} (1965) 309--319.

\ref [ALM] L. Alsed\`a, J. Llibre, and M. Misiurewicz, {\sl Combinatorial
Dynamics and Entropy in Dimension One\/}, (World Scientific ,
Singapore, 1993).

\ref [BCMM] C. Bernhardt, E. Coven, M. Misiurewicz, and I. Mulvey,
``Comparing periodic orbits of maps of the interval'',  {\sl Trans. Amer.
Math. Soc. \bf 333} (1992) 701-707.

\ref [Be] A.F. Beardon, {\sl A Primer on Riemann Surfaces}, London
Mathematical Society Lecture Notes Series {\bf 78} (Cambridge University
Press, Cambridge, 1984).

\ref [BK] L. Block and J. Keesling, ``Computing topological entropy
of maps of the interval with three monotone
pieces'', {\sl J. Stat. Phys. \bf 66} (1992) 755--774.

\ref [BL] A. M. Blokh and M. Lyubich, `` Non-existence of wandering
intervals and structure of topological attractors of one-dimensional
dynamical
systems II, The smooth case'', {\sl Erg. Th. and Dynam. Sys. \bf 9} (1989)
751--758.

\ref [BMT] K. Brucks, M. Misiurewicz and C. Tresser,
``Monotonicity properties of the family of trapezoidal maps'', {\sl
Commun. Math. Phys. \bf 137} (1991) 1--12.

\ref [Bo] R. Bowen, {\sl On Axiom A Diffeomorphisms\/}, Proc. Reg. Conf.
Math. {\bf 35}, (1978).

\ref [Br] M. Brown, ``A proof of the generalized Schoenflies theorem'',
{\sl Bull. Amer. Math. Soc. \bf 66} (1960) 74-76. (See also: M. Brown,
``The monotone union of open $n$-cells is an open $n$-cell'', {\sl Proc.
Amer. Math. Soc. \bf 12} (1961) 812--814.)

\ref [BR] V. Baladi and D. Ruelle, ``An extension of the
theorem of Milnor and Thurston on the zeta functions of interval maps'',
{\sl Ergodic Theory and Dynam. Sys. \bf 14} (1994) 621--632.

\ref [BST] N. J. Balmforth, E. A. Spiegel and C. Tresser, ``The topological
entropy of one-dimen\-sional maps: approximation and bounds'',
{\sl Phys. Rev. Lett. \bf 80} (1994) 80--83.

\ref [D] A. Douady, ``Topological entropy of unimodal maps:
Monotonicity for quadratic polynomials'', pp. 65-87 of
{\sl Real and Complex Dynamical Systems}, (B. Branner and P. Hjorth Eds.)
(Kluwer, Dordrecht, 1995).

\ref [DG] S. P. Dawson and C. Grebogi, ``Cubic
maps as models of two-dimensional antimonotonicity'', {\sl Chaos,
Solitons \& Fractals \bf 1} (1991), 137--144.

\ref [DGK] S. P. Dawson, C. Grebogi and  H. Ko{\c c}ak,
``A geometric mechanism for antimonotonicity in scalar maps with two
critical points'', {\sl Phys. Rev. E \bf 48} (1993) 1676--1682.

\ref [DGKKY] S. P. Dawson, C. Grebogi, I. Kan, H. Ko\c cak and
J. A. Yorke, ``Antimonotonicity: inevitable reversals of period
doubling cascades'', {\sl Phys. Lett.A \bf 162} (1992) 249--254.

\ref [DGMT] S. P. Dawson, R. Galeeva, J. Milnor and C.
Tresser, ``A monotonicity conjecture for real
cubic maps'', pp. 165-183 of
{\sl Real and Complex Dynamical Systems}, (B. Branner and P. Hjorth Eds.)
(Kluwer, Dordrecht, 1995).

\ref [DH1] A. Douady and J. Hubbard, ``A proof of Thurston's topological
characterization of rational maps'',
{\sl Acta Math. \bf 171} (1993) 263--297.

\ref [DH2] A. Douady and J. H. Hubbard,``Etude dynamique des
poly\-n\^omes quadratiques complexes'', {\sl I} (1984) \& {\sl II}
(1985), {\sl Publ. Mat. d'Orsay}.

\ref [E] A. Epstein, ``Algebraic dynamics: contraction, finiteness, and
transversality principles", manuscript in preparation. (See also ``Towers
of Finite Type Complex Analytic Maps'', Thesis, CUNY, 1993.)

\ref [F] P. Fatou, ``Sur les \'equations fonctionnelles, II'', {\sl
Bull. Soc. Math. France \bf 48} (1920) 33--94.

\ref [FK] H.M. Farkas and I. Kra, {\sl Riemann Surfaces}, Graduate Texts
in Mathematics {\bf 71} (Springer Verlag, New York, Second Edition, 1992).
%\ref [Ga] R. Galeeva, ``Kneading sequences for piecewise linear bimodal
%maps'', to appear.

\ref [FT] B. Friedman and C. Tresser, ``Comb structure in
hairy boundaries: some transition problems for circle maps'',
{\sl Phys. Lett. \bf 117A} (1986) 15--22.

\ref [GS] J. Graczyk and G. \'Swi\c atek, ``Generic hyperbolicity in the
logistic family'', {\sl Annals of Math. \bf 146} (1997) 1--52.

\ref [Gu] J. Guckenheimer, {\it in} {\sl Dynamical Systems},
C.I.M.E. Lectures (J. Guckenheimer, J. Moser and S. Newhouse),
 Progress in Mathematics {\bf 8}, (Birkhauser, New York, 1980).

\ref [Ha] G. R. Hall, ``A $C^\infty$ Denjoy counterexample''
{\sl Ergodic Theory and Dynamical Systems \bf 1} (1981) 261--272.

\ref [He] C. Heckman, ``Monotonicity and the Construction of Quasiconformal
Conjugacies in the Real Cubic Family'', Thesis, Stony Brook, 1996.

\ref [HKC] L.P. Hurd, J. Kari, and K. Culik, ``The topological entropy of
cellular automata is uncomputable", {\sl Erg. Th. and Dynam. Sys. \bf 12}
(1992) 255--265.

\ref [K] V. Kaloshin,  ``Generic Diffeomorphisms with Superexponential
Growth of Number of Periodic Orbits'', Stony Brook I.M.S. Preprint 1999\#2.

\ref [Ka] A. Katok, ``Lyapunov exponents, entropy and periodic orbits of
diffeomorphisms'', {\sl Pub. Math. IHES \bf 51} (1980) 137--173.

\ref [KKY] I. Kan, H. Ko\c cak and
J. A. Yorke, ``Antimonotonicity: Concurrent creation and annihilation of
periodic orbits'', {\sl Annals of Math. \bf 136} (1992) 219--252.

% \ref [KN] Y. Komori and K. Nishizawa, ??

\ref [L1] M. Lyubich, ``Non-existence of wandering intervals and structure
of topological attractors of one-dimensional dynamical systems I,
The case of negative Schwarzian derivative'', {\sl Erg. Th. and Dynam. Sys.
 \bf 9} (1989) 737--750.

%\ref [L2] M. Lyubich, ``Geometry of quadratic polynomials: moduli,
%rigidity, and local connectivity'', Stony Brook I.M.S. Preprint 1993\#9.

\ref [L2] M. Lyubich, ``Dynamics of quadratic polynomials, {\sl I} and {\sl
 II}'', {\sl Acta Math. \bf 178} (1997) 185--247, 247--297.

\ref [M1] J. Milnor, ``Remarks on iterated cubic maps'', {\sl
Experimental Math. \bf 1} (1992) 5--24.

\ref [M2] J. Milnor, ``Hyperbolic components in Spaces of Polynomial
Maps (with an appendix by A.~Poirier)'', Stony Brook I.M.S.
Preprint 1992\#3.

\ref [M3] J. Milnor, ``On
cubic polynomials with periodic critical point'', in preparation.

\ref [MaT1] R. S.~MacKay and C.~Tresser,``Boundary of topological chaos
for bimodal maps of the  interval'', {\sl J. London
Math. Soc. \bf 37} (1988), 164--181.

\ref [MaT2] R. S.~MacKay and C.~Tresser,
``Some flesh on the skeleton: the bifurcation
structure of bimodal maps'', {\sl Physica \bf 27D} (1987) 412--422.

\ref [Mc1] C. McMullen, ``Automorphisms of rational maps'', pp. 31--60 of
{\sl Holomorphic Functions and Moduli I\/}, ed. Drasin, Earle,
Gehring, Kra \& Marden; MSRI Publ.$\,10$, (Springer, New York, 1988).

\ref [Mc2] C. McMullen, {\sl Complex Dynamics and Renormalization}, Ann.
Math. Studies {\bf 135}, (Princeton University Press, Princeton, 1994).

\ref [Mis1] M. Misiurewicz, ``On non-continuity of topological entropy'',
{\sl Bull. Ac Pol. Sci., Ser. Sci. Math. Astr. Phys. \bf 19} (1971)
319--320.

\ref [Mis2] M. Misiurewicz, ``Horseshoes for mappings of an interval'',
{\sl Bull. Ac Pol. Sci., Ser. Sci. Math. Astr. Phys. \bf 27} (1979)
167--169.

\ref [Mis3] M. Misiurewicz, ``Jumps of entropy in one dimension'',
{\sl Fund. Math. \bf 132} (1989) 215--226.

\ref [Mis4] M. Misiurewicz, ``Continuity of entropy revisited'', {\sl
Dynamical
systems and applications,
World Sci. Ser. Appl. Anal. \bf 4} (1995) 495--503.

\ref [MMS] M. Martens, W. de Mello, and S. van Strien,
``Julia-Fatou-Sullivan theory for real one-dimensional dynamics'',
Acta Math. {\bf 168} (1992) 273--318.

\ref [MN] M. Misiurewicz and Z. Nitecki, `{\sl Combinatorial
Patterns for Maps of the Interval}, Memoirs of the A.M.S. {\bf 456} (1991).

\ref [MSS] N. Metropolis, M. L. Stein, and P. R. Stein
``On finite limit sets for transformations on the unit interval'',
{\sl J. Comb. Theory \bf 15} (1973) 25--44.

\ref [MSz] M. Misiurewicz and W. Szlenk, ``Entropy of piecewise monotone
mappings'', {\sl Studia Math. \bf 67} (1980) 45--63 (Short version:
{\sl Ast{\'e}risque \bf 50} (1977) 299--310).

\ref [MTh] J. Milnor and W. Thurston, ``On iterated maps of the
interval'', {\it in}  {\sl Springer Lecture Notes \bf 1342} (1988),
465--563.

\ref [Mu] P. Mumbr{\'u}, {\sl{Estructura Peri{\`o}dica i
Entropia Topol{\`o}gica de les Aplicacions Bimodals}} Thesis,
Universitat Aut{\`o}noma de Barcelona, 1987.

\ref [MvS1] W. de Mello and S. van Strien,
``A structure theorem in one-dimensional dynamics'',
{\sl Annals of Math. \bf 129} (1989) 519--546.

\ref [MvS2] W. de Melo and S. Van Strien, {\sl One
Dimensional Dynamics}, (Springer Verlag, Berlin, 1993).

\ref [My] P. J. Myrberg, ``Iteration der reellen Polynome zweiten
Graden'', {\sl Ann. Acad. Sci. Fennic. \bf 256A} (1958) 1--10,
{\bf 268A} (1959) 1--13, and {\bf 336A} (1963) 1--18.

\ref [N] S. Newhouse, ``Continuity properties of
entropy'', {\sl Annals of Math.
\bf 129} (1989) 215--235 and {\bf 131} (1990) 409--410.

\ref [NN], K. Nishizawa and A. Nojiri, ``Center curves in the moduli space
of the real cubic maps'', {\sl Proc. Japan Acad. Ser. A Math. Sci. \bf 69}
(1993) 179--184. See also: ``Algebraic geometry of center curves in the
moduli space of the cubic maps'',  {\sl Proc. Japan Acad. Ser. A Math. Sci.
 \bf 70} (1994) 99--103.

\ref [Po1] A. Poirier, ``On post critically finite polynomials, Part II,
Hubbard Trees'', Stony Brook I.M.S. Preprint 1993\#7. (Thesis, Stony Brook
1993.)

\ref [Po2] A. Poirier, "Realizing reduced schemata", Appendix to J. Milnor,
 ``Hyperbolic components in Spaces of Polynomial Maps, Stony Brook I.M.S.
Preprint 1992\#3.

\ref [Pr] C. Preston, ``What you need to know to
knead'', {\sl Advances Math. \bf 78} (1989) 192--252.

\ref [Re] M. Rees, ``A minimal positive entropy homeomorphism of the
$2$-torus'', {\sl J. London Math. Soc. (2) \bf 23} (1981) 537--550.

\ref [Ro] J. Rothschild, {\sl On the Computation of
Topological Entropy}, Thesis, CUNY, 1971.

\ref [RS] J. Ringland and M. Schell, ``Genealogy and bifurcation
skeleton for cycles of the iterated two-extremum map of the interval'',
{\sl SIAM J. Math. Anal. \bf 22} (1991) 1354--1371.

\ref [RT] J. Ringland and C. Tresser, ``A genealogy for finite
kneading sequences of bimodal maps of the interval'',
{\sl Trans. Amer. Math. Soc. \bf 347} (1995) 4599--4624.

\ref [SI] A.N. Sharkovskii and A.F. Ivanov, ``$C^\infty$-mappings
of an interval with attracting cycles of arbitrary large periods'',
{\sl Ukrain. Mat. Zh. \bf 35} (1983) 455--458.

\ref [St] J. Stimson, {\sl Degree Two Rational Maps with a Periodic
Critical Point}, Thesis, Univ. Liverpool 1993.

%\ref [Su] D. Sullivan, Private communication.  ???
%\ref [Su] D. Sullivan, ``Bounds, quadratic
%differentials, and renormalization conjectures'', pp. 417--466 of
%{\sl A.M.S. Centennial
%Publications, Vol. {\bf 2}, Mathematics into the Twenty-first Century},
%Am. Math. Soc. 1992. ??

\ref [Sw] G. \'Swi\c atek, ``Hyperbolicity is dense in the real
quadratic family'', Stony Brook I.M.S. Preprint 1992\#10.

\ref [Ts] M. Tsujii, ``A simple proof for monotonicity of entropy in the
quadratic family'',  Preprint Hokkaido University, 1998.

\ref [Y] Y. Yomdin, ``Volume growth and entropy'', {\sl Isr. J. Math. \bf
57} (1987) 285--300. (See also ``$C^k\!$-resolution of semialgebraic
mappings'', ibid. pp. 301--317.)
\bigskip

\vfill

\halign{\qquad\qquad #\hfil &\qquad \qquad #\hfil\cr
John Milnor & Charles Tresser\cr
Institute for Mathematical Sciences & I.B.M. \cr
SUNY &  P.O. Box 218\cr
Stony Brook NY 11794-3651 & Yorktown Heights NY 10598\cr
\noalign{\ss}
jack\@math.sunysb.edu & tresser\@watson.ibm.com\cr}
\end

%%%%%%%%%%%%%%%%
\in John Milnor\out
\in Institute for Mathematical Sciences\out
\in SUNY\out
\in Stony Brook NY 11794-3651\out
\in jack\@math.sunysb.edu\bigskip

\in Charles Tresser\out
\in I.B.M.\out
\in P.O. Box 218\out
\in Yorktown Heights NY 10598\out
\in tresser\@watson.ibm.com
\vfil
\eject
\end

   [11] 94m:58195 Ohba, Sachio; Nagano, Yoshihiko; Nishizawa, Kiyoko Chaos
and Newton's
method. (Japanese) Rep. Fac. Engrg. Kanagawa Univ. No. 31 (1993), 1--7.
(Reviewer: Ladislav
Andrey) 58F23 (65L20)

[1] 97d:58121 Alexander, James C.; Hunt, Brian R.; Kan, Ittai; Yorke, James
 A. Intermingled
basins for the triangle map. Ergodic Theory Dynam. Systems 16 (1996), no.
4, 651--662. (Reviewer:
Alec Norton) 58F12

   [6] 94k:58082 Kan, Ittai Open sets of diffeomorphisms having two
attractors, each with an
everywhere dense basin. Bull. Amer. Math. Soc. (N.S.) 31 (1994), no. 1,
68--74. (Reviewer:
Shigenori Matsumoto) 58F12 (58F30)
[11] 93k:58140 Alexander, J. C.; Yorke, James A.; You, Zhiping; Kan, I.
Riddled basins. Internat.
J. Bifur. Chaos Appl. Sci. Engrg. 2 (1992), no. 4, 795--813. (Reviewer: Ian
 Melbourne) 58F12
(34D08 34D45)

In the
case of the cubic family we make an additional modification as follows.
It follows from Lemma 7.3 that each strictly bimodal
bone $B^{\roman{cub}}_\pm(\O)$
consists of a simple arc, possibly together with some disjoint
simple closed curves. We will use the notation $A_\pm(\O)$
for this simple arc, called the {\bit essential bone}, and use the term
{\bit bone-loops} for the closed curves (which we conjecture to not
exist). By the $n\!$-{\bit essential skeleton\/}  $E^{\roman{cub}}_n$
we will mean the union of all of these simple arcs $A_\pm(\O)$,
with period $2\le p\le n$ when $\boldsig=(+-+)$
and with period $3\le p\le n$ when $\boldsig=(-+-)$,
together with $\partial P^{\,\roman{cub}}$.